\documentclass[11pt]{article}
\def\l{\left}\def\r{\right}
\newcommand{\ba}{\noindent $\begin{array}}
\newcommand{\ea}{\end{array}$}
\newcommand{\be}{\begin{equation}}
\newcommand{\ee}{\end{equation}}
\newcommand{\bd}{\begin{displaymath}}
\newcommand{\ed}{\end{displaymath}}
\newcommand{\beq}{\begin{eqnarray*}}
\newcommand{\eeq}{\end{eqnarray*}}
\newcommand{\beqn}{\begin{eqnarray}}
\newcommand{\eeqn}{\end{eqnarray}}
\newcommand{\epc}{\hspace{1pc}}

\usepackage{amsfonts,latexsym,lscape}
\newtheorem{theorem}{Theorem}[section]
\newtheorem{prop}{Proposition}[section]
\newtheorem{lemma}{Lemma}[section]
\newtheorem{corollary}{Corollary}[section]

\newfont{\Bb}{msbm10 scaled\magstep1}

\def\sqr#1#2{{\vcenter{\hrule height .#2pt
      \hbox{\vrule width .#2pt height#1pt \kern#1pt\vrule width.#2pt}
                       \hrule height.#2pt}}}
\def\square{\mathchoice\sqr54\sqr54\sqr{2.1}3\sqr{1.5}3}
\def\qed{\hfill$\square$}
\def\S{{\cal S}}
\def\diag{\,{\rm diag}\,}
\def\argmin{\,{\rm argmin}\,}

\title{The Rate of Convergence of
the Augmented  Lagrangian Method for a Nonlinear Semidefinite Nuclear Norm Composite Optimization Problem\footnote{The research of Liwei Zhang  is supported by the National Natural
Science Foundation of China under project grant No.11571059, No.11731013 and No.91330206.}}

\author{Liwei Zhang\footnote{School of Mathematical Sciences, Dalian University of Technology, Dalian 116024, China (e-mail:lwzhang@dlut.edu.cn).},\quad Yule Zhang\footnote{School of Mathematical Sciences, Dalian University of Technology, Dalian 116024, China (e-mail:zyl@dlut.edu.cn).},\quad
and \quad  Jia Wu\footnote{School of Mathematical Sciences, Dalian University of Technology, Dalian 116024, China (e-mail:wujia@dlut.edu.cn).}}

\begin{document}
\maketitle

\begin{abstract}
We propose two basic assumptions, under which  the rate of convergence of the augmented Lagrange method for a class of composite optimization problems is estimated. We analyze the rate of  local convergence of the augmented Lagrangian method for a  nonlinear semidefinite nuclear norm composite  optimization problem by verifying these two basic assumptions.
     Without requiring strict complementarity,
we prove that, under the constraint nondegeneracy condition and the strong  second order sufficient condition,
 the rate of convergence is linear and the ratio constant is  proportional to $1/c$, where $c$ is the penalty parameter  that exceeds a threshold $\overline{c}>0$. The analysis is based on variational analysis about the proximal mapping of the nuclear norm and the projection operator onto the cone of positively semidefinite symmetric matrices.

\vskip 12 true pt \noindent \textbf{Key words}: Composite optimization, nonlinear semidefinite nuclear norm composite optimization, rate of
convergence, the augmented
Lagrangian method,  variational analysis.
\end{abstract}

\section{Introduction}
\setcounter{equation}{0}


Nuclear norm optimization problems  have seen many applications in engineering and science.
They arise  from the convex relaxation of a rank minimization
problem with noisy data in many machine learning and compressed
sensing applications such as dimensionality reduction,
matrix classification, multi-task learning and matrix
completion, as well as in theoretical applications from mathematics (\cite{Fazel2002},\cite{BAStephen2011},
\cite{Srebro2004},\cite{CTao2009},\cite{LLRabinovich1995}).  A  proximal point algorithmic framework was developed in  \cite{Liu2012} for solving convex nuclear norm optimization problems and  numerical
results  show that the proposed  proximal point algorithms perform quite well in comparison to
several recently proposed state-of-the-art algorithms.
For non-convex nonlinear programming and non-convex semidefinite programming, related to  proximal point algorithms, the augmented Lagrange method is regarded as an effective numerical method. It is quite natural to consider the augmented Lagrange method for  the non-convex  nuclear norm composite optimization problem and study its theoretical properties.
 In the general setting, the augmented
   Lagrangian method can be used to
solve the following composite optimization problem

\vspace{\baselineskip} \noindent (COP)
\vspace{-\baselineskip}
\[
\min \ f(x) + \theta (F(x)) \quad  {\rm s.t.} \quad  h(x)=0,\,\, g(x) \in K\, ,
\]
where $f:\Re^n\mapsto \Re$,$F: \Re^n \mapsto {\cal Z}$, $h:\Re^n \mapsto \Re^m$ and
$g:\Re^n \mapsto {\cal Y}$ are twice continuously differentiable mappings,  $\theta: {\cal Z} \rightarrow \Re \cup \{+\infty\}$ is a proper lower semicontinuous convex function,
 ${\cal Z}$ and ${\cal Y}$ are  finite-dimensional real  Hilbert spaces
equipped with  scalar product $\langle \cdot, \cdot \rangle$ and
 induced norm  $\|\cdot\|$, and $K$ is a closed
convex cone in ${\cal Y}$.

 Let $c>0$
be a parameter. The augmented Lagrangian function with the penalty
parameter $c$ for  problem (COP) is defined as (with no composite term, see \cite[Section
11.K]{RW98})
\begin{equation}\label{al}
\begin{array}{ll}
L_c(x,Y, \mu,\lambda):= & f(x)+ \theta_c (F(x)+Y/c)-\displaystyle\frac{\|Y\|^2}{2c}\\[8pt]
 & +\langle \mu, h(x) \rangle+\displaystyle \frac{c}{2}\|h(x)\|^2+ \frac{1}{2c}\left[ \|
\Pi_{K^*}(\lambda-c g(x))\|^2-\|\lambda\|^2\right]\, ,
\end{array}
\end{equation}
where $(x,Y,\mu, \lambda) \in \Re^n \times {\cal Z} \times \Re^m\times {\cal Y}$ and
$\Pi_{K^*}(\cdot)$ denotes the metric projection operator onto the
set $K^*$($K^*$ is the dual cone of $K$), $\theta_c=e_{1/c}\theta$ and  $[e_{\tau}\theta](\cdot)$ is the Moreau-Yosida regularization of $\theta$ defined by
\begin{equation}\label{eq:morea}
[e_{\tau}\theta](Z)=\inf_{Z' \in \mathbb {\cal Z}} \left\{\theta(Z')+\displaystyle \frac{1}{\tau}\|Z'-Z\|^2\right\}.
\end{equation}

 The augmented Lagrangian method for solving (COP) can be stated
as follows.   Let $c_0>0$  be given. Let
 $(Y^0, \mu^0,\lambda^0) \in {\cal Z} \times \Re^m \times  K^*$ be the initial estimated Lagrange multiplier. At the $k$th iteration,
determine $ x^k$  by minimizing
$L_{c_k}(x,Y^k,\mu^k,\lambda^k)\, $,
compute $(Y^{k+1},\mu^{k+1},\lambda^{k+1})$ by
\[
\left \{ \begin{array}{l}
 Y^{k+1}:={\rm D}\theta_{c_k}(F(x^k)+Y^k/c)^*,\\[4pt]
 \mu^{k+1}:=\mu^k +c h(x^k),\\[4pt]
{\lambda}^{k+1}:=\Pi_{K^*}({\lambda}^k-c_kg(x^k))\, ,
\end{array}
\right.
\]
and update $c_{k+1} $ by
\[
 c_{k+1} :=c_k \quad  {\rm or} \quad c_{k+1}:=\kappa c_k
\]
according to certain rules, where $\kappa >1$ is a given positive number.
 In the case when the sequence of parameters $\{c_k\}$  satisfies
   $c_k\rightarrow +\infty$,  the global convergence of
the augmented Lagrangian method can be  discussed similarly as in \cite{B82}.  In this paper,
instead of considering global convergence properties, we consider
the rate of  convergence of the augmented Lagrangian
method for (COP) when $c_k$ has a finite limit, namely the case in
which $c_k \equiv c$ for all sufficient large $k$.  For
simplicity in our analysis, for $k$ sufficiently large, we
choose $x^k$ as an exact local solution of
$L_{c}(\cdot,Y^k,\mu^k,\lambda^k)$.

The augmented Lagrangian method was proposed by Hestenes
\cite{H69} and Powell \cite{P72} for solving  equality
constrained nonlinear programming problems  and was generalized by
Rockafellar  \cite{R73a} to
  nonlinear programming problems with both equality and inequality constraints.
  For  convex  programming,  Rockafellar \cite{R73a}  established a
saddle point theorem in terms of the augmented Lagrangian  and Rockafellar  \cite{R73b}
 proved
the global convergence of the augmented Lagrangian method for any positive penalty parameter.

For  nonlinear programming,
 the study about the rate of convergence of the augmented Lagrangian method is quite complete. For the equality constrained  problem,
Powell offered a proof  in \cite{P72}  showing that if the linear independence constraint qualification
and the second-order
sufficient condition are  satisfied, then
 the augmented Lagrangian method can  converge locally at a linear rate.
Bertsekas \cite[Chapter 3]{B82}  established an important result on the linear rate of convergence
of the augmented Lagrangian method for  nonlinear programming when the strict complementarity
condition is assumed, in which the ratio constant is proportional
to $1/c$. On the other hand, without assuming the strict
complementarity condition,  Conn et al.
\cite{CGToint91}, Contesse-Becker\cite{Contesse-Becker93}, and Ito
and Kunisch \cite{IKunisch90} derived  linear convergence rate
for the augmented Lagrangian method.

For  nonlinear semidefinte programming, without requiring strict complementarity,
Sun et al. \cite{SSZhang2008} proved that, under the constraint nondegeneracy condition and the strong  second order sufficient condition,
 the rate of convergence of the augmented Lagrangian method  is linear and the ratio constant is  proportional to $1/c$, where $c$ is the penalty parameter  that exceeds a threshold $\overline{c}>0$.  Moreover, Sun et al. \cite{SSZhang2008}  used  a direct way to derive the same linear rate of convergence  under the strict complementarity condition.

The main objective of this paper is to study, without assuming the
strict complementarity,  the rate of convergence of the augmented
Lagrangian method for solving the nonlinear semidefinite nuclear norm composite optimization problem

\vspace{\baselineskip} \noindent (SDNOP)
\vspace{-\baselineskip}
\[
\min \ f(x) + \theta(F(x))\quad
 {\rm s.t.} \quad  h(x)=0\, ,
 \ g(x) \in  {\cal
S}^p_+\, ,
\]
where  $\theta (X)=\|X\|_*$ is the nuclear norm function of $X \in {\cal S}^q$ (for simplicity, here we only consider the nuclear norm of a symmetric matrix),
 ${\cal S}^p_+$ is  the cone of all  positive semidefinite
 matrices in ${\cal S}^p$, the linear space of all $p$ by $p$
symmetric matrices in $\Re^{p\times p}$.

 The organization of  this paper is as follows.
In  Section \ref{general-discussions}, we develop a general theory
on  the rate of convergence of the augmented Lagrangian method for
a class of composite optimization problems under two basic
assumptions. In  Section \ref{preliminaries}, we discuss
 variational properties of the projection over the cone of symmetric positively semidefinite matrices and the proximal mapping of the nuclear norm, and the second-order optimality conditions for nonlinear semidefinite nuclear norm composite optimization problem.
Section \ref{NLSDP-case} is devoted to applying the theory
developed in Section \ref{general-discussions} to nonlinear
semidefinite nuclear norm composite optimization problem. Finally, we give our conclusions in
Section \ref{final-section}.

\section{General discussions on the rate of
convergence} \label{general-discussions} \setcounter{equation}{0}

In this section, we always assume that  the cone $K$ presented in the
optimization  problem (COP) is a closed convex  cone and that $\Pi_{K^*}(\cdot)$ is semismooth everywhere, where
$K^*$ is the dual cone of $K$, i.e.,
\[
K^*:=\{v \in {\cal Y} \,|\, \langle v, z \rangle \geq 0, \ \ \forall\,  z
\in K\}.
\]
The cones $\Re^p_+$,  ${\cal S}_+^p$, $\mbox{epi}\|\cdot\|_2$ and $\mbox{epi}\|\cdot\|_*$  satisfy these assumptions, where $\|\cdot\|_2$ and $\|\cdot\|_*$ stand for the spectral norm of a matrix and the nuclear norm of a matrix, respectively.
Moreover we always assume that ${\rm D}\theta_c(\cdot)$ is semismooth everywhere, where $\theta_c(\cdot)=e_{1/c}\theta(\cdot)$ and  $[e_{\tau}\theta](\cdot)$ is the Moreau-Yosida regularization of $\theta$ defined by
(\ref{eq:morea}).

A feasible point $x\in \Re^n$ to (COP)  is called a stationary point if there exists $(Y,\mu,\lambda)\in {\cal Z} \times \Re^m \times {\cal Y}$
such that the following Karush-Kuhn-Tucker (KKT)
 condition is satisfied at $(x, Y,\mu, \lambda)$:
\begin{equation}\label{eq:stationarity}
\nabla _x L(x,Y,\mu, \lambda) =0,\,  Y \in \partial \theta (F(x)), h(x)=0, g(x)\in
K,\, \lambda \in K^*\,  {\rm and} \, \langle
g(x),\lambda \rangle =0,
\end{equation}
where the Lagrangian function $L: \Re^n \times {\cal Z} \times \Re^m \times {\cal Y} \mapsto \Re$ is defined as
 \[
L(x, Y,\mu, \lambda): = f(x) +\langle Y, F(x) \rangle +\langle \mu, h(x) \rangle-\langle
\lambda, g(x)\rangle.
\]
Any point $(x, Y, \mu,\lambda)\in \Re^n \times {\cal Z} \times \Re^m \times {\cal Y}$  satisfying
(\ref{eq:stationarity}) is  named as a KKT point and the
corresponding point $(Y, \mu,\lambda)$ is called a Lagrange multiplier
at $x$. Let ${\cal M}(x)$   be the set of all
Lagrangian multipliers  at $x$.

 Let $c>0$ and  $\overline{ x} $ be a
stationary point of (COP), namely ${\cal M}(\overline{x})\ne \emptyset$. Since $f, F,
h$, and $g$ are assumed to be twice continuously differentiable,
we know from (\ref{al}), \cite{Zarantonello71} and Chapter 2 of \cite{RW98} that
 the augmented Lagrangian  function  $L_c( \cdot)$ is continuously
differentiable and for any $(x, Y,\mu,\lambda)\in \Re^n \times {\cal Z} \times \Re^m
\times {\cal Y}$,
\begin{equation}\label{compFx}
\begin{array}{ll}
 \nabla_x L_c(x,Y,\mu,\lambda)= & \nabla f(x) + {\rm D}F(x)^*{\rm D} \theta_c(F(x)+Y/c)^*\\[4pt]
 & +{\cal J} h(x)^T (\mu+ch(x)) - {\rm D} g(x)^* \Pi_{K^*}(\lambda-c g(x)).
 \end{array}
\end{equation}
Therefore, from (\ref{eq:stationarity}),  we have $\nabla_x L_c(\overline{x},Y,\mu,\lambda)=0$ for any
$(Y,\mu,\lambda) \in {\cal M} (\overline{x})$.

 For any $(x, Y,\mu,\lambda) \in \Re^n \times {\cal Z}  \times \Re^m \times {\cal Y}$, let
\[
\begin{array}{l}
\Phi_c(x,Y,\mu,\lambda):={\rm D}F(x)^*{\rm D} \theta_c(F(x)+Y/c)^*,\\[4pt]
\Psi_c(x,Y,\mu,\lambda):=  {\rm D}g(x)^* \Pi_{K^*}(\lambda-cg(x)).\end{array}
\]
Let $(x, Y,\mu,\lambda) \in \Re^n \times {\cal Z}  \times \Re^m \times {\cal Y}$. Then
  from the semismoothness of ${\rm D}\theta_c(\cdot)$ and $\Pi_{K^*}(\cdot)$  we obtain that
for any $(\Delta x, \Delta Y,\Delta \mu, \Delta \lambda) \in \Re^n \times {\cal Z} \times \Re^m\times {\cal Y}$,
\begin{equation}\label{composite2-app}
\begin{array}{l}
\quad \partial_B\Phi_c(x,Y,\mu,\lambda) (\Delta x, \Delta Y, \Delta \mu,
\Delta \lambda)
\\[2mm]
=  {\rm D}^2  F(x) (\Delta x) {\rm D} \theta_c(F(x)+Y/c)^*+ {\rm D}F(x)^*
\partial_B [{\rm D} \theta_c]^*(F(x)+Y/c)({\rm D}F(x)\Delta x+\Delta Y/c),\\[2mm]
\quad \partial_B\Psi_c(x,Y,\mu,\lambda) (\Delta x, \Delta Y, \Delta \mu,
\Delta \lambda)
\\[2mm]
=  {\rm D}^2  g(x) (\Delta x) \Pi_{K^*}(\lambda-cg(x)) + {\rm D}g(x)^*
\partial_B \Pi_{K^*} (\lambda-cg(x))(\Delta \lambda-c {\rm D}g(x) \Delta x).
\end{array}
\end{equation}
{} From  (\ref{compFx}) and the definition of $\Psi_c(\cdot)$ we
know that
\[
\begin{array}{l}
\partial_B (\nabla_x L_c)(x,Y,\mu,\lambda)= \\[2mm]
(\nabla^2 f(x),0, 0,0)+\left(\displaystyle\sum_{i=1}^m
(\mu_i+ch_i(x))
\nabla^2 h_i(x)+c {\cal J} h(x)^T{\cal J}h(x), 0, {\cal J} h(x)^T,0\right)\\[4mm]
 +\partial_B \Phi_c(x,Y,\mu,\lambda)
\partial_B \Psi_c(x,Y,\mu,\lambda),
\end{array}
$$
which implies that for any $\Delta x \in \Re^n$,
\begin{equation}\label{comequation3.14}
\begin{array}{l}
 \left( \pi_x \partial_B (\nabla_x L_c)
(x,Y,\mu,\lambda) \right) (\Delta x)\\[2mm]
 =  \nabla^2_{xx} L(x,{\rm D} \theta_c(F(x)+Y/c)^*,
\mu+ch(x), \Pi_{K^*}(\lambda-c g(x))) (\Delta x)
\\[2mm]
\quad +{\rm D}F(x)^*
\partial_B [{\rm D} \theta_c]^*(F(x)+Y/c){\rm D}F(x)(\Delta x)\\[2mm]
  \quad + c {\cal J} h(x)^T{\cal J} h(x) (\Delta x)
   +c {\rm D}g(x)^* \partial_B \Pi_{K^*}(\lambda-c g(x)){\rm D} g(x) (\Delta x)\, ,
\end{array}
\end{equation}
where
\[
\begin{array}{l}
\nabla^2_{xx} L(x,{\rm D} \theta_c(F(x)+Y/c)^*,
\mu+ch(x), \Pi_{K^*}(\lambda-c g(x))) (\Delta x)\\[4pt]
=\nabla^2 f(x)(\Delta x)+{\rm D}^2F(x)(\Delta x){\rm D} \theta_c(F(x)+Y/c)^*\\[4pt]
\quad +
{\rm D}^2h(x)(\Delta x)(\mu+ch(x))
-{\rm D}^2g(x)(\Delta x)\Pi_{K^*}(\lambda-c g(x)).
\end{array}
\]
Let $(\overline{Y},\overline{\mu},\overline{\lambda}) \in {\cal M}(\overline{x})$
be a Lagrange multiplier at  $\overline{x}$. For any linear operators $W_1:{\cal Z} \mapsto {\cal Z}$, $W_2: {\cal Y}
\mapsto {\cal Y}$,  let
\begin{equation}\label{comac}
\begin{array}{ll}
{\cal A}_c(\overline{Y},\overline{\mu},\overline{\lambda},W_1,W_2):= & \nabla^2_{xx}
L(\overline{x},\overline{Y},\overline{\mu},\overline{\lambda})
+{\rm D}F(\overline x)^*W_1{\rm D}F(\overline x)\\[6pt]
& + c  {\cal J}h(\overline{x})^T{\cal J} h(\overline{x})
+c {\rm D}g(\overline{x})^*W_2{\rm D}g(\overline{x}).
\end{array}
\end{equation}
Then for any $\Delta x \in \Re^n$,
\begin{equation}\label{compix}
\begin{array}{l}
\left( \pi_x \,\partial_B (\nabla_x
L_c)(\overline{x},\overline{Y},\overline{\mu},\overline{\lambda}\right) (\Delta
x) \\[10pt]
=\left\{{\cal A}_c(\overline{Y},\overline{\mu},\overline{\lambda},W_1,W_2) (\Delta x):
\begin{array}{l}
W_1\in \partial_B [{\rm D} \theta_c]^*(F(\overline x)+\overline Y/c)\\[4pt]
W_2 \in
\partial_B \Pi_{K^*}(\overline{\lambda}-c g(\overline{x}))
\end{array}
\right\}.
\end{array}
\end{equation}

\vskip 7 true pt
 Next, we make two basic  assumptions for the
constrained optimization composite optimization problem (COP). The first one is about the
positive definiteness of ${\cal A}_c(\overline{Y},\overline{\mu},\overline{\lambda},\cdot,\cdot)$.

\vskip 7 true pt \noindent {\bf Assumption B1}.  We assume that
$(\overline{Y},\overline{\mu},\overline{\lambda})$ is the unique Lagrange
multiplier at $\overline{x}$, i.e.,  ${\cal
M}(\overline{x})=\{(\overline{Y},\overline{\mu},\overline{\lambda})\}$ and that
there exist two positive numbers  $c_0$ and $\underline{\eta}$ such
that for any $c \geq c_0$ and any $W_1\in \partial_B [{\rm D} \theta_c]^*(F(\overline x)+\overline Y/c)$,
 $W_2 \in \partial_B
\Pi_{K^*}(\overline{\lambda}-c g(\overline{x}))$,
\[\left\langle  d ,
{\cal A}_{c}(\overline{Y},\overline{\mu},\overline{\lambda},W_1,W_2) d\right\rangle
\geq\underline{\eta}\left\langle d ,d \right\rangle,  \quad \forall \, d \in
\Re^n.
\]

 Assumption B1 is related to the sufficient optimality conditions for  the
constrained composite optimization problem (COP). It  will be shown in
Proposition \ref{pronlsdp}  that, under the constraint
nondegeneracy condition and the strong second order sufficient
condition (they will be clarified in Section \ref{preliminaries}), Assumption B1 is valid  for  (SDNOP).

\vskip 7 true pt
 Let
$\overline{y}:=(\overline{Y},\overline{\mu},\overline{\lambda})$. Then $\nabla_x
L_c(\overline{x}, \overline{y})=0$.  Let $c_0$ and $\underline{\eta}$ be two positive
numbers defined in Assumption B1 and $c\ge c_0$ be a positive
number.
Since by  (\ref{compix})  and  Assumption B1,  every element in
$\pi_x
\partial_B(\nabla_x L_c)(\overline{x}, \overline{y})$ is positive
definite, we know from  the implicit function theorem for semismooth
functions developed in \cite{S01},
that  there exist an open neighborhood ${\cal O}_{\overline y}$ of
$\overline{y}$ and a locally Lipschitz  continuous function
$x_c(\cdot)$ defined on ${\cal O}_{\overline y}$ such that for any
$y\in {\cal O}_{\overline y}$, $\nabla_x L_c( x_c(y), y)=0$.
Furthermore, since ${\rm D} \theta_c(\cdot)$  and $\Pi_{K^*}(\cdot)$ are assumed to be semismooth
everywhere, $x_c(\cdot)$ is semismooth (strongly semismooth if
$\nabla^2f, {\rm D}^2F, {\rm D}^2 g$, and ${\rm D}^2 h$ are locally Lipschitz
continuous, and both ${\rm D} \theta_c(\cdot)$ and $\Pi_{K^*}(\cdot)$ are strongly semismooth everywhere) at
any point in ${\cal O}_{\overline y}$. Moreover, there exist two
positive numbers $\varepsilon >0$ and $\delta_0>0$ (both depending
on $c$) such that for any $x \in
\mathbb{B}_{\varepsilon}(\overline{x})$ and $y\in
\mathbb{B}_{\delta_0} (\overline{y}) :=\{y\in {\cal Z}\times \Re^m \times {\cal Y}\, | \,
\|y-\overline{y}\|< \delta _0\} \subset {\cal O}_{\overline y}$,
every element in $\pi_x
\partial_B (\nabla_x L_c)(x,y)$ is positive definite.
Thus, for any $y\in \mathbb{B}_{\delta_0} (\overline{y})$,
$x_c(y)$ is the unique minimizer of
 $L_c(\cdot,y)$  over
$\mathbb{B}_{\varepsilon}(\overline{x})$, i.e.,
\begin{equation}\label{def-xc}
\{x_c(y)\} = \argmin \Big\{ L_c(x,y) \,|\, x \in
\mathbb{B}_{\varepsilon}(\overline{x}) \Big \}.
\end{equation}
 Summarizing the above discussions, we obtain the
following proposition.
\begin{prop}\label{prop:general-discussions}
Suppose that Assumption B1 is  satisfied. Let $c\ge c_0$. Then there
exist two positive numbers $\varepsilon >0$ and $\delta_0>0$ $($both
depending on $c$$)$ and a locally Lipschitz continuous function
$x_c(\cdot)$, given by (\ref{def-xc}),  defined on the open ball
$\mathbb{B}_{\delta_0}(\overline{y})$ such that the following
conclusions hold:
\begin{description}

\item[(i)]  The function $x_c(\cdot)$ is semismooth  at any point
in $\mathbb{B}_{\delta_0}(\overline{y})$.

\item[(ii)] If $\nabla^2f, {\rm D}^2F, {\rm D}^2 g$, and ${\rm D}^2 h$ are
locally Lipschitz continuous, ${\rm D} \theta_c(\cdot)$ and $\Pi_{K^*}(\cdot)$  are strongly
semismooth everywhere, then  $x_c(\cdot)$ is strongly semismooth
at any point in $\mathbb{B}_{\delta_0}(\overline{y})$.

\item[(iii)] For any $x \in
\mathbb{B}_{\varepsilon}(\overline{x})$ and $y\in
\mathbb{B}_{\delta_0} (\overline{y})$, every element in $\pi_x
\partial_B (\nabla_x L_c)(x,y)$ is positive definite.

\item [(iv)] For any   $y\in \mathbb{B}_{\delta_0}
(\overline{y})$, $x_c(y)$ is the unique optimal solution to
 \[ \min \ L_c(x, y)
\quad {\rm s.t.} \ x\in {\mathbb B}_{\varepsilon}(\overline{x})\, .
\]
\end{description}
\end{prop}

Let $\vartheta_c:{\cal Z} \times \Re^m \times {\cal Y} \mapsto \Re$ be defined as
\begin{equation}\label{dual}
  \vartheta_c(Y,\mu,\lambda):=\min_{x \in
\mathbb{B}_{\varepsilon}(\overline{x})} L_c(x,Y,\mu,\lambda), \quad
(Y,\mu,\lambda)\in {\cal Z} \times \Re^m \times {\cal Y}\,  .
\end{equation}
Since for each fixed $x\in X$,  $L_c(x, \cdot)$ is a concave
function,  we have that $\vartheta_c(\cdot)$ is also a concave function. By using
the fact that  for any $y\in \mathbb{B}_{\delta_0}
(\overline{y})$, $x_c(y)$ is the unique minimizer of
$L_c(\cdot,y)$ over $\mathbb{B}_{\varepsilon}(\overline{x})$, we
have
\[
\vartheta_c(y)=L_c(x_c(y),y) \, , \quad  y\in
\mathbb{B}_{\delta_0} (\overline{y})\, .
\]
For any $y\in \mathbb{B}_{\delta_0} (\overline{y})$ with
$y=(Y,\mu,\lambda)\in {\cal Z} \times \Re^m \times {\cal Y}$, let
\begin{equation}\label{comnota1}
\left(\begin{array}{c}
Y_c(y)\\
\mu_c(y)
\\
\lambda_c(y)
\end{array} \right):=
\left(\begin{array}{c}
{\rm D} \theta_c(F(x_c(y))+Y/c)^*\\
 \mu+ch(x_c(y)) \\
\Pi_{K^*}(\lambda-c g(x_c(y)))\end{array}\right) \, .
\end{equation}
Then we have
\begin{equation}\label{comnabla0}
\nabla_x L(x_c(y),Y_c(y),\mu_c(y),\lambda_c(y)) =\nabla_x
L_c(x_c(y),y)=0\, , \quad y\in \mathbb{B}_{\delta_0}
(\overline{y})\, .
\end{equation}

\begin{prop}\label{comdcproperty1}
Suppose that Assumption B1 is  satisfied. Let $c\ge c_0$. Then the
concave  function $\vartheta_c(\cdot)$ defined by $(\ref{dual})$
is continuously differentiable on $\mathbb{B}_{\delta_0}
(\overline{y})$ with
\begin{equation}\label{comdcjacobi}
{\rm D}\vartheta_c (y)^*=\left (
\begin{array}{c}
-c^{-1}Y+c^{-1}{\rm D}\theta_c(F(x_c(y))+Y/c))^*\\
h(x_c(y))\\
-c^{-1}\lambda+c^{-1}\Pi_{K^*}(\lambda-c g(x_c(y)))
\end{array}
\right )\, , \quad y =(Y,\mu,\lambda) \in \mathbb{B}_{\delta_0}
(\overline{y})\, .
\end{equation}
Moreover, ${\rm D}\vartheta_c (\cdot)$ is semismooth at any point
in $\mathbb{B}_{\delta_0} (\overline{y})$. It is
  strongly semismooth at any point in
$\mathbb{B}_{\delta_0}(\overline{y})$ if  $\nabla^2f, {\rm D}^2F, {\rm D}^2 g$, and ${\rm D}^2 h$ are locally Lipschitz continuous, and ${\rm D} \theta_c(\cdot)$ and
$\Pi_{K^*}(\cdot)$  are strongly semismooth everywhere.
\end{prop}

\noindent {\bf Proof.}   Let $y =(Y,\mu,\lambda) \in
\mathbb{B}_{\delta_0} (\overline{y})$. Then from (\ref{comnabla0})
and \cite[Theorem 2.6.6]{C83} we have for any $(\Delta Y,\Delta\mu,
\Delta\lambda)\in {\cal Z} \times \Re^m \times {\cal Y}$ that
\[
\begin{array}{ll}
\partial \vartheta_c(y)(\Delta Y,\Delta\mu,
\Delta\lambda)&=
{\cal J}_x L_c(x_c(y),y)(\partial
x_c(y)(\Delta Y,\Delta\mu,
\Delta\lambda))\\
&\\
& \quad  +{\rm D}_{Y}L_c(x_c(y),y)(\Delta Y)+
{\cal
J}_{\mu}L_c(x_c(y),y)(\Delta\mu)+{\rm D}_{\lambda}L_c(x_c(y),y)(\Delta \lambda)\\
&\\
&=
\langle -c^{-1}Y,\Delta Y \rangle +c^{-1}{\rm D}\theta_c(F(x_c(y))+Y/c)(\Delta Y)\\
&\\
& \quad  +\langle h(x_c(y)), \Delta\mu \rangle -c^{-1} \langle \lambda,
\triangle\lambda \rangle+\langle c^{-1} \Pi_{K^*}(\lambda-c
g(x_c(y))),\Delta\lambda\rangle\, .
\end{array}
\]
Thus,  $\partial \vartheta_c (y)(\Delta Y,\Delta\mu,
\Delta\lambda)$ is a
singleton for each  $(\Delta Y,\Delta\mu,
\Delta\lambda)\in {\cal Z} \times \Re^m \times {\cal Y}$.
This implies that $\partial \vartheta_c(y)$ is a singleton.
Therefore, $ \vartheta_c(\cdot)$ is Fr\'{e}chet-differentiable at
$y$ and ${\rm D} \vartheta_c(y)$ is
 given by  (\ref{comdcjacobi}). The continuity of ${\rm D}
 \vartheta_c(\cdot)$ follows from the continuity of
 $x_c(\cdot)$.

The properties on the (strong) semismoothness of ${\rm D}
\vartheta_c (\cdot)$ at $y$  follows directly from (\ref
{comdcjacobi}) and Proposition \ref{prop:general-discussions}.
\qed

\vskip 7 true pt
For any $c\ge c_0$ and $\Delta y :=(\Delta Y,\Delta\mu,
\Delta\lambda)\in {\cal Z}\times \Re^m \times {\cal Y}$,  define
\begin{equation}\label{comcalV}
\begin{array}{l}
\overline{{\cal V}}_c(\Delta y):=\\[8pt]
\left\{\left
 [
 \begin{array}{c}
 c^{-1}W_1{\rm D}F(\overline x)\\
{\cal J} h (\overline{x})\\
 -W_2{\rm D}g(\overline{x})
\end{array}
\right
 ]\right.{\cal A}_c(\overline{y},W_1,W_2)^{-1}\left[-c^{-1}{\rm D}F(\overline x)^*W_1(\Delta Y)\right.\\[5mm]
 \quad \quad \quad \quad \quad \quad \quad \quad \quad \quad\quad \left.-{\cal J}
h(\overline{x})^T\Delta\mu+{\rm D}g(\overline{x})^*W_2(\Delta\lambda)\right]\\[6mm]
 +\left.\left
 (
 \begin{array}{c}
 -c^{-1}\Delta Y+c^{-2}W_1(\Delta Y)\\
 0\\
-c^{-1}\Delta \lambda+ c^{-1}W_2( \Delta \lambda)\end{array} \right
 ) \, \Big |\,  \begin{array}{l}
W_1\in \partial_B [{\rm D} \theta_c]^*(F(\overline x)+\overline Y/c)\\[4pt]
W_2 \in
\partial_B \Pi_{K^*}(\overline{\lambda}-c g(\overline{x}))
\end{array}\right\}.
\end{array}
\end{equation}
Since by Assumption B1,  ${\cal A}_c(\overline{y},W_1,W_2)$ is positive
definite for any $W_1\in \partial_B [{\rm D} \theta_c]^*(F(\overline x)+\overline Y/c)$,$W_2 \in
\partial_B \Pi_{K^*}(\overline{\lambda}-c g(\overline{x}))$, $\overline{{\cal V}}_c(\cdot)$ is well defined. The next
proposition shows that
$\overline{{\cal V}}_c(\cdot)$ is an outer approximation to  $\partial_B [{\rm D}
\vartheta_c]^*(\overline{y})(\cdot)$.

\begin{prop}\label{comthta}
Suppose that Assumption B1 is satisfied. Let $c\ge c_0$. Then for
any $\Delta y :=(\Delta Y,\Delta\mu,
\Delta\lambda)\in {\cal Z}\times \Re^m \times {\cal Y}$,
\begin{equation}\label{combdiff}
\partial_B[{\rm D} \vartheta_c]^*(\overline{y}) (\Delta y )\subseteq
\overline{{\cal V}}_c (\Delta y)\, .
\end{equation}
\end{prop}

\noindent {\bf Proof.}  Choose $\Delta y :=(\Delta Y,\Delta\mu,
\Delta\lambda)\in {\cal Z}\times \Re^m \times {\cal Y}$.  From Proposition
\ref{comdcproperty1}, we know  that ${\rm D}\vartheta_c(\cdot)$ is
semismooth at any  point $y \in
\mathbb{B}_{\delta_0}(\overline{y})$. Let ${\cal
D}_{{\rm D}\vartheta_c}$ denote the set of all
Fr\'{e}chet-differentiable points of ${\rm D}\vartheta_c(\cdot)$ in
$\mathbb{B}_{\delta_0}(\overline{y})$. Then for any $y=(Y,\mu,\lambda) \in {\cal D}_{{\rm D}\vartheta_c}$, we have
\begin{equation}\label{comgrad}
\begin{array}{l}
 {\rm D}^2 \vartheta_c(y)(\Delta y)\\[6mm]
  = \left (
\begin{array}{c}
-c^{-1}\Delta y+c^{-1}[{\rm D}\theta_c]^*(F(x_c(y))+Y/c); {\rm D}F(x_c(y))(x_c)'(y;\Delta y)+\Delta Y/c)\\[6pt]
{\cal J} h(x_c(y)) (x_c)'(y;\Delta y)
\\[6pt]
-c^{-1}\Delta\lambda+c^{-1}\Pi^\prime_{K^*}\Big(\lambda-cg(x_c(y));\Delta\lambda-c
{\rm D} g(x_c(y))(x_c)'(y;\Delta y)\Big)
\end{array}
\right )\, .
\end{array}
\end{equation}

Let $y\in \mathbb{B}_{\delta_0}(\overline{y})$. Now, we derive the
formula for $(x_c)^\prime(y;\Delta y)$.  From (\ref{comnabla0})
and (\ref{comnota1})  we have
\begin{equation}\label{comdiff0}
\begin{array}{lcl}
0&=& \nabla^2_{xx}
L(x_c(y),Y_c(y),\mu_c(y),\lambda_c(y))(x_c)^\prime(y;\Delta y) +c{\cal J}
h(x_c(y))^T{\cal J} h(x_c(y))(x_c)^\prime(y;\Delta y)\\[2mm]
& & +{\cal J}h(x_c(y))^T (\Delta\mu) -{\rm D}g(x_c(y))^*
\Pi^\prime_{K^*}\Big(\xi-cg(x_c(y));\Delta\xi-c{\rm D}
g(x_c(y))(x_c)^\prime(y;\Delta y)\Big)\\[2mm]
&&+{\rm D}F(x_c(y))^*[{\rm D}\theta_c]^{*\prime}(F(x_c(y))+Y/c;{\rm D}F(x_c(y))(x_c)^\prime(y;\Delta y)+\Delta Y/c).
\end{array}
\end{equation}
Since ${\rm D}\theta_c(\cdot)$ and $\Pi_{K^*}(\cdot)$ are semismooth everywhere, there exist $\widehat{W}_1\in
\partial_B [{\rm D}\theta_c]^*(F(x_c(y))+Y/c)$ and $\widehat{W}_2\in
\partial_B \Pi_{K^*}(\lambda-cg(x_c(y)))$ such that
\begin{equation}\label{comderiv}
\begin{array}{l}
[{\rm D}\theta_c]^{*\prime}(F(x_c(y))+Y/c;{\rm D}F(x_c(y))(x_c)^\prime(y;\Delta y)+\Delta Y/c)\\[3mm]
\quad \quad \quad \quad \quad \quad=\widehat{W}_1({\rm D}F(x_c(y))(x_c)^\prime(y;\Delta y)+\Delta Y/c),\\[3mm]
\Pi^\prime_{K^*}\Big(\lambda-cg(x_c(y));\Delta\lambda-c{\cal J} g(x_c(y))
(x_c)^\prime(y;\Delta y)\Big)\\[3mm]
\quad \quad \quad \quad \quad \quad=\widehat{W}_2 (\Delta\lambda-c {\rm D}
g(x_c(y)) (x_c)^\prime(y;\Delta y)).
\end{array}
\end{equation}
 For any $W_1\in
\partial_B [{\rm D}\theta_c]^*(F(x_c(y))+Y/c)$ and $W_2 \in
\partial_B \Pi_{K^*}(\lambda-cg(x_c(y)))$, let
\[
\begin{array}{lcl}
{\cal A}_c(y,W_1,W_2):&= & \nabla^2_{xx} L(x_c(y),Y_c(y),\mu_c(y),\lambda_c(y))
+ {\rm D}F(x_c(y))^*W_1{\rm D}F(x_c(y))\\[2mm]
 & & +c {\cal J}h(x_c(y))^T{\cal J}
h(x_c(y))+c {\rm D}g(x_c(y))^*W_2{\rm D} g(x_c(y))\, .
\end{array}
\]
 From (\ref{comequation3.14}) and the definition of
${\delta_0}$, ${\cal A}_c(y,W_1,W_2)$ is positive definite for any $W_1\in
\partial_B [{\rm D}\theta_c]^*(F(x_c(y))+Y/c)$ and $W_2 \in
\partial_B \Pi_{K^*}(\lambda-cg(x_c(y)))$.
 Then from (\ref{comdiff0}) and (\ref{comderiv}) we obtain that
\begin{equation}\label{comxderi}
\begin{array}{ll}
(x_c)'(y;\Delta y) & ={\cal A}_c(y,\widehat{W}_1,\widehat{W}_2 )^{-1}\left(-{\rm D}F(x_c(y))^*\widehat W_1(\Delta Y/c)\right.\\[6mm]
& \quad \quad \quad \quad \quad  \quad \quad \quad \quad \quad  \left. -{\cal J}
h(x_c(y))^T (\Delta\mu)+{\rm D}g(x_c(y))^*
\widehat{W}_2(\Delta\lambda)\right).
\end{array}
\end{equation}
Therefore, we have from (\ref{comxderi}) and (\ref{comgrad}) that
for any  $y=(Y,\mu,\lambda) \in {\cal D}_{{\rm D}\vartheta_c}$,

\[
\begin{array}{l}
{\rm D}^2 \vartheta_c(y) (\Delta y) \in
\left\{\left
 [
 \begin{array}{c}
 c^{-1}W_1{\rm D}F(x_c(y))\\
{\cal J} h (x_c(y))\\
 -W_2{\rm D}g(x_c(y))
\end{array}
\right
 ]\right.{\cal A}_c(y,W_1,W_2)^{-1}\left[-c^{-1}{\rm D}F(x_c(y))^*W_1(\Delta Y)\right.\\[4mm]
 \quad \quad \quad \quad \quad\quad \quad \quad \quad \quad \quad \quad \quad \quad \quad \quad \quad \quad \quad \quad \left.
 -{\cal J}
h(x_c(y))^T\Delta\mu+{\rm D}g(x_c(y))^*W_2(\Delta\lambda)\right]\\[6mm]
 \quad \quad \quad \quad \quad \quad \quad +\left.\left
 (
 \begin{array}{c}
 -c^{-1}\Delta Y+c^{-2}W_1\Delta Y\\
 0\\
-c^{-1}\Delta \lambda+ c^{-1}W_2( \Delta \lambda)\end{array} \right
 ) \, \Big |\,  \begin{array}{l}
W_1\in \partial_B [{\rm D} \theta_c]^*(F(x_c(y))+ Y/c)\\[4pt]
W_2 \in
\partial_B \Pi_{K}(\lambda-c g(x_c(y)))
\end{array}\right\}.
\end{array}
\]
 which, together with the continuity of $x_c(\cdot)$ and the upper semicontinuity of
$\partial_B \Pi_{K^*}(\cdot)$, implies that  $V(\Delta y) \in
\overline{{\cal V}}_c (\Delta y)$ for any $V \in
\partial_B [{\rm D} \vartheta_c]^*(\overline{y})$.
 Consequently, (\ref{combdiff}) holds.
\qed

 The second  basic assumption required in this
section is stated as below.

\noindent {\bf Assumption B2}. There exist positive  numbers
$\overline{c}\ge c_0$,
 $\mu_0>0$, $\varrho_0>0$,  and $\gamma>1$ such that
 for any $c \geq {\overline{c}}$ and $\Delta y\in {\cal Z} \times \Re^m \times {\cal Y}$,
\begin{equation}\label{eq:direction-added}
\|(x_c)^\prime(\overline{y};\Delta y)\|\le \varrho_0 \|\Delta
y\|/c
\end{equation}
and
 \begin{equation}\label{comimport110}
\left \langle V(\Delta y)+c^{-1} \Delta y, \, \Delta y \right
\rangle \in \mu_0  \left[-1, 1\right]\| \Delta y\|^2/c^{\gamma}\,
\quad \forall\, V( \Delta y ) \in \overline{{\cal V}}_c(\Delta
y)\, .
\end{equation}

 It will be shown in Proposition \ref{thdcnlsdp} that
 Assumption B2 is
valid for (SDNOP)  when the constraint nondegeneracy condition and
the strong second order sufficient condition are satisfied.

Let $C$ be a closed convex set in ${\cal Y}$. It  follows from
 \cite{Zarantonello71} that the metric projector $\Pi_C(\cdot)$
is Lipschitz continuous with the Lipschitz modulus  $1$. Then  for
any $y\in {\cal Y}$, $\partial \Pi_C(y)$ is well defined and it has the following variational properties.
\begin{lemma} \label{Lem:projector-pros} {\rm \cite[Proposition
1]{MSZhao05}}
Let $C\subseteq {\cal Y}$ be a closed convex set. Then, for any $y\in {\cal Y}$
and $V \in \partial \Pi_C(y)$, it holds that
\begin{description}
\item[{(i)}] $V$ is  self-adjoint.

\item[{(ii)}] $\langle d, Vd \rangle \geq 0, \quad  \forall \, d
\in {\cal Y}$.

\item[{(iii)}] $\langle Vd, d - Vd \rangle \geq 0, \quad  \forall
\, d \in {\cal Y}.$
\end{description}
\end{lemma}

Under Assumptions B1 and  B2, we are ready to give the main result
on the rate of convergence of  the augmented Lagrangian method for
the composite optimization problem (COP).

\begin{theorem}\label{comthconv1} Suppose that  $K$ is an nonempty closed convex cone
and that ${\rm D} \theta_c(\cdot)$ and $\Pi_{K^*}(\cdot)$ are semismooth everywhere.  Let Assumptions
B1 and B2 be satisfied. Let $c_0,$ $\underline{\eta},$
$\overline{c}$,  $ \mu_0$, $\varrho_0$, and $\tau$ be the positive
numbers defined in these assumptions. Define
\[
\varrho_1:= 2\varrho_0 \quad {\rm and} \quad \varrho_2:=4\mu_0.
\]
Then for any  $c \geq {\overline{c}}$,  there exist two positive
numbers $\varepsilon$ and  $\delta $ $($both depending on $c$$)$
such that for any $(Y,\mu, \lambda)
\in\mathbb{B}_{\delta}(\overline Y, \overline{\mu}, \overline{\lambda})$, the
problem
 \begin{equation}\label{Problem-op}
   \min \ L_c(x, Y,\mu, \lambda)
\quad {\rm s.t.} \ x\in {\mathbb{B}}_{\varepsilon}(\overline{x})\,
\end{equation}
has a unique solution denoted $x_c(Y,\mu, \lambda)$. The function
$x_c(\cdot,\cdot, \cdot)$ is locally Lipschitz continuous on
$\mathbb{B}_{\delta}(\overline Y, \overline{\mu}, \overline{\lambda})$ and is
semismooth at any point in $\mathbb{B}_{\delta}(\overline Y, \overline{\mu}, \overline{\lambda})$, and for
any $(Y,\mu, \lambda) \in
\mathbb{B}_{\delta}(\overline Y, \overline{\mu}, \overline{\lambda})$, we have
\begin{equation}\label{comest111}
\|x_c(Y,\mu, \lambda)-\overline{x}\| \leq \varrho_1 \|(Y,\mu,
\lambda)-(\overline Y, \overline{\mu}, \overline{\lambda})\|/{c}
\end{equation}
and
\begin{equation}\label{comest112}
\|(Y_c(Y, \mu, \lambda), \mu_c(Y, \mu, \lambda), \lambda_c(Y, \mu, \lambda))- (\overline Y, \overline{\mu}, \overline{\lambda})\| \leq \varrho_2 \|(Y,\mu,
\lambda)-(\overline Y, \overline \mu, \overline \lambda)\|/{c^{\gamma-1}}\, ,
\end{equation}
where $Y_c(Y, \mu, \lambda)$, $\mu_c(Y, \mu, \lambda)$ and $\lambda_c(Y, \mu, \lambda)$ are
defined by $(\ref{comnota1})$, i.e.,
\[
\begin{array}{lll}
Y_c(Y, \mu, \lambda)&=&{\rm D} \theta_c(F(x_c(y))+Y/c)^*,\\
\mu_c(Y, \mu, \lambda)&=&\mu+ch(x_c(y)),
\\
\lambda_c(Y, \mu, \lambda)&=&\Pi_{K^*}(\lambda-c g(x_c(y))).
\end{array}
\]
\end{theorem}

\noindent {\bf Proof.} Let $c\ge {\overline{c}}$. {}From
Proposition \ref{prop:general-discussions} we have already known
that  there exist two positive numbers $\varepsilon >0$ and
$\delta_0>0$ (both depending on $c$) and a locally Lipschitz
continuous function $x_c(\cdot, \cdot,\cdot)$ defined on
$\mathbb{B}_{\delta_0}(\overline Y, \overline{\mu}, \overline{\lambda})$ such
that the function $x_c(\cdot,\cdot, \cdot)$ is semismooth at any point
in $\mathbb{B}_{\delta_0}(\overline Y, \overline{\mu}, \overline{\lambda})$ and
for any $(Y, \mu,\lambda)\in \mathbb{B}_{\delta_0} (\overline Y, \overline{\mu}, \overline{\lambda})$,
$x_c(Y,\mu, \lambda)$ is the unique  solution to
(\ref{Problem-op}).

Denote $y:=(Y, \mu, \lambda)\in {\cal Z} \times \Re^m\times {\cal Y}$. Since $x_c(\cdot)$ is
locally Lipschitz continuous on
$\mathbb{B}_{\delta_0}(\overline{y})$ and is directionally
differentiable at $\overline{y}$, by \cite{Shapiro90} we know that
$x_c(\cdot)$ is Bouligand-differentiable at $\overline{y}$, i.e.,
$x_c(\cdot)$ is directionally differentiable at $\overline{y}$ and
\[
\lim _{y\to \overline{y}} \frac{\|x_c(y)-x_c(\overline{y})
-(x_c)^\prime(y;y-\overline{y}) \|}{\|y-\overline{y}\|} =0\, .
\]
By Proposition \ref{comdcproperty1},   ${\rm D}
\vartheta_c(\cdot)$ is semismooth at $\overline{y}$, and thus is
also Bouligand-differentiable at $\overline{y}$. Then there exists
$\delta \in (0, \delta_0]$ such that for any $y \in
\mathbb{B}_{\delta}(\overline{y})$,
\begin{equation}\label{comest113-A}
\|x_c(y)-x_c(\overline{y})-(x_c)^\prime(\overline{y};y-\overline{y})\|
\leq\varrho_0 \|y-\overline{y}\|/{c}
\end{equation}
and
\begin{equation}\label{comest113-B}
\|{\rm D} \vartheta_c(y)-{\rm D} \vartheta_c(\overline{y})-({\rm D}
\vartheta_c)^\prime(\overline{y};y-\overline{y})\|\leq
\mu_0\|y-\overline{y}\|/c^{\gamma}\, .
\end{equation}

Let $y:=(Y,\mu, \lambda)\in  \mathbb{B}_{\delta}(\overline{y})$ be an
arbitrary point. {}From (\ref{eq:direction-added}),
(\ref{comest113-A}), and the fact that $x_c(\overline{y})
=\overline{x}$, we have
\[
\|x_c(y)-\overline{x}\| \leq
\|(x_c)^\prime(\overline{y};y-\overline{y})\|+\varrho_0\|y-\overline{y}\|
{c} = \varrho_1 \|y-\overline{y}\|/{c}\, ,
\]
which,  shows that    (\ref{comest111}) holds.

 Since ${\rm D} \vartheta_{c}(\cdot)$ is semismooth
 at $\overline{y}$, there exists an element $V
\in
\partial_B [{\rm D} \vartheta_{c}]^*(\overline{y})$ such that
$({\rm D}
\vartheta_{c})^{*\prime}(\overline{y};y-\overline{y})=V(y-\overline{y})$.
By using the fact that $V$ is self-adjoint (see Lemma
\ref{Lem:projector-pros}), we know from (\ref{comimport110}) in
Assumption B2 and Proposition \ref{comthta} that
\begin{equation}\label{eq:squareroot}
\|V(y-\overline{y})+c^{-1} (y-\overline{y})  \|
\le 3 \mu_0\|y-\overline{y}\|/{c^\gamma} \, .
\end{equation}
Therefore,  we have from (\ref{comest113-B}) and
(\ref{eq:squareroot})
\[
\begin{array}{l}
  \quad \|y+c({\rm D}\vartheta_{c})^*(y)-\overline{y}\|
\\[2mm]
 = c \|  ({\rm D}\vartheta_{c})^*(y)-({\rm D}\vartheta_{c})^*(\overline{y})-({\rm D}\vartheta_{c})^{*\prime}(\overline{y};y-\overline{y}) + ({\rm D}\vartheta_{c})^{*\prime}(\overline{y};y-\overline{y})+{c}^{-1}(y-\overline{y})\|
\\[2mm]
 \le {c} \|({\rm D}\vartheta_{c})^*(y)-({\rm D}\vartheta_{c})^*(\overline{y})-({\rm D}\vartheta_{c})^{*\prime}(\overline{y};y-\overline{y})\|
+c\|V(y-\overline{y})+c^{-1} (y-\overline{y})\|\\[2mm]
  \leq  \mu_0 \|y-\overline{y}\|/c^{\gamma-1} +
3{\mu_0}\|y-\overline{y}\|/c^{\gamma-1} =
{\varrho_2}\|y-\overline{y}\|/c^{\gamma-1}\,  ,
\end{array}
\]
which, together with (\ref{comdcjacobi}) and the definitions of
$Y_c(Y,\mu,\lambda)$, $\mu_c(Y,\mu,\lambda)$ and $\lambda_c(Y,\mu,\lambda)$, proves
(\ref{comest112}).
 The proof is completed.
\qed

\vskip 7 true pt Under Assumptions B1 and B2, Theorem
\ref{comthconv1} shows that if for all $k$ sufficiently large with
$c_k\equiv c$ larger than a threshold and if $(x^k,Y^k,
\mu^k,\lambda^k)$ is sufficiently close to $(\overline{x},\overline Y, \overline \mu,
\overline \lambda)$, then the augmented Lagrangian
method can locally be regarded as the gradient ascent method
applied to the dual problem
\[
    \max \ \vartheta_c(Y,\mu,\lambda)
  \quad
  {\rm s.t.}\
(Y,\mu,\lambda)\in {\cal Z}\times  \Re^m \times {\cal Y}\,
\]
with a constant step-length $c$, i.e.,  for all $k$ sufficiently
large
\[
\left(
\begin{array}{c}
y^{k+1}\\
\mu^{k+1} \\
\lambda^{k+1} \end{array}\right)
 = \left(
 \begin{array}{c}
 Y^k\\
 \mu^{k} \\
\lambda^{k} \end{array}\right) + c{\rm D}
{\vartheta_c}(Y^k, \mu^k,\lambda^k)^*.
\]
In
Section \ref{NLSDP-case}, we shall check, under what kind of conditions, Assumptions B1 and B2
imposed in this section can be satisfied by  the nonlinear
semidefinite nuclear norm composite optimization problem.

\section{Variational analysis  for SDNOP}\label{preliminaries}
 \setcounter{equation}{0}
 For studying the rate of convergence of the augmented Lagrange method for the nonlinear semidefinite nuclear norm composite optimization problem (SDNOP), we have to provide some variational properties of $\Pi_{{\cal S}^p_+}(\cdot)$ and $\|\cdot\|_*$, and the second-order optimality conditions for (SDNOP).
\subsection{Variational properties of  $\Pi_{{\cal S}^p_+}(\cdot)$ and $\|\cdot\|_*$}
Since there exists an  nonlinear semidefinite constraint in Problem (SDNOP), we  need more
properties  about the tangent cone of the cone ${\cal S}^p_+$  and  the B-subdifferential of the metric
projector  $\Pi_{{\cal S}_+^p}(\cdot)$ over ${\cal S}_+^p$. Let ${\cal O}^p$ be the set of all $p\times p$ orthogonal matrices.  For a given matrix $M\in\mathcal{S}^p$, there exists $P\in\mathcal{O}^p$ such that
\begin{equation}\label{spec-11}
M=P\Lambda(M)P^T,
\end{equation}
where $\Lambda(M)={\rm
diag}(\lambda_1(M),\lambda_2(M),\ldots,\lambda_p(M))$ and
$\lambda_1(M)\geq\lambda_2(M)\geq\ldots\geq\lambda_p(M)$ are eigenvalues of $M$.
We denote the set of such $P$ in the eigenvalue decomposition by ${\cal O}(M)$.
 Let $\overline{M}\in {\cal S}^p$ and $\overline{M}_+:=
\Pi_{{\cal S}_+^p}(\overline{M})$. Suppose that $\overline{M}$ has
the following spectral decomposition
\begin{equation} \label{eq:orthogonal-decomposition}
\overline{M}  =  \overline{P} \Lambda \overline{P} ^T,
\end{equation}
where $P\in {\cal O}(\overline M)$  and $\Lambda $ is the diagonal matrix of eigenvalues of
$\overline{Z}$. Then
\[
\overline{M}_+  =  P  \Lambda_+P ^T,
\]
where $\Lambda_+$ is the diagonal matrix whose diagonal entries
are the nonnegative parts of the respective diagonal entries of
$\Lambda$ \cite{Higham88, Tseng98}.
Define three index sets of positive, zero, and negative
eigenvalues of $\overline{M}$, respectively, as
\[
\alpha  :=  \{  i \, | \, \lambda_i  > \, 0  \}, \epc \beta  := \{
i \, | \, \lambda_i \, = \, 0  \}, \epc \gamma  :=  \{  i \, | \,
\lambda_i \, < \, 0 \}.
\]
 Write
\[
\Lambda  =  \left[ \begin{array}{ccc} \Lambda_{\alpha} & 0 & 0 \\
[5pt] 0 & 0 & 0 \\ [5pt] 0 & 0 & \Lambda_{\gamma}
\end{array} \right] \epc \mbox{and} \epc
P  =  [ \begin{array}{ccc}P_{\alpha} &
P_{\beta} & P_{\gamma}
\end{array} ]
\]
with $P_{\alpha} \in \Re^{p \times |\alpha|}$,
$P_{\beta} \in \Re^{p \times |\beta|}$, and
$P_{\gamma} \in \Re^{p \times |\gamma| }$.    Let
$\Theta $ be any matrix in ${\cal S}^p$ with entries
\begin{equation}\label{eq:Theta-form}
\left\{ \begin{array}{ll}
 \Theta_{ij}  =\displaystyle
{\frac{\max\{\lambda_i,0\} + \max\{\lambda_j,0\}}{ | \, \lambda_i
\, | +| \, \lambda_j \, |}}& {\rm if}\ (i, j) \notin \beta\times
\beta\, ,
\\[2mm]
 \Theta_{ij}  \in [0,1] & {\rm if}\ (i, j) \in
\beta\times \beta\, .
\end{array} \right.
\end{equation}
 The projection operator $\Pi_{{\cal S}^p_+}(\cdot)$ is
directionally differentiable everywhere in ${\cal S}^p$
\cite{BCS99} and
 is a strongly semismooth
matrix-valued function \cite{SS02}.  For any $ H\in {\cal S}^p$,
we have
\begin{equation} \label{eq:dd-projection}
\Pi_{{\cal S}_+^p} ^{ \prime}(\overline{M};H) =  P
\left[
\begin{array}{ccc}
P_\alpha^T HP_ {\alpha}   &
P_\alpha^T HP_ {\beta} & \Theta_{\alpha
\gamma} \circ P_\alpha^T HP_ {\gamma}
\\ [7pt]
P_\beta^T HP_ {\alpha} &  \Pi_{{\cal
S}_+^{|\beta|}} ( P_\beta^T HP_ {\beta}   )
 & 0 \\ [7pt]
 P_\gamma^T HP_ {\alpha} \circ {\Theta}_{\gamma \alpha} & 0 & 0
\end{array}
\right]  P^T\, ,
\end{equation}
where  $``\circ"$ denotes the Hadamard product \cite{SS02}.
When $\beta =\emptyset$, $\Pi_{{\cal S}^p_+}(\cdot)$ is
Fr\'{e}chet-differentiable at $\overline{M}$ and
(\ref{eq:dd-projection}) reduces to the classical result:
\begin{equation} \label{eq:fd-projection}
{\cal J} \Pi_{{\cal S}_+^p} (\overline{M})H = P \left[
\begin{array}{cc}
P_\alpha^T HP_ {\alpha}  & \Theta_{\alpha
\gamma} \circ P_\alpha^T HP_ {\gamma}
  \\ [7pt]
 P_\gamma^T HP_ {\alpha} \circ {\Theta}_{\gamma \alpha} &  0
\end{array}
\right]  P^T\, \quad \forall \, H\in {\cal S}^p\, .
\end{equation}
 The tangent cone of ${\cal S}_+^p$ at $\overline{M}_+$,  denoted ${\cal
T}_{{\cal S}_+^p} (\overline{M}_+)$,  can be completely
characterized as follows
\[
{\cal T}_{{\cal S}_+^p}  (\overline{M}_+) =\{B\in {\cal S}^p\, |
\, B = \Pi_{{\cal S}_+^p} ^{ \prime}(\overline{M}_+;B)\} =\{B\in
{\cal S}^p\, |\, [P_\beta \ P_\gamma]^T B
[P_\beta \ P_\gamma] \succeq 0\}\, .
\]
The lineality space of ${\cal T}_{{\cal
S}_+^p} (\overline{M}_+)$, i.e., the largest linear space in
${\cal T}_{{\cal S}_+^p} (\overline{M}_+)$, denoted by ${\rm lin}
\left ({\cal T}_{{\cal S}_+^p} (\overline{M}_+)\right)$, takes the
following form:
\[
{\rm lin} \left ( {\cal T}_{{\cal S}_+^p}  (\overline{M}_+)
\right) =\{B\in {\cal S}^p\, |\, [P_\beta \
P_\gamma]^T B [P_\beta \
P_\gamma] =0\}\, .
\]
The  critical cone  of ${\cal S}^p_+$ at $\overline{M} \in {\cal
S}^p$ associated with the problem of finding the metric projection
of $\overline{M}$ onto ${\cal S}^p_+$ (i.e., $\overline{M}_+$) is
defined as \cite[Section 5.3]{BS00}
\[
{\cal C}_{{\cal S}_+^p}(\overline{M}) \, : = \, {\cal T}_{{\cal
S}_+^p}  (\overline{M}_+) \cap \{B \in {\cal S}^p\, |\, \langle B,
\overline{M}_+-\overline{M}\rangle =0\}\, .
\]
Thus, it holds that
\[
{\cal C}_{{\cal S}_+^p}(\overline{M}) \, = \, \left\{  B  \in {\cal
S}^p \, \Big{|} \, P_{\beta}^T B P_{\beta} \succeq
0, \ P_{\beta }^TB P_{\gamma} = 0, \
P_{\gamma }^T B P_{\gamma  } =  0 \right \}.
\]
The affine hull of ${\cal C}_{{\cal S}_+^p}(\overline{M})$, denoted by
${\rm aff}({\cal C}_{{\cal S}_+^p}(\overline{M}))$, can then be
written as
\begin{equation}\label{eq:aff-critical-cone}
{\rm aff}\left({\cal C}_{{\cal S}_+^p}(\overline{M})\right)=\left\{
B  \in {\cal S}^p \, | \, P_{\beta }^T B
P_{\gamma}  = 0, \ P_{\gamma}^T
BP_{\gamma } =  0 \right\}.
\end{equation}

The following lemma on $\partial _B\Pi_{{\cal S}_+^p}
(\overline{M})$ is part of \cite[Proposition 4]{S05}, which is
based on \cite[Lemma 11]{PSS03}.

\begin{lemma}\label{lemma:projector-Jacobian-cone-pre}  Let  $\Theta\in {\cal S}^p$
satisfy $(\ref{eq:Theta-form})$. Then $W \in
\partial _B\Pi_{{\cal S}_+^p} (\overline{M})$ if and only if there
exists $W_0 \in
\partial _B\Pi_{{\cal S}_+^{|\beta|}} (0)$  such that
\[
W(H)=
 P
\left[
\begin{array}{ccc}
P_\alpha^T HP_ {\alpha}   &
P_\alpha^T HP_ {\beta} & \Theta_{\alpha
\gamma} \circ P_\alpha^T HP_ {\gamma}
\\ [7pt]
P_\beta^T HP_ {\alpha} & W_0(
P_\beta^T HP_ {\beta} )
 & 0 \\ [7pt]
 P_\gamma^T HP_ {\alpha} \circ {\Theta}_{\gamma \alpha} & 0 & 0
\end{array}
\right]  P^T
 \quad \forall \, H\in {\cal S}^p\, .
\]
\end{lemma}
From the definition of $\partial _B\Pi_{{\cal
S}_+^{|\beta|}} (0)$ and (\ref {eq:fd-projection}) we know that if
 $W_0\in \partial _B\Pi_{{\cal S}_+^{|\beta|}}
(0)$, then there exist matrices $Q\in {\cal O}^{|\beta|}$ and  $\Omega\in
{\cal S }^{|\beta|}$ with entries $\Omega_{ij}\in [0,1]$  such
that
\[
W_0 (D)= Q(\Omega \circ (Q^T DQ)) Q^T \quad \forall\, D \in {\cal
S}^{|\beta|}\, .
\]
For an extension to the above result, see \cite[Lemma 4.7]{CQT03}.
By using Lemma \ref{lemma:projector-Jacobian-cone-pre} we
obtain the following useful lemma, which does not need   further
explanation.
\begin{lemma}\label{lemma:projector-Jacobian-cone}  For any $W\in \partial
_B\Pi_{{\cal S}_+^p} (\overline{M})$, there exist two matrices
$P\in {\cal O}(\overline{M})$ and  $\Theta \in {\cal S}^p$ satisfying
$(\ref{eq:Theta-form})$ such that
\[
W(H) = P\left( \Theta \circ (P^THP)\right )P^T\quad \forall \, H
\in {\cal S}^p\, .
\]
\end{lemma}



For discussions on the nuclear norm function we  need more
properties  about the first and second-order directional derivatives of $\theta$ and the sub-differential of its proximal mapping.  For a given matrix $X\in\mathcal{S}^q$, there exists $Q\in\mathcal{O}^q$ such that
\begin{equation}\label{X-spec}
X=Q\Lambda(X)Q^T,
\end{equation}
where $\Lambda(X)={\rm
diag}(\lambda_1(X),\lambda_2(X),\ldots,\lambda_q(X))$ and
$\lambda_1(X)\geq\lambda_2(X)\geq\ldots\geq\lambda_q(X)$ are eigenvalues of $X$.
We denote the set of such $Q$ in the eigenvalue decomposition by ${\cal O}(X)$.

Let
$\varpi_1>\varpi_2>\ldots>\varpi_r$ be the distinct eigenvalues of $X$. Define
$$
\begin{array}{l}
a_k:=\{i\,|\,\lambda_i(X)=\varpi_k\},\quad k=1,\ldots,r.
\end{array}
$$
Partition $Q$  as
$Q=[Q_{a_1}\,\,\,Q_{a_2}\,\,\,\ldots\,\,\,Q_{a_r}]$,
where $Q_{a_k}=(q_i:i \in a_k)$ and $Q_{a_k}\in\Re^{q\times |a_k|}$, $k=1,\ldots,r$.

For a given
$H\in\mathcal{S}^q$ and $k\in\{1,\ldots,r\}$, suppose that
$Q^T_{a_k}HQ_{a_k}\in\Re^{|a_k|\times|a_k|}$
has the following spectral decomposition:
$$
(Q^k)^T(Q^T_{a_k}HQ_{a_k})Q^k=\diag(\xi^k_1,\ldots,\xi^k_{|a_k|}),
$$
where $ Q^k\in\mathcal{O}^{|a_k|}(Q^T_{a_k}HQ_{a_k})$
and $\xi^k_i=\lambda_{i}(Q^T_{a_k}HQ_{a_k})$,
$i=1,\ldots,|a_k|$.
Let $\eta^k_1,\ldots,\eta^k_{n_k}$ be the distinct eigenvalues of
$Q^T_{a_k}HQ_{a_k}$ and define
$$
\begin{array}{l}
b^k_j:=\{i\,|\,\xi^k_i=\eta^k_j,\quad i=1,\ldots,|a_k|\}, \quad
j=1,\ldots,n_k.
\end{array}
$$
For simplicity, we denote $\displaystyle
\kappa_i:=\sum^{i}_{j=1}|a_j|$ $(\kappa_0:=0)$,
$\displaystyle \kappa^{(k)}_i=\sum^{i}_{j=1}|b^k_j|$
($\kappa^{(k)}_{0}:=0$), and define the following mappings:
$$
\begin{array}{lll}
&& \nu:\{1,\ldots,n\}\rightarrow\{1,\ldots,r\},\quad
\nu(i)=k,\,\,\, {\rm if} \,\,\,i\in a_k, \\
& & l:\{1,\ldots,n\}\rightarrow\mathbb{N},\quad
l(i)=i-\kappa_{\nu(i)-1}, \\
&& \omega: \{1,\ldots,n\}\rightarrow\mathbb{N},\quad
\omega(i)=j,\,\,\,{\rm if} \,\,\,l(i)\in b^{\nu(i)}_j, \\
&& l':\{1,\ldots,n\}\rightarrow\mathbb{N},\quad
l'(i)=l(i)-\kappa^{\nu(i)}_{\mu(i)-1}.
\end{array}
$$
Then, for $i'\in  b^k_j$, its corresponding index $i \in \{1,\ldots,q\}$ is expressed as $i=\kappa^{(k)}_{j-1}+i'+\kappa_{k-1}$.
Define
\[
\begin{array}{l}
\widehat Q_{a_ka_k}=Q_{a_k}^THQ_{a_k},k=1,\ldots, r,\\[4pt]
\widehat V_k(H,W)=Q_{a_k}^T[W-2H(X-\varpi_k I)^{\dag}H]Q_{a_k},k=1,\ldots, r.
\end{array}
\]
Then we have from \cite[Theorem 3.1]{ZZX2013} that
\[
\begin{array}{l}
\lambda_i'(X;H)=\lambda_{l(i)}(\widehat H_{a_{\nu(i)}a_{\nu(i)}}),\, i=1,\ldots,q,\\[4pt]
\lambda_i''(X;H,W)=\lambda_{l'(i)}\left({Q^{\nu(i)}_{b^{\nu(i)}_{\omega(i)}}}^T\widehat V_{\nu(i)}(H,W)Q^{\nu(i)}_{b^{\nu(i)}_{\omega(i)}}\right),\, i=1,\ldots,q.
\end{array}
\]
Assume that there exists an integer $s_0$ satisfying  $1 \leq s_0 \leq  N_s$ and $\eta^s_{s_0}=0$. Let $b^s_+=b^s_1\cup \cdots\cup b^s_{s_0-1}$, $b^s_0=b^s_{s_0}$ and $b^s_-=b^s_{s_0+1}\cup \cdots\cup b^s_{N_s}$.
 Then we obtain the following proposition about the directional derivative and the second-order directional derivative of $\theta (X)$.
 \begin{lemma}\label{1-2-dd}
 Under the above notations, one has
 he directional derivative of $\theta$ at $X$ along $H$ is expressed as
 \begin{equation}\label{eq:1d}
 \theta'(X;H)=\displaystyle \sum_{i=1}^{s-1}{\rm Tr}(\widehat H_{a_ia_i})-
 \displaystyle \sum_{i=s+1}^r{\rm Tr}(\widehat H_{a_ia_i})+\|\widehat H_{a_sa_s}\|_*.
 \end{equation}
 and the second-order directional derivative of $\theta$ at $X$ along $(H,W)$ is expressed as
\begin{equation}\label{eq:2d}
\begin{array}{ll}
\theta''(X;H,W)
&=\displaystyle \sum_{i=1}^{s-1}{\rm Tr}(\widehat V_i(H,W))-
 \displaystyle \sum_{i=s+1}^r{\rm Tr}(\widehat V_i(H,W))+{\rm Tr}({Q^s_{b^s_+}}^T\widehat V_s(H,W)Q^s_{b^s_+})\\[4pt]
 & \quad -
 {\rm Tr}({Q^s_{b^s_-}}^T\widehat V_s(H,W)Q^s_{b^s_-})
 +
 \|{Q^s_{b^s_0}}^T\widehat V_s(H,W)Q^s_{b^s_0}\|_*.
 \end{array}
\end{equation}
\end{lemma}
{\bf Proof}. For $\theta(X)=\|X\|_*$, the nuclear norm of a symmetric matrix in $X \in {\cal S}^q$, it is the  spectral function corresponding to the symmetric function
$$
\varsigma (z)=\sum_{j=1}^q |z_i|, z=(z_1,\ldots, z_q)^T \in \Re^q,
$$
namely $\theta(X)=\|X\|_*=[\varsigma \circ \lambda](X)$.

Let $\overline z \in \Re^q$. We define
\[
I_+(\overline z)=\{i: \overline z_i >0\},\, I_0(\overline z)=\{i:\overline z=0\},\,
I_-(\overline z)=\{i: \overline z_i <0\}
\]
and
\[
\begin{array}{l}
I_{0+}(\overline z,\Delta z)=\{i \in I_0(\overline z): \Delta z_i >0\},\\[4pt]
 I_{00}(\overline z,\Delta z)=\{i \in I_0(\overline z): \Delta z_i =0\},\\[4pt]
 I_{0-}(\overline z,\Delta z)=\{i \in I_0(\overline z): \Delta z_i<0\}.
 \end{array}
\]
Then the directional derivative of $\varsigma$ at $\overline z$ along $\Delta z$ is
\[
\varsigma'(\overline z;\Delta z)=\displaystyle \sum_{i \in I_+(\overline z)} \Delta z_i-\displaystyle \sum_{i \in I_-(\overline z)} \Delta z_i +\displaystyle \sum_{i \in I_0(\overline z)}| \Delta z_i |
\]
and the second-order parabolic directional derivative at $\overline z$ along $\Delta z$ and $\Delta w$ is
\[
\varsigma''(\overline z;\Delta z,\Delta w)=\displaystyle \sum_{i \in I_+(\overline z)\cup I_{0+}(\overline z,\Delta z)} \Delta w_i-\displaystyle \sum_{i \in I_-(\overline z)\cup I_{0-}(\overline z,\Delta z)} \Delta w_i +\displaystyle \sum_{i \in I_{00}(\overline z,\Delta z)}| \Delta w_i |.
\]
 Then, from the chain rules of directional derivatives  (see Chapter 2 of \cite{BS00}), we obtain
 $$
 \theta'(X;H)=\varsigma'(\lambda (X);\lambda'(X;H))=\displaystyle \sum_{i=1}^{s-1}{\rm Tr}(\widehat H_{a_ia_i})-
 \displaystyle \sum_{i=s+1}^r{\rm Tr}(\widehat H_{a_ia_i})+\|\widehat H_{a_sa_s}\|_*
 $$
 and
$$
\begin{array}{ll}
\theta''(X;H,W) &=\varsigma'(\lambda (X);\lambda'(X;H),\lambda''(X;H,W))\\[6pt]
&=\displaystyle \sum_{i=1}^{s-1}{\rm Tr}(\widehat V_i(H,W))-
 \displaystyle \sum_{i=s+1}^r{\rm Tr}(\widehat V_i(H,W))+{\rm Tr}({Q^s_{b^s_+}}^T\widehat V_s(H,W)Q^s_{b^s_+})\\[4pt]
 & \quad -
 {\rm Tr}({Q^s_{b^s_-}}^T\widehat V_s(H,W)Q^s_{b^s_-})
 +
 \|{Q^s_{b^s_0}}^T\widehat V_s(H,W)Q^s_{b^s_0}\|_*.
 \end{array}
$$
The proof is completed.\hfill $\Box$\\
By direct calculation, we may obtain the following conclusion.
\begin{prop}\label{prop-psi}
 Let $\psi (W)=\theta''(X;H,W)$, then
 \begin{equation}\label{eq:conj}
 \begin{array}{l}
  \psi^*(Y)=\\[10pt]
 \left
 \{
 \begin{array}{ll}
 \begin{array}{l}
 2\displaystyle \sum_{i=1}^{s-1} {\rm Tr}(Q_{a_i}^TH(X-\varpi_iI)^{\dag}HQ_{a_i})\\
 +2\displaystyle \sum_{i=s+1}^{r} {\rm Tr}(Q_{a_i}^TH(X-\varpi_iI)^{\dag}HQ_{a_i})\\[4pt]
 + 2{\rm Tr}({Q^s_{b^s_+}}^TQ_{a_s}^THX^{\dag}HQ_{a_s}Q^s_{b^s_+})\\[4pt]
 +2{\rm Tr}({Q^s_{b^s_-}}^T Q_{a_s}^THX^{\dag}HQ_{a_s}Q^s_{b^s_-})\\[4pt]+2\langle
 {Q^s_{b^s_0}}^T\widehat Y_{a_sa_s}Q^s_{b^s_0}, {Q^s_{b^s_0}}^TQ_{a_s}^THX^{\dag}HQ_{a_s}Q^s_{b^s_0}
 \rangle
 \end{array} & \mbox{if} \left(\begin{array}{l}
 \widehat Y_{a_ia_i}=I_{|a_i|},\\
 \quad \quad \mbox{for }  1\leq i\leq s-1, \\
  \widehat Y_{a_ia_i}=-I_{|a_i|},\\
  \quad \quad \mbox{for } s+1\leq i\leq r,\\
  {Q^s_{b^s_+}}^T \widehat Y_{a_sa_s}Q^s_{b^s_+}=I_{|b^s_+|},\\
  {Q^s_{b^s_-}}^T \widehat Y_{a_sa_s}Q^s_{b^s_-}=I_{|b^s_-|},\\
  \|{Q^s_{b^s_0}}^T\widehat Y_{a_sa_s}Q^s_{b^s_0}\|_2\leq 1,
  \end{array}\right.\\[26pt]
   0 & \mbox{otherwise}.
  \end{array}
  \right.
  \end{array}
   \end{equation}
\end{prop}

Now we characterize elements in  $\partial \theta (X)$ for $X \in {\cal S}^q$. If follows from Page 121 of Borwin and Lewis (2006) \cite{BLewis2006}, for the given $X \in {\cal S}^q$ with the spectral decomposition (\ref{X-spec}),  that $Y \in \partial [\varsigma \circ \lambda]( X)$, or $Y \in \partial \theta (X)$  if and only if there exists $w \in \partial \varsigma (\lambda ( X))$
$$
Y=Q{\rm Diag}\,(w) Q^T,
$$
where $X$ has the spectral decomposition $ X=Q{\rm Diag}(\lambda ( X))Q^T$. Define the following three index sets:
$$
a=\{i: \lambda_i( X)>0\},\,\, b=\{i: \lambda_i( X)=0\},\,\, c=\{i: \lambda_i( X)<0\},
$$
or alternatively $a=a_1 \cup\cdots\cup a_{s-1}$, $b=a_s$ and $c=a_{s+1}\cup \cdots\cup a_r$.
Then, $w \in \partial \varsigma (\lambda (X))$ has the following property
$$
w_a =\textbf{1}_{|a|}, w_{c}=-\textbf{1}_{|c|} \mbox{ and }  -\textbf{1}_{|b|}\leq w_{b}\leq \textbf{1}_{|b|},
$$
and
$$
Y=Q{\rm Diag}\,(w) Q^T=Q_{a}Q_{a}^T+Q_{b}{\rm Diag}\,(w_{b}) Q^T_{b}-Q_{c}Q_{c}^T.
$$
For the index set $b$, we partition it as follows
$b=b_L \cup b_S \cup b_U$:
$$
 b_L=\{i\in b: w_i=-1\}, b_S=\{i\in b: -1< w_i< 1\},b_U=\{i\in b: w_i=1\}.
$$
Then  $Y \in \partial \theta (X)$ can be expressed as
\begin{equation}\label{eq:dec-Y}
Y=Q{\rm Diag}\,(w) Q^T=Q_{a\cup_{b_U}}Q_{a\cup_{b_U}}^T+Q_{b_S}{\rm Diag}\,(w_{b_S}) Q^T_{b_S}-Q_{c\cup_{b_L}}Q_{c\cup_{b_L}}^T
\end{equation}
and for $Z=X+Y$,
\begin{equation}\label{eq:cDecop}
Z=[Q_a\,Q_{b_U}\,Q_{b_S}\,Q_{b_L}\,Q_{c}]\left
[
\begin{array}{ccccc}
\Lambda_a+I_{|a|} &&&&\\[4pt]
& I_{|b_U|} &&&\\[4pt]
&&{\rm Diag}\,(w_{b_S})&&\\[4pt]
&&&-I_{|b_L|}&\\[4pt]
&&&&\Lambda_c-I_{|c|}
\end{array}
\right
]\left
[\begin{array}{c}
Q_a^T\\[4pt]
Q_{b_U}^T\\[4pt]
Q_{b_S}^T\\[4pt]
Q_{b_L}^T\\[4pt]
Q_{c}^T
\end{array}
\right].
\end{equation}
The critical cone of $\theta$ at $Z$ associated with $Y \in  \partial  \theta (X)$ is defined by
\begin{equation}\label{eq:directional-c}
{\cal C}_\theta(Z)=\{H \in {\cal S}^q: \theta'(X;H)= \langle Y, H\rangle\}.
\end{equation}
The next lemma gives an characterization of the critical cone ${\cal C}_\theta$.
\begin{lemma}\label{lem:critic-char}
Let $X,Y,Z \in {\cal S}^q$, $Z=X+Y$ satisfies $Y \in \partial  \theta (X)$. Then $H \in {\cal C}_\theta (Z)$ if and only if
\begin{equation}\label{eq:critical-Prox}
{\cal C}_{\theta}(Z)=\left\{H \in {\cal S}^q:\begin{array}{l}
 Q^{T}_{b_S}H[Q_{b}]=0,
 Q^{T}_{b_U}H[Q_{b_L}]=0\\[4pt]
   Q^{T}_{b_U}HQ_{b_U}\in {\cal S}^{|b_U|}_+,
          Q^{T}_{b_L}HQ_{b_L}\in {\cal S}^{|b_L|}_-
          \end{array}
          \right
          \}.
\end{equation}
\end{lemma}
{\bf Proof.} Noting that
$$
\partial \theta (X)=\{Q_{a}Q_{a}^T+Q_{b}W_{b} Q^T_{b}-Q_{c}Q_{c}^T:W_{b} \in {\cal S}^{|b|}, \|W_b\|_2 \leq 1\},
$$
 where $\|W_b\|_2$ denotes the spectral norm of $W_b$. And the directional derivative of $\theta$ at $X$ along $H$ is
 $$
 \theta'(X;H)=\langle Q_{a}Q_{a}^T-Q_{c}Q_{c}^T, H \rangle+\|Q_{b}^T[H] Q_{b}\|_*.
 $$
Noting that $B$ has the expression
$$
B=Q_{a}Q_{a}^T+Q_{b}{\rm Diag}\, (w_{b}) Q^T_{b}-Q_{c}Q_{c}^T
$$
where $w_b \in \Re^{|b|}$ satisfies $\|w_b\|_{\infty} \leq 1$. Then $\theta'(A;H)= \langle B, H\rangle$ is equivalent to
\begin{equation}\label{eq:hlp10}
\|Q_{b}^T[H] Q_{b}\|_*=\langle Q_{b}^T[H] Q_{b},{\rm Diag}\, (w_{b})  \rangle.
\end{equation}
From Fan's inequality one has
$$
\langle Q_{b}^T[H] Q_{b},{\rm Diag}\, (w_{b})  \rangle
\leq \lambda (Q_{b}^T[H] Q_{b})^Tw_{b},
$$
 which implies, from (\ref{eq:hlp10}), for $\lambda (Q_{b}^T[H] Q_{b})=(\lambda_1 (Q_{b}^T[H] Q_{b}),\ldots, \lambda_{|b|} (Q_{b}^T[H] Q_{b}))^T$ with $\lambda_1 (Q_{b}^T[H]  Q_{b})\geq \cdots \geq \lambda_{|b|} (Q_{b}^T[H] Q_{b})$,  that
 $$
 \langle Q_{b}^T[H] Q_{b},{\rm Diag}\, (w_{b})  \rangle
= \lambda (Q_{b}^T[H] Q_{b})^Tw_{b}=\|Q_{b}^T[H] Q_{b}\|_*.
 $$
 Then $Q_{b}^T[H] Q_{b}$ and ${\rm Diag}\, (w_{b})$ admit a simultaneous ordered eigenvalue decomposition, and thus we can check that $H$ satisfies
 $$
 \begin{array}{l}
 Q^{T}_{b_S}H[Q_{b}]=0,
 Q^{T}_{b_U}H[Q_{b_L}]=0,
   Q^{T}_{b_U}HQ_{b_U}\in {\cal S}^{|b_U|}_+,
          Q^{T}_{b_L}HQ_{b_L}\in {\cal S}^{|b_L|}_-.
          \end{array}
          $$
   The proof is completed. \hfill $\Box$\\
 \begin{corollary}\label{coro-an-equality}
 Let $X,Y,Z \in {\cal S}^q$, $Z=X+Y$ satisfies $Y \in \partial  \theta (X)$. Then $H \in {\cal C}_\theta (Z)$ if and only if
 \[
 \|Q_{a_s}^THQ_{a_s}\|_*=\langle Q_{a_s}^TYQ_{a_s}, Q_{a_s}^THQ_{a_s} \rangle.
 \]
 \end{corollary}

 \begin{prop}\label{prop-sigmaterm}
 Let $X,Y,Z,H \in {\cal S}^q$, $Z=X+Y$ satisfies $Y \in \partial  \theta (X)$ and $H \in {\cal C}_{\theta}(Z)$.  Then
 \begin{equation}\label{eq:conj1}
 \begin{array}{ll}
\psi^*(Y)= &
 2\displaystyle \sum_{i=1}^{s-1} {\rm Tr}(Q_{a_i}^TH(X-\varpi_iI)^{\dag}HQ_{a_i})
 +2\displaystyle \sum_{i=s+1}^{r} {\rm Tr}(Q_{a_i}^TH(X-\varpi_iI)^{\dag}HQ_{a_i})\\[14pt]
 & + 2{\rm Tr}[(Q_{b_U}^THX^{\dag}HQ_{b_U})]
 +2{\rm tr}[(Q_{b_L}^THX^{\dag}HQ_{b_L})]+2 \langle Q_{b_S}^T YQ_{b_S},Q_{b_S}^THX^{\dag}HQ_{b_S} \rangle.
 \end{array}
 \end{equation}
 \end{prop}
 {\bf Proof}. Since $H \in {\cal C}_{\theta}(Z)$, we have from (\ref{eq:critical-Prox}) that there exist $\widehat Q^s_U \in
 {\cal O}(Q^{T}_{b_U}HQ_{b_U})$ and $\widehat Q^s_L \in
 {\cal O}(Q^{T}_{b_L}HQ_{b_L})$ such that
 \[
 Q^s=\left
 [
 \begin{array}{ccc}
 \widehat Q^s_U & 0 &0\\[2pt]
 0 & I_{|b_S|} & 0\\[2pt]
 0 & 0 & \widehat Q^s_L
 \end{array}
 \right
 ].
 \]
 Then $Q^s_{b^k_+}$ and $Q^s_{b^k_+}$ can be expressed as
 \[
 Q^s_{b^k_+}=\left
 [
 \begin{array}{c}
 \widehat Q^s_U\\
 0\\
 0
 \end{array}
 \right
 ] \mbox{  and } Q^s_{b^k_-}=\left
 [
 \begin{array}{c}
 0\\
 0\\
 \widehat Q^s_L
 \end{array}
 \right
 ].
 \]
 Then we obtain (\ref{eq:conj1}) from (\ref{eq:conj}). \hfill $\Box$
 \begin{corollary}\label{coro-simple-version}
 Let $X,Y,Z,H \in {\cal S}^q$, $Z=X+Y$ satisfies $Y \in \partial  \theta (X)$ and $H \in {\cal C}_{\theta}(Z)$.  Then
 \begin{equation}\label{eq:conj1n}
\psi^*(Y)
 =2\displaystyle \sum_{i=1}^{r} \langle Q^T_{a_i}YQ_{a_i}, Q_{a_i}^TH(X-\varpi_iI)^{\dag}HQ_{a_i} \rangle,
 \end{equation}
 or alternatively
 \begin{equation}\label{eq:conj1ns}
\psi^*(Y)
 =2\displaystyle \sum_{i \ne s}\displaystyle \frac{1}{\varpi_i} \langle Q^T_{b_s}YQ_{b_s}-I_{|b_s|}, Q_{b_s}^THQ_{a_i}Q_{a_i}^THQ_{b_s} \rangle.
 \end{equation}
\end{corollary}

We now discuss the differential of $[e_{\tau}\theta](X)$  for $\theta(X)=\|X\|_*$, where  $[e_{\tau}\theta](X)$ is the Moreau-Yosida regularization defined by (\ref{eq:morea}). Let  proximal mapping of $\theta$ be defined by
$$
[{\cal P}_{\tau}\theta] (X)={\rm argmin}_{X' \in {\cal S}^q} \left\{\theta(X')+\displaystyle \frac{1}{\tau}\|X'-X\|^2\right\}.
$$
 For simplicity, we use ${\cal P}\theta$ to denote ${\cal P}_1\theta$. Then $[e_{\tau}\theta](X)$ is the spectral function corresponding to the  Moreau-Yosida regularization $e_{\tau}\varsigma$, namely
$$
[e_{\tau}\theta](X)=[e_{\tau}\varsigma \circ \lambda](X).
$$
It follows from \cite{Lewis95} or \cite{WDSun14} that
$$
[{\cal P}_{\tau}\theta](X)=P{\rm Diag}\,([{\cal P}_{\tau}\varsigma](\lambda(X)))P^T
$$
or $[{\cal P}_{\tau}\theta]( X)$ is the L\"{o}wner operator associated with $p_{\tau}(t)=[t-\tau]_+-[-t-\tau]_+$, namely
$$
[{\cal P}_{\tau}\theta]( X)=Q{\rm Diag}\left(p_{\tau}(\lambda_1( X)),\cdots,p_{\tau}(\lambda_q( X))\right) Q^T.
$$
Let $X$ have $r$ distinct eigenvalues, among them there are $r_1$ positive distinct eigenvalues and $r-r_1$ negative distinct eigenvalues and zero  eigenvalues:
$$
\varpi_1>\varpi_2 > \cdots > \varpi_{r_1}>\varpi_{r_1+1}=0> \varpi_{r_1+2} > \cdots > \varpi_r.
$$
Define
$$
a_k=\{i: \lambda_i(X)=\varpi_k\},\,\, k=1,\ldots, r
$$
and the first divided difference matrix at $X$ along $H \in{\cal S}^q$ as follows for $k,l=1,\ldots, r$,
\begin{equation}\label{eq:ptdir}
(p_{\tau}^{[1]}(\Lambda (X), Q^THQ))_{a_ka_l}:=
\left
\{
\begin{array}{ll}
\displaystyle \frac{p_{\tau}(\varpi_k)-p_{\tau}(\varpi_l)}{\varpi_k-\varpi_l} Q_{a_k}^THQ_{a_l} & \mbox{if } k \ne l,\\[16pt]
\Psi_k(Q_{a_k}^THQ_{a_k}) & \mbox{if }k = l.
\end{array}
\right.
\end{equation}
where $\Psi_k (\cdot)$ is the L\"{o}wner operator with respect to $\psi_k (\cdot)=p_{\tau}'(\varpi_k;\cdot)$.
Then the directional derivative of ${\cal P}_{\tau}\theta$ at $X$ along $H \in{\cal S}^q$ is expressed as
\begin{equation}\label{eq:drD-prox}
[{\cal P}_{\tau}\theta]'(X;H)=Q [p_{\tau}^{[1]}(\Lambda (X), Q^THQ)]Q^T.
\end{equation}
\begin{prop}\label{lem:proxSub}
 Let $X,Y,Z,H \in {\cal S}^q$, $Z_{\tau}=X+\tau Y$ satisfies $Y \in \partial  \theta (X)$. Then ${\cal P}_{\tau}\theta$ is strongly semismooth at $Z_{\tau}$ and
for $W\in \partial {\cal P}_{\tau}\theta (Z_{\tau})$, there exist $W_{b_U}\in \partial \Pi_{{\cal S}^{|b_U|}_+}(0)$ and $W_{b_L}\in \partial \Pi_{{\cal S}^{|b_L|}_-}(0)$  such that
\begin{equation}\label{eq:VH}
\begin{array}{l}
W(H)=\\[8pt]
Q\left
[
\begin{array}{ccccc}
\widehat H_{aa} & \widehat H_{ab_U} & \widehat H_{ab_S} \circ (\Omega_{\tau})_{ab_S}
&\widehat H_{ab_L} \circ (\Omega_{\tau})_{ab_L} & \widehat H_{ac} \circ (\Omega_{\tau})_{ac}\\[8pt]
\widehat H_{b_Ua} & W_{b_U}(\widehat H_{b_Ub_U}) & 0 & 0 & \widehat H_{b_Uc} \circ (\Omega_{\tau})_{b_Uc}\\[8pt]
\widehat H_{b_Sa} \circ (\Omega_{\tau})_{b_Sa}  & 0 & 0 & 0 & \widehat H_{b_Sc} \circ (\Omega_{\tau})_{b_Sc}\\[8pt]
\widehat H_{b_La}\circ (\Omega_{\tau})_{b_La} & 0 & 0 & W_{b_L}(\widehat H_{b_Lb_L})& \widehat H_{b_Lc} \\[8pt]
\widehat H_{ca}\circ (\Omega_{\tau})_{c a} & \widehat H_{cb_U}\circ (\Omega_{\tau})_{c b_U} & \widehat H_{cb_S}\circ (\Omega_{\tau})_{c b_S} & \widehat H_{cb_L}& \widehat H_{cc}
\end{array}
\right
]Q^T,
\end{array}
\end{equation}
where  $\widehat H=Q^THQ$,
$$
(\Omega_{\tau})_{ij}=[p_{\tau}^{[1]}(\Lambda (Z_{\tau}))]_{ij}, (i,j) \in a \times [b_S \cup b_L \cup c]
\mbox{ or } (i,j) \in c \times [ b_U \cup b_S].
$$
In other words,  $\nabla e_{\tau}\theta$ is strongly semismooth at $Z_{\tau}$ and
for $V \in \partial \nabla e_{\tau}\theta (Z_{\tau})$, there exist $V_{b_U}\in \partial \Pi_{{\cal S}^{|b_U|}_-}(0)$ and $V_{b_L}\in \partial \Pi_{{\cal S}^{|b_L|}_+}(0)$  such that
\begin{equation}\label{eq:VH}
\begin{array}{l}
\tau V(H)=\\[8pt]
Q\left
[
\begin{array}{ccccc}
 0 & 0 & \widehat H_{ab_S} \circ (\Delta_{\tau})_{ab_S}
&\widehat H_{ab_L} \circ (\Delta_{\tau})_{ab_L} & \widehat H_{ac} \circ (\Delta_{\tau})_{ac}\\[8pt]
0 & V_{b_U}(\widehat H_{b_Ub_U}) & \widehat H_{b_Ub_S} & \widehat H_{b_Ub_L} & \widehat H_{b_Uc} \circ (\Delta_{\tau})_{b_Uc}\\[8pt]
\widehat H_{b_Sa} \circ (\Delta_{\tau})_{b_Sa}  &\widehat H_{b_Sb_U} & \widehat H_{b_Sb_S} & \widehat H_{b_Sb_L} & \widehat H_{b_Sc} \circ (\Delta_{\tau})_{b_Sc}\\[8pt]
\widehat H_{b_La}\circ (\Delta_{\tau})_{b_La} & \widehat H_{b_Lb_U} & \widehat H_{b_Lb_S} & V_{b_L}(\widehat H_{b_Lb_L})& 0 \\[8pt]
\widehat H_{ca}\circ (\Delta_{\tau})_{c a} & \widehat H_{cb_U}\circ (\Delta_{\tau})_{c b_U} & \widehat H_{cb_S}\circ (\Delta_{\tau})_{c b_S} & 0 & 0
\end{array}
\right
]Q^T,
\end{array}
\end{equation}
where  $\widehat H=Q^THQ$,
\begin{equation}\label{eq:Delta}
(\Delta_{\tau})_{ij}=1-[p_{\tau}^{[1]}(\Lambda (Z_{\tau}))]_{ij}, (i,j) \in a \times [b_S \cup b_L \cup c]
\mbox{ or } (i,j) \in c \times [ b_U \cup b_S].
\end{equation}
\end{prop}

\subsection{Optimality conditions for (SDNOP)}
This subsection is devoted to studying optimality conditions for the following nonlinear
semidefinite nuclear norm composite optimization problem
\vspace{\baselineskip}\noindent (SDNOP)
 \vspace{-\baselineskip}
 \[
     \min \ f(x)+\theta (F(x))\quad
 {\rm s.t.} \  h(x)=0\, ,
 \  g(x) \in {\cal S}^p_+\, ,
\]
 where $\theta (X)=\|X\|_*$ is the nuclear norm function of $X \in {\cal S}^q$, $f:\Re^n\mapsto \Re$, $F: \Re^n \mapsto {\cal S}^q$, $h:\Re^n \mapsto \Re^m$ and
$g:\Re^n \mapsto {\cal S}^p$ are twice continuously differentiable functions.
Obviously, Problem (SDNOP) is a special case of (COP) with ${\cal Z}:={\cal S}^q$, $\theta (X):=\|X\|_*$,
  ${\cal Y}:={\cal S}^p$ and $K:={\cal S}^p_+$.
The  Lagrange  function  for (SDNOP) is
\[L(x,Y,\mu,\Gamma)=f(x)+\langle Y, F(x)\rangle+
\langle \mu, h(x)\rangle- \langle \Gamma, g(x)\rangle\,\, (x,
Y,\mu, \Gamma)\in \Re^n \times {\cal S}^q \times \Re^m \times {\cal S}^p.
\]
Then for any $(x, Y,\mu, \Gamma)\in \Re^n \times {\cal S}^q \times \Re^m \times {\cal
S}^p$,
\[\nabla_xL(x,Y,\mu,\Gamma) =\nabla f(x) +{\rm D}F(x)^*Y+ {\cal J} h(x)^T\mu
-{\rm D}g(x)^* \Gamma.
\]
If $x$ is a stationary point, the set of Lagrange multipliers at $x$ is defined by
\[
\Lambda (x)=\left\{(Y,\mu, \Gamma)\in {\cal S}^q \times \Re^m \times {\cal
S}^p:\begin{array}{l}
\nabla_x L(x,Y,\mu,\Gamma)=0, \\
Y \in \partial \theta(F(x)), \Gamma \in -{\cal N}_{{\cal S}^p_+}(g(x))
\end{array}
 \right\}.
\]
When discussing optimality conditions, we need some constraint qualifications. We say that
Robinson  constraint
qualification holds at $\overline{x}$ if
\[
\left(\begin{array}{c}
{\cal J}h(\overline{ x})
\\
{\rm D}g(\overline{ x})
\end{array} \right)\Re^n
 +\left(\begin{array}{c}
 \{0\}
\\
\left({\cal T}_{ {\cal S}^p_+}(g(\overline{ x})) \right)
\end{array}\right)  = \left(\begin{array}{c}
\Re^m\\
{\cal S}^p
\end{array}
\right).
\]

The critical cone of  Problem (SDNOP) at $x$ is defined by
\[
{\cal C}(x)=\{d\in {\cal T}_{\Phi}(x): \nabla f(x)^Td+\theta'(F(x);{\rm D}F(x)d)\leq 0\}.
\]
We can easily derive the following necessary optimality conditions and second-order sufficient optimality conditions.
\begin{prop}\label{nec-th}
If $\overline x \in \Phi$ is a local minimizer around which $f,F,h$ and $g$ are twice continuously differentiable and Robinson constraint qualification holds at $\overline x$. Then
\begin{itemize}
\item[(1)] $\Lambda (\overline x)$ is non-empty, compact and convex.
\item[(2)] For any $d \in {\cal C}(\overline x)$,
\[
\displaystyle \sup_{y \in \Lambda (\overline x)} \left\{\left\langle d, \nabla^2_{xx}L(\overline{x},y)d \right \rangle-\psi^*(\overline Y)+
 2\left \langle \overline \Gamma,[{\rm D}g(\overline x)d]g(\overline x)^\dag [{\rm D}g(\overline x)d]\right
\rangle\right\}
\geq 0,
\]
where $\psi (W)=\theta''(F(\overline x);{\rm D}F(\overline x)d,W)$.
\end{itemize}
\end{prop}
\begin{prop}\label{suf-th}
 Let $\overline x$ be a feasible point  around which $f,F,h$ and $g$ are twice continuously differentiable. Suppose the following conditions hold:
\begin{itemize}
\item[(1)] $\Lambda (\overline x)$ is non-empty;
\item[(2)] For any $d \in {\cal C}(\overline x)\setminus \{0\}$,
\[
\displaystyle \sup_{y \in \Lambda (\overline x)}\left \{\left\langle d, \nabla^2_{xx}L(\overline{x},y)d \right \rangle-\psi^*(\overline Y)+
 2\left \langle \overline \Gamma,[{\rm D}g(\overline x)d]g(\overline x)^\dag [{\rm D}g(\overline x)d]\right
\rangle\right\}
> 0,
\]
where $\psi (W)=\theta''(F(\overline x);{\rm D}F(\overline x)d,W)$.
\end{itemize}
Then the second-order growth condition holds at $\overline x$.
\end{prop}
Now we list our two assumptions for Problem (SDNOP), which will be used in the next section to derive Assumptions B1 and B2.\\
\noindent {\bf Assumption (sdnop-A1)}\cite{CDing2017}. The constraint
nondegeneracy condition holds at $\overline{x}$:
\begin{equation}\label{eq:nondegeracy-SDP}
\left(\begin{array}{c}
{\rm D}F(\overline x)\\
{\cal J}h(\overline{ x})
\\
{\rm D}g(\overline{ x})
\end{array} \right)\Re^n
 +\left(\begin{array}{c}
 {\cal T}^{{\rm lin}}(F(\overline x))\\
 \{0\}
\\
{\rm lin} \left({\cal T}_{ {\cal S}^p_+}(g(\overline{ x})) \right)
\end{array}\right)  = \left(\begin{array}{c}
{\cal S}^q\\
\Re^m\\
{\cal S}^p
\end{array}
\right),
\end{equation}
where
\begin{equation}\label{eq:Tlin}
{\cal T}^{{\rm lin}}(\overline X)=\{H \in {\cal S}^q: \theta'(\overline X;H)=-\theta'(\overline X;-H)\}=\{H \in {\cal S}^q: Q_{b}^THQ_{b}=0\}.
\end{equation}

Assumption (sdnop-A1) is the analogue to the linear independence
constraint qualification for nonlinear programming, which implies that ${\cal
M}(\overline{x})$ is a singleton \cite[Proposition 4.50]{BS00}.

\noindent {\bf Assumption (sdnop-A2)} The strong second order
sufficient condition holds at $\overline{x}$ :
\[
\left \langle d, \nabla^2_{xx}L(\overline{x},\overline Y,
\overline{\mu},\overline{\Gamma})d \right \rangle-\psi^*(\overline Y)+
 2\left \langle \overline \Gamma, [{\rm D}g(\overline x)d]g(\overline x)^\dag [{\rm D}g(\overline x)d]\right
\rangle
>0, \quad  \forall \, d \in {\rm app}(\overline Y,\overline{\mu},\overline{\Gamma}) \setminus \{0\}\, ,
\]
where $\psi (W)=\theta''(F(\overline x);{\rm D}F(\overline x)d,W)$ and
\begin{equation}\label{equation3.5}
\mbox{app} (\overline Y, \overline{\mu},\overline{\Gamma}):=\l\{d\in \Re^n:
\begin{array}{l} {\cal
J} h(\overline{x})d=0, \, {\rm D}F(\overline x)d \in \mbox{aff}({\cal C}_{\theta}(F(\overline x)+\overline Y))\\
{\rm D} g(\overline{x})d \in \mbox{aff}({\cal C}_{\S^p_+}(g(\overline{x})-\overline{\Gamma}))
\end{array}\r\}.
\end{equation}
From the expressions ${\cal C}_{\theta}$ and ${\cal C}_{\S^p_+}$, we obtain the following expression of ${\rm app }(\overline Y, \overline{\mu},\overline{\Gamma})$:
\begin{equation}\label{equation3.5a}
{\rm app }(\overline Y, \overline{\mu},\overline{\Gamma})=\left\{d\in \Re^n:
\begin{array}{l}
  Q_{b_S}^T({\rm D}F(\overline x)d)Q_b=0, Q_{b_U}^T({\rm D}F(\overline x)d)Q_{b_L}=0\\[6pt]
 {\overline P}_{\alpha}^T({\cal J}
g(\overline{x})d){\overline P}_{\alpha}=0, {\overline
P}_{\alpha}^T({\cal J} g(\overline{x})d){\overline P}_{\beta}=0,{\cal J}
h(\overline{x}) d =0
\end{array}
\right\}.
\end{equation}

\vskip 7 true pt
 At the end of this subsection, we list two technical results coming from \cite{SSZhang2008}, which will be used in the next section.
\begin{lemma}\label{lemma3.2nlp}\cite[Lemma 7]{SSZhang2008}
Let $\phi: {\cal X} \mapsto \Re$ be continuous and positive homogeneous of degree two:
\[
\phi (t d)=t^2 \phi (d),  \quad \forall \, t \geq 0 \ {\rm and} \ d
\in {\cal X} \, .
\]
Suppose that there exists a positive number $\eta_0>0$ such that for
any $d$ satisfying ${\cal L}d=0$, one has $\phi (d) \geq \eta_0
\|d\|^2$,  where ${\cal L}: {\cal X}  \mapsto {\cal Y} $ is a given linear operator.
Then there exist positive numbers $\underline{\eta} \in (0, \eta_0]$
and $c_0
>0$ such that
\[
\phi (d)+c_0 \langle {\cal L}d, {\cal L}d\rangle \geq
\underline{\eta} \langle d, d\rangle,  \quad  \forall \, d \in {\cal X} \, .
\]
\end{lemma}
\begin{lemma}\label{phiineq}\cite[Lemma 8]{SSZhang2008}
Let $a,b,c,$ and $c_0$ be four positive scalars  with $c\ge c_0$.
Let
\begin{equation}\label{phicom}
\psi(t;c,a,b,c_0):=a- \frac{1}{c}t+\frac{t^2}{b+(c-c_0)t}\, ,
\quad t \in [0,1]\, .
\end{equation}
 Then, for any  $c \geq  \max \big\{c_0,(b-c_0)^2/c_0\big\}$,
 $\psi (\cdot;c,a,b,c_0)$ is a convex function on $[0,1]$,
  \begin{equation}\label{k4com}
\min_{t \in [0,1]} \psi(t;c,a,b,c_0)=a-\frac{1}{c}
\frac{b}{(\sqrt{c}+\sqrt{c_0})^2}\, ,
\end{equation}
and
\begin{equation}\label{u4com}
\max_{t \in [0,1]}  \psi(t;c,a,b,c_0)= \max
\Big\{\psi(0;c,a,b,c_0), \, \psi(1; c,a,b,c_0)\Big\} \, .
\end{equation}
\end{lemma}
\section{ On the augmented Lagrange method  for SDNOP
}\label{NLSDP-case}
 \setcounter{equation}{0}

 This section is devoted to studying the rate of convergence of the augmented Lagrange method  for Problem (SDNOP).
 Let $(\overline{x}, \overline Y, \overline{\mu},
\overline{\Gamma}) \in \Re^n \times {\cal S}^q \times \Re^m \times {\cal S}^p$ be a given
KKT point. Then, $(\overline{x},\overline Y, \overline{\mu},
\overline{\Gamma})$ satisfies
\begin{equation}\label{nlsdp:KKT}
\nabla_xL(\overline{x},\overline Y, \overline{\mu},
\overline{\Gamma})=0, \ \overline Y \in \partial \theta (F(\overline x),\  h(\overline{x})=0, \  \overline{\Gamma} \succeq 0,
\  g(\overline{x})\succeq 0\,\, {\rm and}\,\, \langle
\overline{\Gamma}, g(\overline{x})\rangle =0.
\end{equation}
Let $\overline X=F(\overline x)$ and $\overline Y\in \partial \theta (\overline X)$.
 Define the following three index sets:
$$
a=\{i: \lambda_i(\overline X)>0\},\,\, b=\{i: \lambda_i(\overline  X)=0\},\,\, c=\{i: \lambda_i(\overline  X)<0\},
$$
then
\[
\overline X=Q  \left[ \begin{array}{ccc} \Lambda_a & 0 & 0 \\
[5pt] 0 & 0 & 0 \\ [5pt] 0 & 0 & \Lambda_c
\end{array} \right]Q^T \epc \mbox{and  } Q \in {\cal O}(\overline X) \mbox{ with }
Q \, = \, [ \begin{array}{ccc}Q_a &
Q_{b} & Q_{c}
\end{array} ]
\]
with $Q_{a} \in \Re^{q \times |a|}$,
$Q_{b} \in \Re^{q \times |b|}$, and
$Q_{c} \in \Re^{q \times |c|}$.
Then there exists $w \in \partial \varsigma (\lambda (\overline  X))$ satisfying $\overline Y=Q{\rm Diag}(w)Q^T$ and $w$ has the following relations
$$
w_a =\textbf{1}_{|a|}, w_{c}=-\textbf{1}_{|c|} \mbox{ and }  -\textbf{1}_{|b|}\leq w_{b}\leq \textbf{1}_{|b|}.
$$
For the index set $b$, we partition it as follows
$b=b_L \cup b_S \cup b_U$:
$$
 b_L=\{i\in b: w_i=-1\}, b_S=\{i\in b: -1< w_i< 1\},b_U=\{i\in b: w_i=1\}.
$$
Then $\overline Y$ can be expressed as follows:
\begin{equation}\label{eq:dec-Ya}
\overline Y=
\left (Q_{a\cup_{b_U}} \quad Q_{b_S} \quad Q_{c\cup_{b_L}} \right)
 \left[ \begin{array}{ccc}
 I_{|a\cup_{b_U}|} & 0& 0\\
 0 & {\rm Diag}\,(w_{b_S}) & 0 \\
 0 & 0 & I_{|c\cup_{b_L}|}
 \end{array}
 \right
 ]
 \left
 (
 \begin{array}{c}
 Q_{a\cup_{b_U}}^T\\[4pt]
   Q^T_{b_S}\\[4pt]
   Q_{c\cup_{b_L}}^T
   \end{array}
  \right)
\end{equation}
with $Q_{a\cup_{b_U}} \in \Re^{q \times |a\cup_{b_U}|}$,
$Q_{b_S} \in \Re^{q \times |b_S|}$, and
$Q_{c\cup_{b_L}} \in \Re^{q \times |c\cup_{b_L}|}$.

Let $\overline{M}: =\overline{\Gamma}-g(\overline{x})$. Suppose that
$\overline{M}$ has the  spectral decomposition as in
(\ref{eq:orthogonal-decomposition}), i.e, $\overline{M}  =  P \Lambda
P^T$.
Define three index sets of positive, zero, and negative
eigenvalues of $\overline{M} $, respectively, as
\[
\alpha  :=  \{  i \, | \, \lambda_i  > \, 0  \}, \epc \beta  := \{
i \, | \, \lambda_i \, = \, 0  \}, \epc \gamma  :=  \{  i \, | \,
\lambda_i \, < \, 0 \}.
\]
 Write
\[
\Lambda  =  \left[ \begin{array}{ccc} \Lambda_{\alpha} & 0 & 0 \\
[5pt] 0 & 0 & 0 \\ [5pt] 0 & 0 & \Lambda_{\gamma}
\end{array} \right] \epc \mbox{and} \epc
P \, = \, [ \begin{array}{ccc} P_{\alpha} &
P_{\beta} & P_{\gamma}
\end{array} ]
\]
with $P_{\alpha} \in \Re^{p \times |\alpha|}$,
$P_{\beta} \in \Re^{p \times |\beta|}$, and
$P_{\gamma} \in \Re^{p \times |\gamma| }$. {}From
(\ref{nlsdp:KKT}), we know that $\overline{\Gamma} g(\overline{x}) =
g(\overline{x}) \overline{\Gamma}=0$. Thus, we have
\[
 \overline{\Gamma} = P
\left[
\begin{array}{ccc} \Lambda_{\alpha} & 0 & 0 \\ [5pt] 0 & 0 & 0 \\
[5pt] 0 & 0 & 0
\end{array} \right] P^T \, ,  \quad
g(\overline{x}) = P \left[ \begin{array}{ccc} 0& 0 & 0
\\ [5pt] 0 & 0 & 0
\\ [5pt] 0 & 0 & -\Lambda_{\gamma}
\end{array} \right] P^T
\]
\begin{equation}\label{eq:Zc}
\overline{\Gamma}- t g(\overline{x}) =P\left[
\begin{array}{ccc} \Lambda_{\alpha} & 0 & 0 \\ [5pt] 0 & 0 & 0 \\
[5pt] 0 & 0 & t\Lambda_{\gamma}
\end{array} \right] P^T \, .
\end{equation}
For $\overline X=F(\overline x)$, let
\begin{equation}\label{eq:two-nus}
\begin{array}{ll}
\underline{\nu}_{a,b_S} :=\displaystyle\min_{i\in a, 1\leq j\leq |b_S|}\displaystyle  \frac{1-(w_b)_j}{\lambda_i(\overline X)}
&  \overline{\nu}_{\alpha,\gamma} :=\max_{i\in a, 1\leq j\leq |b_S|} \displaystyle  \frac{1-(w_b)_j}{\lambda_i(\overline X)};\\[4mm]
\underline{\nu}_{a,b_L} :=\displaystyle \min_{i\in a} \displaystyle \frac{2}{\lambda_i(\overline X)}
&  \overline{\nu}_{a,b_L} :=\max_{i\in a} \displaystyle \frac{2}{\lambda_i(\overline X)};\\[4mm]
\underline{\nu}_{a,c} :=\displaystyle\min_{i\in a, j\in c}\displaystyle \frac{2}{\lambda_i(\overline X)-\lambda_j(\overline X)}
& \overline{\nu}_{a,c} :=\max_{i\in a, j\in c}\displaystyle \frac{2}{\lambda_i(\overline X)-\lambda_j(\overline X)};\\[4mm]
\underline{\nu}_{c,b_U} :=\displaystyle\min_{i\in c}\displaystyle \frac{2}{-\lambda_i(\overline X)}
& \overline{\nu}_{c,b_U} :=\displaystyle\max_{i\in c}\displaystyle \frac{2}{-\lambda_i(\overline X)};\\[4mm]
\underline{\nu}_{c,b_S} :=\displaystyle\min_{i\in c, 1\leq j\leq |b_S|}\displaystyle \frac{1+(w_{b_S})_j}{-\lambda_i(\overline X)}
&  \overline{\nu}_{c,b_U} :=\displaystyle\max_{i\in c, 1\leq j\leq |b_S|}\displaystyle \frac{1+(w_{b_S})_j}{-\lambda_i(\overline X)};\\[4mm]
\underline{\nu}_{\alpha,\gamma} :=\displaystyle \min_{i\in\alpha, j\in \gamma} \lambda_i
/|\lambda_j| & \overline{\nu}_{\alpha,\gamma} :=\displaystyle \max_{i\in\alpha, j\in
\gamma} \lambda_i / |\lambda_j|
\end{array}
\end{equation}
and
\begin{equation}\label{nu0}
\begin{array}{l}
\underline{\nu}_0=\min\{\underline{\nu}_{a,b_S},\underline{\nu}_{a,b_L},\underline{\nu}_{a,c},
\underline{\nu}_{c,b_U},\underline{\nu}_{c,b_S},\underline{\nu}_{\alpha,\gamma}\};\\[4mm]
\overline{\nu}_0=\min\{\overline{\nu}_{a,b_S},\overline{\nu}_{a,b_L},\overline{\nu}_{a,c},
\overline{\nu}_{c,b_U},\overline{\nu}_{c,b_S},\overline{\nu}_{\alpha,\gamma}\}.
\end{array}
\end{equation}

For a given symmetric
matrix $M$, we use  ${\rm vec}(M)$  to denote the
 vector
obtained by  stacking up all the columns of a given matrix $M$ and
${\rm svec}(M)$ to denote the vector obtained by  stacking up all
the columns of the upper triangular part of $M$.

Let  $Q\in {\cal O}(\overline X)$ with $Q=[Q_a\,\, Q_{b_U}\,\, Q_{b_S} \,\, Q_{b_L} \,\, Q_c]$. For index sets
$\chi, {\chi'} \in \{a, b_U,b_S,b_L,c\}$, let
\[
B_{(\chi,\chi')}(Q):=\Big(\mbox{vec}(Q_{\chi}^T{\cal J}_{x_1}
F(\overline{x})Q_{\chi'}) \,\, \cdots \,\,
\mbox{vec}(Q_{\chi}^T{\cal J}_{x_n}
F(\overline{x})Q_{\chi'})\Big)\,
\]
and
\[
{\widehat B}_{(\chi,\chi)}(Q):=\Big(\mbox{svec}(Q_{\chi}^T{\cal
J}_{x_1} F(\overline{x})Q_{\chi} )\,\, \cdots \,\,
\mbox{svec}(Q_{\chi}^T{\cal J}_{x_n} F(\overline{x})Q_{\chi}
)\Big).
\]
Let $P\in {\cal O}(g(\overline x))$ with $
P=[ P_\alpha\ P_\beta \
 P_\gamma]$. For index sets
$\chi, {\chi'} \in \{\alpha,\beta,\gamma\}$, let
\[
C_{(\chi,\chi')}(P):=\Big(\mbox{vec}(P_{\chi}^T{\cal J}_{x_1}
g(\overline{x})P_{\chi'}) \,\, \cdots \,\,
\mbox{vec}(P_{\chi}^T{\cal J}_{x_n}
g(\overline{x})P_{\chi'})\Big)\,
\]
and
\[
{\widehat C}_{(\chi,\chi)}(P):=\Big(\mbox{svec}(P_{\chi}^T{\cal
J}_{x_1} g(\overline{x})P_{\chi} )\,\, \cdots \,\,
\mbox{svec}(P_{\chi}^T{\cal J}_{x_n} g(\overline{x})P_{\chi}
)\Big).
\]
 Define
\[
n_1:=m+|b|(|b|+1)/2 \, ,  \
n_2:=n_1+(|\alpha|+|\beta|)(|\alpha|+|\beta|+1)/2 \, , \ n_3:=n-n_2\,
,
\]
and
\[ A(Q,P):=\left ( \begin{array}{c}
 {\cal J} h(\overline{x})\\
 {\widehat B}_{(b_U,b_U)}(Q)\\
 B_{(b_U,b_S)}(Q)\\
 B_{(b_U,b_L)}(Q)\\
 {\widehat B}_{(b_S,b_S)}(Q)\\
 B_{(b_S,b_L)}(Q)\\
 {\widehat B}_{(b_L,b_L)}(Q)\\
 -{\widehat C}_{(\alpha,\alpha)}(P)\\
-{\widehat C}_{(\beta,\beta)}(P)\\
- C_{(\alpha,\beta)}(P)
 \end{array}
 \right
 )\, .
\]
Suppose that Assumption (sdnop-A1) holds. Then by
(\ref{eq:nondegeracy-SDP}) in Assumption (sdnop-A1) we know that
$A(Q,P)$ is of full row rank. Let $A(Q,P)$ have the following
singular value decomposition:
\begin{equation}\label{singvalu}
 A(Q,P)=U[\Sigma (Q,P)\,\,\,\,\, 0]R^T\, ,
\end{equation}
where $U \in \Re^{n_2 \times n_2}$ and $R \in \Re^{n \times n}$ are
orthogonal matrices, $\Sigma(Q,P)=\mbox{Diag}
\Big(\sigma_1(A(Q,P)),\cdots,$  $ \sigma_{n_2}(A(Q,P))\Big)$, and
$\sigma_1(A(Q,P))\geq \sigma_2(A(Q,P)) \geq $ $ \cdots $  $\geq
\sigma_{n_2}(A(Q,P))>0$ are the singular values of $A(Q,P)$. It should
be pointed out here that $U$ and $R$ also depend on $(Q,P)$. But for the
sake of notational simplification, we drop the argument $(Q,P)$ from $U$
and $R$ in our analysis below.

Let
\[
\underline{\sigma}:=\min\left\{1, \min_{Q \in {\cal O}(\overline X),P \in {\cal O}(\overline M)}  \min_{1
\leq i \leq n_2} \sigma_i^{-2}(A(Q,P))\right\}
\]
 and
 \[
\overline{\sigma}:=\max\left\{1, \max_{Q \in {\cal O}(\overline X),P \in {\cal O}(\overline M)} \max_{1
\leq i \leq n_2}\sigma_i^{-2}(A(Q,P))\right \}\, .
\]
Then, since ${\cal O}(\overline X)$ and ${\cal O}(\overline M)$ are  compact sets and $\Sigma(Q,P)$ changes
continuously with  respect to $(Q,P)$,  both $\underline{\sigma}$ and
$\overline{\sigma}$ are finite positive numbers.
Define
\[
C=\left
[
\begin{array}{c}
 B_{(a,b_S)}\\
  B_{(a,b_L)}\\
   B_{(a,c)}\\ B_{(c,b_U)}\\
   B_{(c,b_S)}\\
- C_{(\alpha,\gamma)}\\
 \end{array}
 \right].
 \]
 Thus there exist
 numbers $\underline{\nu} \geq 0$ and $\overline{\nu}>0$ such
that for any $Q\in {\cal O}(\overline X)$, $P\in {\cal O}(\overline M)$ and  $s\in \Re^{|a||b_S|+|a||b_L|+|a||c|+|c||b_U|+|\alpha||\gamma|}$,
\begin{equation}\label{eq:nu-lower-upper}
\underline{\nu} \|s\|^2 \leq \max\left\{ \left\langle s,
\widetilde{C}(Q,P)(\widetilde{C}^{\,
T}(Q,P)) s \right\rangle, \left\langle s,
C C^Ts \right\rangle
\right\}
 \le \overline{\nu}\|s\|^2\, ,
\end{equation}
where
\[\widetilde{C}(Q,P):=C\widetilde{R}
\quad {\rm and} \quad    \widetilde{R}:=R\left [
\begin{array}{cc}
\Sigma(Q,P)^{-1}U^T & 0 \\ 0 & I_{n_3}
\end{array}
\right ].
\]
 When no
ambiguity arises, we often  drop $Q$ and $P$ from  $A(Q,P)$, $B_{(\chi,\chi')}(Q)$, $\widehat{B}_{(\alpha,\gamma)}(Q)$.
$C_{(\chi,\chi')}(Q)$,  and $\widehat{C}_{(\alpha,\gamma)}(P)$.
Let $c>0$ and $W_1\in \partial_B [{\rm D} \theta_c]^*(F(\overline x)+\overline Y/c)$, there exist
matrices $Q\in {\cal O}(F(\overline x))$ and $\Delta_{1/c} \in {\cal S}^q$ such that
\begin{equation}\label{eq:w1h}
W_1(H_1)= Q\left( \Delta_{1/c} \circ (Q^TH_1Q)\right )Q^T,\quad \forall \, H_1
\in {\cal S}^q\, .
\end{equation}
with the entries of $\Delta_{\tau}$ being given by
\begin{equation}\label{Shc}
\left\{ \begin{array}{ll}
 (\Delta_{\tau})_{ij}  =1-[p_{\tau}^{[1]}(\Lambda (Z_{\tau}))]_{ij}, & (i,j) \in a \times [b_S \cup b_L \cup c]
\mbox{ or } (i,j) \in c \times [ b_U \cup b_S]
\\[2mm]
 (\Delta_{\tau})_{ij}   \in [0,1], &  (i, j) \in
b_U\times b_U  \mbox{ or } (i,j) \in b_L\times b_L.
\end{array} \right.
\end{equation}
It can be easily verified, for $\overline X=F(\overline x)$, that
\begin{equation}\label{eq:dDeltac}
[{\Delta_{1/c}}]_{ij}=\left
 \{
\begin{array}{ll}
\quad \quad 0 & (i,j)\in (a\times a \cup b_U)\cup (c \times b_L\cup c),\\[12pt]
\displaystyle \frac{c^{-1}(1-(w_{b_S})_j)}{\lambda_i(\overline X)+c^{-1}(1-(w_{b_S})_j)} & (i,j)\in a\times \{1,\ldots,|b_S|\},\\[12pt]
\displaystyle \frac{2c^{-1}}{\lambda_i(\overline X)+2c^{-1}} & (i,j)\in a\times b_L,\\[12pt]
\displaystyle \frac{2c^{-1}}{\lambda_i(\overline X)-\lambda_j(\overline X)+2c^{-1}} & (i,j)\in a\times c,\\[12pt]
\displaystyle \frac{2c^{-1}}{-\lambda_j(\overline X)+2c^{-1}} & (i,j)\in b_U\times c,\\[12pt]
\displaystyle \frac{c^{-1}((w_{b_S})_i+1)}{c^{-1}((w_{b_S})_i+1)-\lambda_j(\overline X)} & (i,j)\in \{1,\ldots,|b_S|\}\times c.
\end{array}
\right.
\end{equation}
Let $c>0$ and  $W_2 \in \partial_B \Pi_{{\cal
S}^p_+}(\overline{\Gamma}-cg(\overline{x}))$.  Define $\lambda_c \in
\Re^p$ as
 \[
 (\lambda_c)_i := \left\{
 \begin{array}{ll} \lambda_i & {\rm if} \
i\in \alpha \cup \beta\, ,
\\
c \lambda_i &  {\rm if} \ i\in \gamma\, .
\end{array} \right.
\]
Then  it follows from Lemma \ref{lemma:projector-Jacobian-cone}
that there exist two matrices  $Q\in {\cal O}(\overline M)$ and
$\Theta_c \in {\cal S}^p$ such that
\begin{equation}\label{wh}
W_2(H_2) = P\left( \Theta_c \circ (P^TH_2P)\right )P^T,\quad \forall \, H_2
\in {\cal S}^p\,
\end{equation}
with the entries of $\Theta_c$ being given by
\begin{equation}\label{Sh}
\left\{ \begin{array}{ll}
 (\Theta_c)_{ij}  =\displaystyle
{\frac{\max\{(\lambda_c)_i,0\} +  \max\{(\lambda_c) _j,0\}}{ | \,
(\lambda_c)_i \, | + | \, (\lambda_c)_j \, |}}& {\rm if}\ (i, j)
\notin \beta\times \beta\, ,
\\[2mm]
 (\Theta_c)_{ij}  \in [0,1] & {\rm if}\ (i, j) \in
\beta\times \beta\, .
\end{array} \right.
\end{equation}
For index sets $\chi, {\chi'} \in \{a,b_U,b_S,b_L,c\}$, we
introduce the following notation:
\[
(\Delta_{\tau})_{(\chi,\chi')}=\mbox{Diag} \left
(\mbox{vec}((\Delta_{\tau})_{\chi\chi'} )\right)\, , \ \ (\widehat
\Delta_{\tau})_{(\chi,\chi)}=\mbox{Diag} \left
(\mbox{svec}((\Delta_{\tau})_{\chi\chi} \circ E_{\chi\chi})\right),
\]
where  $``\circ"$ is the Hadamard product and $E$ is a matrix in
${\cal S}^q$ with entries being given by
\[
E_{ij} :=\left \{
\begin{array}{ll}
1 & {\rm if} \ i=j\, ,
\\
{2} & {\rm if} \ i\ne j\, .
\end{array}
\right.
\]
For index sets $\chi, {\chi'} \in \{\alpha,\beta,\gamma\}$, we
introduce the following notation:
\[
(\Theta_c)_{(\chi,\chi')}=\mbox{Diag} \left
(\mbox{vec}((\Theta_c)_{\chi\chi'} )\right)\, , \ \ ({\widehat
\Theta}_c)_{(\chi,\chi)}=\mbox{Diag} \left
(\mbox{svec}((\Theta_c)_{\chi\chi} \circ E'_{\chi\chi})\right)\, ,
\]
where   $E'$ is a matrix in
${\cal S}^p$ with entries being given by
\[
E'_{ij} :=\left \{
\begin{array}{ll}
1 & {\rm if} \ i=j\, ,
\\
{2} & {\rm if} \ i\ne j\, .
\end{array}
\right.
\]
 Let
\[
 D_c:=\left [
\begin{array}{ccccc}
I_{m} & 0 & 0 & 0 &0\\
0 & \Sigma_c & 0 & 0 &0\\
0 & 0 &({\widehat \Theta_c})_{(\alpha,\alpha)} & 0 & 0\\
0 & 0& 0& ({\widehat \Theta_c})_{(\beta,\beta)} &0  \\
 0 & 0 & 0& 0& 2I_{|\alpha||\beta|}
\end{array}
\right ],
\]
 where
 \[
  \Sigma_c=\left
  [
  \begin{array}{ccc}(\widehat \Delta_{1/c})_{b_Ub_U} & 0 & 0\\
 0 & 2I_{m_0} & 0\\
 0 & 0 &  ( \widehat \Delta_{1/c})_{b_Lb_L}
 \end{array} \right ]
 \]
 with $m_0=|b_U|(|b_S|+|b_L|)+|b_S|((|b_S|+1)/2+|b_L|)$.\\

 Let ${\cal A}_c(\overline Y,\overline{\mu}, \overline{\Gamma},W_1,W_2)$ be defined as
 (\ref{comac}) for the semidefinite nuclear norm composite optimization problem (SDNOP), i.e,
\[
\begin{array}{ll}
{\cal A}_c(\overline Y,\overline{\mu}, \overline{\Gamma},W_1,W_2)=& \nabla^2_{xx}
L(\overline{x},\overline Y,\overline{\mu}, \overline{\Gamma})\\[3mm]
&+c  {\cal  J}
h(\overline{x})^T{\cal J} h(\overline{x}) +{\rm D}F(\overline x)^*W_1{\rm D}F(\overline x)+
c{\rm D}g(\overline{x})^*W_2{\rm D} g(\overline{x}).
\end{array}
\]
A compact formula for ${\cal A}_c(\overline Y,\overline{\mu}, \overline{\Gamma},W_1,W_2)$ is given in the next lemma.
\begin{lemma}\label{acpres}
The matrix ${\cal A}_c(\overline Y,\overline{\mu}, \overline{\Gamma},W_1,W_2)$ can be
expressed equivalently as
\begin{equation}\label{acy}
\begin{array}{ll}{\cal A}_c(\overline Y,\overline{\mu}, \overline{\Gamma},W_1,W_2) =&  \nabla^2_{xx}
L(\overline{x},\overline Y,\overline{\mu}, \overline{\Gamma})+c \left( {\cal J} h(\overline{x})^T
{\cal J} h(\overline{x})+2 B_{(b_U,b_S)}^{\, T} B_{(b_U,b_S)}\right.\\[2mm]
& \left.+2 B_{(b_U,b_L)}^{\, T} B_{(b_U,b_L)}
+\widehat B_{(b_S,b_S)}^{\, T}\widehat B_{(b_S,b_S)}+2 B_{(b_S,b_L)}^{\, T} B_{(b_S,b_L)}\right.\\[2mm]
&\left.+\widehat B_{(b_U,b_U)}^{\, T}(\widehat \Delta_{1/c})_{(b_U,b_U)}\widehat B_{(b_U,b_U)}+\widehat B_{(b_L,b_L)}^{\, T}(\widehat \Delta_{1/c})_{(b_L,b_L)}\widehat B_{(b_L,b_L)}\right.\\[2mm]
&\left.+2 B_{(a,b_S)}^{\, T}(\Delta_{1/c})_{(a,b_S)} B_{(a,b_S)}+2 B_{(a,b_L)}^{\, T}(\Delta_{1/c})_{(a,b_L)} B_{(a,b_L)}\right.\\[2mm]
&\left. +2 B_{(a,c)}^{\, T}(\Delta_{1/c})_{(a,c)} B_{(a,c)}+2 B_{(c,b_U)}^{\, T}(\Delta_{1/c})_{(c,b_U)} B_{(c,b_U)}\right.\\[2mm]
&\left.+2 B_{(c,b_S)}^{\, T}(\Delta_{1/c})_{(c,b_S)} B_{(c,b_S)}+{\widehat C}^{\, T}
_{(\alpha,\alpha)}({\widehat
\Theta}_c)_{(\alpha,\alpha)}{\widehat C}_{(\alpha,\alpha)}\right. \\[2mm]
& \left.  +2C_{(\alpha,\beta)}^{\, T} C_{(\alpha,\beta)} + 2
C_{(\alpha,\gamma)}^{\, T}
(\Theta_c)_{(\alpha,\gamma)}C_{(\alpha,\gamma)} +{\widehat
C}_{(\beta,\beta)}^{\, T} ({\widehat
\Theta}_c)_{(\beta,\beta)}{\widehat C}_{(\beta,\beta)}
 \right)\, .
\end{array}
\end{equation}
\end{lemma}

Lemma \ref{acpres} shows that  ${\cal
A}_c(\overline Y,\overline{\mu},\overline{\Gamma},W_1, W_2)$ can be written as
\begin{equation}\label{acyy}
\begin{array}{l}
{\cal
A}_c(\overline Y,\overline{\mu},\overline{\Gamma},W_1, W_2) = \nabla^2_{xx}
L(\overline x,\overline Y,\overline{\mu},\overline{\Gamma})+cA^TD_cA\\[2mm]
\quad \, +2c B_{(a,b_S)}^{\, T}(\Delta_{1/c})_{(a,b_S)} B_{(a,b_S)}+2c B_{(a,b_L)}^{\, T}(\Delta_{1/c})_{(a,b_L)} B_{(a,b_L)}\\[2mm]
\quad \, +2c B_{(a,c)}^{\, T}(\Delta_{1/c})_{(a,c)} B_{(a,c)}+2c B_{(c,b_U)}^{\, T}(\Delta_{1/c})_{(c,b_U)} B_{(c,b_U)}\\[2mm]
 \quad \,+2c B_{(c,b_S)}^{\, T}(\Delta_{1/c})_{(c,b_S)} B_{(c,b_S)}
+2cC_{(\alpha,\gamma)}^{\, T}
(\Theta_c)_{(\alpha,\gamma)}C_{(\alpha,\gamma)}\, .
\end{array}
\end{equation}
For any $c^\prime, c>0$, let
\begin{equation}\label{eq:B-form}
\begin{array}{l}
{\cal B}_{c^\prime, c}(\overline Y,\overline{\mu},\overline{\Gamma},W_1, W_2) = \nabla^2_{xx}
L(\overline x,\overline Y,\overline{\mu},\overline{\Gamma})+c'A^TD_cA\\[2mm]
\quad \, +2c B_{(a,b_S)}^{\, T}(\Delta_{1/c})_{(a,b_S)} B_{(a,b_S)}+2c B_{(a,b_L)}^{\, T}(\Delta_{1/c})_{(a,b_L)} B_{(a,b_L)}\\[2mm]
\quad \, +2c B_{(a,c)}^{\, T}(\Delta_{1/c})_{(a,c)} B_{(a,c)}+2c B_{(c,b_U)}^{\, T}(\Delta_{1/c})_{(c,b_U)} B_{(c,b_U)}\\[2mm]
 \quad \,+2c B_{(c,b_S)}^{\, T}(\Delta_{1/c})_{(c,b_S)} B_{(c,b_S)}
+2cC_{(\alpha,\gamma)}^{\, T}
(\Theta_c)_{(\alpha,\gamma)}C_{(\alpha,\gamma)}\, .
\end{array}
\end{equation}

The following proposition shows that,  under Assumptions
(sdnop-A1) and (sdnop-A2), the basic Assumption B1  made in
Section \ref{general-discussions}
 is satisfied by nonlinear semidefinite nuclear norm composite optimization problem.

\begin{prop}\label{pronlsdp}
Suppose that  Assumptions (sdnop-A1)  and (sdnop-A2) are
satisfied. Then there exist two positive numbers $c_0 $ and
$\underline{\eta} $ such that for any $c \geq c_0$ and $W_1\in \partial_B [{\rm D} \theta_c]^*(F(\overline x)+\overline Y/c)$,
$W_2 \in
\partial_B \Pi_{{\cal S}^p_+}(\overline{\Gamma}-c g(\overline{x}))$,
 \[
\left \langle d,{\cal
A}_c(\overline Y,\overline{\mu},\overline{\Gamma},W_1, W_2)
d\right\rangle\ge \left  \langle d,{\cal B}_{c_0,
c}(\overline Y,\overline{\mu},\overline{\Gamma},W_1, W_2) d\right\rangle   \geq
\underline{\eta}\langle d, d\rangle, \quad \forall \, d\in \Re^n\, .
\]
\end{prop}

\noindent
 {\bf Proof.}  It follows from  Assumption (sdnop-A2) that there exists $\eta_0 >0$ such that
\begin{equation}\label{secnlp9-secnlp9}
\left \langle d, \nabla^2_{xx}L(\overline{x},\overline Y,
\overline{\mu},\overline{\Gamma})d \right \rangle-\psi^*(\overline Y)+
 2\left \langle \overline \Gamma, [{\rm D}g(\overline x)d]g(\overline x)^\dag [{\rm D}g(\overline x)d]\right
\rangle
\geq \eta_0 \|d\|^2
\end{equation}
for all $d \in {\rm app}(\overline Y,\overline{\mu},\overline{\Gamma}) \setminus \{0\}$.
By (\ref{equation3.5}), we obtain
\begin{equation}\label{equation3.20}
\begin{array}{ll}
\mbox{app }(\overline Y,\overline{\mu},\overline{\Gamma})
 &=\left\{d \in \Re^n: \begin{array}{l}
  {\cal J}  h(\overline{x}) d =0,
 B_{(b_U,b_S)}(Q)d=0,
 B_{(b_U,b_L)}(Q)d=0\\[4pt]
 {\widehat B}_{(b_S,b_S)}(Q)d=0,
 B_{(b_S,b_L)}(Q)d=0\\[4pt]
   {\widehat C}_{(\alpha,
\alpha)}(P)d=0 \, ,  \ C_{(\alpha, \beta)} (P)d = 0
\end{array}
\right\}\, .
\end{array}
\end{equation}
Since (\ref{secnlp9-secnlp9}) and (\ref{equation3.20}) hold, by
using Lemma \ref{lemma3.2nlp} with $\phi$ and ${\cal L}$ being
defined by
$$\phi (d):=\langle d, \nabla^2_{xx}L(\overline{x},
\overline Y,\overline{\mu},\overline{\Gamma})d \rangle-\psi^*(\overline Y)+
 2\left \langle \overline \Gamma, [{\rm D}g(\overline x)d]g(\overline x)^\dag [{\rm D}g(\overline x)d]\right
\rangle
$$
 and
$${\cal L}(d):=({\cal J}  h(\overline{x})d;B_{(b_U,b_S)}(Q)d;
 B_{(b_U,b_L)}(Q)d;
 {\widehat B}_{(b_S,b_S)}(Q)d;
 B_{(b_S,b_L)}(Q)d;{\widehat C}_{(\alpha,
\alpha)}(P)d;C_{(\alpha, \beta)}(P)d),$$
 for any   $d\in
\Re^n$, respectively, we know that
 there exist two positive numbers  $c_1$ and $\underline{\eta} \in (0, \eta_0/2]$ such
 that for any $c\ge c_1$,
\begin{equation}\label{eq:upgrade}
\begin{array}{l}
\left \langle d, \nabla^2_{xx}L(\overline{x},\overline Y,
\overline{\mu},\overline{\Gamma})d \right \rangle-\psi^*(\overline Y)+
 2\left \langle \overline \Gamma, [{\rm D}g(\overline x)d]g(\overline x)^\dag [{\rm D}g(\overline x)d]\right
\rangle
\\[2mm]
\quad +c\|B_{(b_U,b_S)}(Q)d\|^2+c\|B_{(b_U,b_L)}(Q)d\|^2+c\|{\widehat B}_{(b_S,b_S)}(Q)d\|^2+c\|B_{(b_S,b_L)}(Q)d\|^2
\\[2mm]
\quad +c\|{\cal J}  h(\overline{x}) d\|^2 +c \|{\widehat
C}_{(\alpha, \alpha)}(P)d\|^2 + c \|C_{(\alpha,
\beta)}(P)d\|^2
 \geq
2 \underline{\eta} \|d\|^2,\quad  \forall \, d \in \Re^n\, .
\end{array}
\end{equation}

Let $c_0 \ge c_1$ be  such that for any $c \geq c_0$,
\begin{equation}\label{equation3.21}
\begin{array}{l}
\displaystyle \max_{1 \leq l \leq n} \|{\cal J}_{x_l} F(\overline{x})\|^2
\sum_{i \in a, 1\leq j \leq |b_S|} \frac{c^{-1}(1-(w_{b_S})_j)^2}{\lambda_i(F(\overline x))(\lambda_i(F(\overline x))+c^{-1}(1-(w_{b_S})_j))} \leq \underline{\eta}/4,\\
\displaystyle \max_{1 \leq l \leq n} \|{\cal J}_{x_l} g(\overline{x})\|^2
\sum_{i \in \gamma, j \in \alpha} \frac{\lambda_j^2}{|\lambda_i|
(\lambda_j+c |\lambda_i|)} \leq \underline{\eta}/4.
\end{array}
\end{equation}
Let $c\ge c_0$ and $W_1\in \partial_B [{\rm D} \theta_c]^*(F(\overline x)+\overline Y/c)$,
$W_2 \in
\partial_B \Pi_{{\cal S}^p_+}(\overline{\Gamma}-c g(\overline{x}))$. Then there exist two
matrices $Q\in {\cal O}(F(\overline x))$ and $P\in {\cal O}(g(\overline x))$  and $\Delta_{1/c} \in {\cal S}^q$ satisfying (\ref{Shc})  and
$\Theta_c \in {\cal S}^p$ satisfying (\ref{Sh}) such that
\[
W_1(H_1)= Q\left( \Delta_{1/c} \circ (Q^TH_1Q)\right )Q^T,\quad \forall \, H_1
\in {\cal S}^q\, .
\]
and
\[
W_2(H_2) = P\left( \Theta_c \circ (P^TH_2P)\right )P^T,\quad \forall \, H_2
\in {\cal S}^p\, .
\]
It is easy to see from (\ref{equation3.21}) that   for any   $c
\ge c_0$ and $d\in \Re^n$ we have for $H_1={\rm D}F(\overline x)d$ and $\widehat H_1=Q^TH_1Q$ that
\[
\begin{array}{l}
-\psi^*(\overline Y)-2c \langle d, [ B_{(a,b_S)}^{\, T}(\Delta_{1/c})_{(a,b_S)} B_{(a,b_S)}]d \rangle\\[2mm]
=2\displaystyle \sum_{i \ne s}\displaystyle \frac{1}{\varpi_i} \langle I_{|b_s|}-Q^T_{b_s}YQ_{b_s}, Q_{b_s}^TH_1Q_{a_i}Q_{a_i}^THQ_{b_s} \rangle-2c \langle d, [ B_{(a,b_S)}^{\, T}(\Delta_{1/c})_{(a,b_S)} B_{(a,b_S)}]d \rangle\\[2mm]
=2 \displaystyle \sum_{i\ne s}\sum_{j=1}^{|b_S|}\frac{1-(w_{b_S})_j}{\varpi_i}\|Q_{a_i}^TH_1(Q_{b_S})_j\|^2-2c \displaystyle \sum_{i=1}^{s-1}\langle Q_{a_i}^TH_1Q_{b_S}\circ \Delta_{a_ib_S},Q_{a_i}^TH_1Q_{b_S} \rangle\\[2mm]
\leq 2\displaystyle \sum_{i=1}^{s-1}\sum_{j=1}^{|b_S|}\frac{1-(w_{b_S})_j}{\varpi_i}\|Q_{a_i}^TH_1(Q_{b_S})_j\|^2-2c \displaystyle \sum_{i=1}^{s-1}\langle Q_{a_i}^TH_1Q_{b_S}\circ \Delta_{a_ib_S},Q_{a_i}^TH_1Q_{b_S} \rangle\\[2mm]
=2\displaystyle \sum_{i=1}^{s-1}\sum_{j=1}^{|b_S|}\frac{1-(w_{b_S})_j}{\varpi_i}\|Q_{a_i}^TH_1(Q_{b_S})_j\|^2-2\displaystyle \sum_{i=1}^{s-1}\sum_{j=1}^{|b_S|}\frac{1-(w_{b_S})_j}{\varpi_i+c^{-1}(1-(w_{b_S})_j)}\|Q_{a_i}^TH_1(Q_{b_S})_j\|^2\\[2mm]
=2\displaystyle \sum_{i\in a}\sum_{j=1}^{|b_S|}\frac{c^{-1}(1-(w_{b_S})_j)^2}{\lambda_i(F(\overline x))(\lambda_i(F(\overline x)+c^{-1}(1-(w_{b_S})_j))}\|Q_i^T[{\rm D}F(\overline x)d](Q_{b_S})_j\|^2\end{array}
\]
\[
\begin{array}{l}
=2\displaystyle \sum_{i\in a}\sum_{j=1}^{|b_S|}\frac{c^{-1}(1-(w_{b_S})_j)^2}{\lambda_i(F(\overline x))(\lambda_i(F(\overline x)+c^{-1}(1-(w_{b_S})_j))}\displaystyle \sum_{l=1}^n[Q_i^T[{\cal J}_{x_l}F(\overline x)](Q_{b_S})_jd_l]^2\\[2mm]
\leq 2\displaystyle \sum_{i\in a}\sum_{j=1}^{|b_S|}\frac{c^{-1}(1-(w_{b_S})_j)^2}{\lambda_i(F(\overline x))(\lambda_i(F(\overline x)+c^{-1}(1-(w_{b_S})_j))}\displaystyle \sum_{l=1}^n\|Q_i\|^2\|[{\cal J}_{x_l}F(\overline x)]\|^2\|(Q_{b_S})_j\|^2d_l^2\\[2mm]
\leq \displaystyle \max_{1 \leq l \leq n} \|{\cal J}_{x_l} F(\overline{x})\|^2
\sum_{i \in a, 1\leq j \leq |b_S|} \frac{c^{-1}(1-(w_{b_S})_j)^2}{\lambda_i(F(\overline x))(\lambda_i(F(\overline x))+c^{-1}(1-(w_{b_S})_j))}\|d\|^2\\
\leq  \underline{\eta}
\|d\|^2/2.
\end{array}
\]
Similarly, we have from (\ref{equation3.21}) that   for any   $c
\ge c_0$ and $d\in \Re^n$  that
\[
\begin{array}[b]{l}
\quad 2\left \langle \overline \Gamma, [{\rm D}g(\overline x)d]g(\overline x)^\dag [{\rm D}g(\overline x)d]\right
\rangle -2c \left\langle d, C^{\, T}_{(\alpha,\gamma)}
(\Theta_c)_{(\alpha,\gamma)}C_{(\alpha,\gamma)} d \right\rangle
\\[2mm]
\displaystyle
  = 2 \left\langle \overline{\Gamma}, [{\rm D}g(\overline{x})d]
  g(\overline{x})^{\dagger}[{\rm D}g(\overline{x})d] \right\rangle
-2c \left\langle d, C^{\, T}_{(\alpha,\gamma)}
(\Theta_c)_{(\alpha,\gamma)}C_{(\alpha,\gamma)} d \right\rangle
\\[2mm]
\displaystyle
  = 2  \sum_{i \in \gamma, j \in
\alpha} \frac{\lambda_j}{ |\lambda_i| } \left(\sum_{l=1}^n P_i^T
{\cal J}_{x_l} g(\overline{x}) P_j d_l\right)^2 -2 c \sum_{i \in
\gamma, j \in \alpha} \frac{\lambda_j}{\lambda_j+c |\lambda_i| }
\sum_{l=1}^n \left(P_i^T {\cal J}_{x_l} g(\overline{x}) P_j
d_l\right)^2
 \\[2mm]\displaystyle
  \leq 2\sum_{i \in \gamma, j \in
\alpha}\left[\frac{\lambda_j^2}{|\lambda_i|(\lambda_j+c
|\lambda_i|)}\sum_{l=1}^n \|{\cal J}_{x_l} g(\overline{x})\|^2
\|P_i\|^2\|P_j\|^2d_l^2\right]
\\[2mm]\displaystyle
\leq 2 \max_{1 \leq l \leq n}\|{\cal J}_{x_l}
g(\overline{x})\|^2\sum_{i \in \gamma, j \in
\alpha}\frac{\lambda_j^2}{|\lambda_i|(\lambda_j+c
|\lambda_i|)} \|d\|^2\\[2mm] \displaystyle\leq  \underline{\eta}
\|d\|^2/2\, ,
\end{array}
\]
Therefore, we have from (\ref{eq:upgrade}),
 for any $c\ge c_0$, that
\begin{equation}\label{eq:eta-one-form}
\begin{array}{l}
\langle d, \nabla^2_{xx}L(\overline{x},\overline Y,
\overline{\mu},\overline{\Gamma})d \rangle+2c \langle d, [ B_{(a,b_S)}^{\, T}(\Delta_{1/c})_{(a,b_S)} B_{(a,b_S)}]d \rangle\\[2mm]
\quad +c_0\|B_{(b_U,b_S)}(Q)d\|^2+c_0\|B_{(b_U,b_L)}(Q)d\|^2+c_0\|{\widehat B}_{(b_S,b_S)}(Q)d\|^2\\[2mm]
+c_0\|B_{(b_S,b_L)}(Q)d\|^2
+2c \left\langle d,
C_{(\alpha,\gamma)}^{\, T}
(\Theta_c)_{(\alpha,\gamma)}C_{(\alpha,\gamma)} d \right\rangle
\\[2mm]
\quad +c_0\|{\cal J}  h(\overline{x}) d\|^2 +c_0 \|{\widehat
C}_{(\alpha, \alpha)}(P)d\|^2 + c_0 \|{C}_{(\alpha,
\beta)}(P)d\|^2
 \geq
\underline{\eta} \|d\|^2, \quad  \forall \, d \in \Re^n\, .
\end{array}
\end{equation}
In view of the expression $(\Delta_{1/c})_{ij}$ from (\ref{eq:dDeltac}) for $(i,j) \in (a \times b_L)\cup (a\times \{1,\ldots,|b_S|\})\cup (a \times c)\cup (b_U \times c)\cup (\{1,\ldots,|b_S|\} \times c)$, we obtain
\[
\begin{array}{rl}
B_{(a,b_L)}^{\, T}(\Delta_{1/c})_{(a,b_L)} B_{(a,b_L)}\succeq 0, &
 B_{(a,c)}^{\, T}(\Delta_{1/c})_{(a,c)} B_{(a,c)}\succeq 0,\\[2mm]
  B_{(c,b_U)}^{\, T}(\Delta_{1/c})_{(c,b_U)} B_{(c,b_U)}\succeq 0,&
 B_{(c,b_S)}^{\, T}(\Delta_{1/c})_{(c,b_S)} B_{(c,b_S)}\succeq 0,\\[2mm]
 \widehat B_{(b_U,b_U)}^{\, T}(\widehat \Delta_{1/c})_{(b_U,b_U)}\widehat B_{(b_U,b_U)} \succeq 0,& \widehat B_{(b_L,b_L)}^{\, T}(\widehat \Delta_{1/c})_{(b_L,b_L)}\widehat B_{(b_L,b_L)}\succeq 0.
 \end{array}
 \]
From this
 and the fact
 that ${\widehat
C}_{(\beta,\beta)}^{\, T} ({\widehat
\Theta}_c)_{(\beta,\beta)}{\widehat C}_{(\beta,\beta)} \succeq 0$,
 we can see that for any $c\ge c_0$,
\[
\left\langle d, {\cal B}_{c_0, c}(\overline Y,\overline{\mu},\overline{\Gamma},W_1, W_2)d \right\rangle \ge \underline{\eta} \|d\|^2,\quad  \forall
\, d \in \Re^n\, .
\]
By noting the fact that \[ {\cal A}_{c}(\overline Y,\overline{\mu},\overline{\Gamma},W_1, W_2) ={\cal B}_{c_0, c}(\overline Y,\overline{\mu},\overline{\Gamma},W_1, W_2) +(c-c_0)A^TD_cA\, ,
\]
we  complete the proof. \qed

Let Assumptions (sdnop-A1) and (sdnop-A2) be  satisfied.  Let the
two positive numbers $c_0$ and  $\underline{\eta}$ be defined as in
Proposition \ref{pronlsdp}. Let $c\ge c_0$. Then, by Propositions
\ref{prop:general-discussions} and  \ref{pronlsdp} and the fact that ${\rm D} \theta_c(\cdot)$
and $\Pi_{{\cal S}^p_+} (\cdot)$ are strongly semismooth everywhere,
there exist two positive numbers $\varepsilon >0$ and $\delta_0>0$
(both depending on $c$) and a locally Lipschitz continuous function
$x_c(\cdot, \cdot,\cdot)$ defined on
$\mathbb{B}_{\delta_0}(\overline Y,\overline{\mu},\overline{\Gamma})$  such that
for any $(Y,\mu, \Gamma)\in \mathbb{B}_{\delta_0}(\overline Y,\overline{\mu},\overline{\Gamma})$, $x_c(Y,\mu, \Gamma)$ is the unique minimizer of
$L_c(\cdot, Y,\mu, \Gamma)$ over
$\mathbb{B}_{\varepsilon}(\overline{x})$ and $x_c(\cdot,\cdot, \cdot)$ is
semismooth at $(Y,\mu, \Gamma)$.  Let $\vartheta_c: {\cal S}^q \times \Re^m \times {\cal
S}^p\mapsto \Re$ be defined as
 (\ref{dual}), i.e.,
\[
  \vartheta_c(Y,\mu,\Gamma):=\min_{x \in
\mathbb{B}_{\varepsilon}(\overline{x})} L_c(x,Y,\mu,\Gamma), \quad
(Y, \mu,\Gamma)\in {\cal S}^q \times  \Re^m \times {\cal S}^p\,  .
\]
Then it holds that
\[
  \vartheta_c(Y,\mu,\Gamma)= L_c (x_c(Y,\mu, \Gamma),Y, \mu, \Gamma)\, , \quad (Y,\mu, \Gamma)
\in \mathbb{B}_{\delta_0}(\overline Y,\overline{\mu},\overline{\Gamma})\, .
\]
Furthermore, it follows  from Propositions \ref{comdcproperty1}
and  \ref{pronlsdp} that the concave function
 $\vartheta_c(\cdot, \cdot,\cdot)$ is continuously differentiable
 on $\mathbb{B}_{\delta_0}(\overline Y,\overline{\mu},\overline{\Gamma})$ with
\[
{\rm D}\vartheta_c (Y,\mu, \Gamma)^*=\left (
\begin{array}{c}
-c^{-1}Y+c^{-1}{\rm D}\theta_c(F(x_c(Y,\mu, \Gamma))+Y/c))^*\\[2mm]
h(x_c(Y,\mu, \Gamma))\\[2mm]
-c^{-1}\Gamma+c^{-1}\Pi_{{\cal S}^p_+}(\Gamma-c g(x_c(Y,\mu, \Gamma)))
\end{array}
\right ), \,y =(Y,\mu,\Gamma) \in \mathbb{B}_{\delta_0}
(\overline Y,\overline{\mu},\overline{\Gamma})\, .
\]
For any $(\Delta Y,\Delta \mu, \Delta \Gamma) \in {\cal S}^q \times \Re^m \times {\cal S}^p$,
let $\overline{{\cal V}}_c (\Delta Y,\Delta \mu, \Delta \Gamma)$ be defined
as in (\ref{comcalV}). By Propositions \ref{comthta} and
\ref{pronlsdp},
 we have for any $(\Delta Y,\Delta \mu, \Delta \Gamma) \in {\cal S}^q \times \Re^m \times {\cal S}^p$ that
\[
\partial_B(\nabla \vartheta_c)(\overline Y,\overline{\mu}, \overline{\Gamma})(\Delta Y,\Delta \mu, \Delta \Gamma)
 \subseteq \overline{{\cal V}}_c (\Delta Y,\Delta \mu, \Delta \Gamma)\, .
\]
Since when $c \rightarrow \infty$,
\[
c[{\Delta_{1/c}}]_{ij}=\left
 \{
\begin{array}{clll}
\quad \quad 0 && & (i,j)\in (a\times a \cup b_U),\\[12pt]
\quad \quad 0 && & (i,j)\in  (c \times b_L\cup c),\\[12pt]
\displaystyle \frac{1-(w_{b_S})_j}{\lambda_i(\overline X)+c^{-1}(1-(w_{b_S})_j)}
& \rightarrow& \displaystyle \frac{1-(w_{b_S})_j}{\lambda_i(\overline X)}
  & (i,j)\in a\times \{1,\ldots,|b_S|\},\\[12pt]
\displaystyle \frac{2}{\lambda_i(\overline X)+2c^{-1}}& \rightarrow&
\displaystyle \frac{2}{\lambda_i(\overline X)}
 & (i,j)\in a\times b_L,\\[12pt]
\displaystyle \frac{2}{\lambda_i(\overline X)-\lambda_j(\overline X)+2c^{-1}} & \rightarrow&
\displaystyle \frac{2}{\lambda_i(\overline X)-\lambda_j(\overline X)}
& (i,j)\in a\times c,\\[12pt]
\displaystyle \frac{2}{-\lambda_j(\overline X)+2c^{-1}} & \rightarrow&
\displaystyle \frac{2}{-\lambda_j(\overline X)} & (i,j)\in b_U\times c,\\[12pt]
\displaystyle \frac{(w_{b_S})_i+1}{c^{-1}((w_{b_S})_i+1)-\lambda_j(\overline X)} & \rightarrow&
\displaystyle \frac{(w_{b_S})_i+1}{-\lambda_j(\overline X)}
& (i,j)\in \{1,\ldots,|b_S|\}\times c.
\end{array}
\right.
\]
where $\overline X=F(\overline x)$, and
\[ \lim_{c\to \infty} c (\Theta_c)_{ij} =\lim_{c\to \infty}
c
\frac{\lambda_i}{\lambda_i+c|\lambda_j|}=\frac{\lambda_i}{|\lambda_j|},
\forall \, (i, j)\in \alpha \times \gamma\, ,
\]
 we know that there exists a positive number $\overline{\eta}$ such that
\begin{equation}\label{mu1nlsdp}
\left\langle d, {\cal B}_{c_0, c}(\overline Y,\overline{\mu},\overline{\Gamma},W_1, W_2)d \right\rangle \le \overline{\eta} \langle d, d\rangle
\quad \begin{array}{l}
\forall\, d\in \Re^n,  c\ge c_0 \,{\rm and}\,\\
W_1\in \partial_B [{\rm D} \theta_c]^*(F(\overline x)+\overline Y/c),\\
W_2 \in
\partial_B \Pi_{S^p_+}(\overline{\Gamma}-c g(\overline{x})).
\end{array}
\end{equation}


Let $c\ge c_0$, $W_1\in \partial_B [{\rm D} \theta_c]^*(F(\overline x)+\overline Y/c)$ and $W_2 \in
\partial_B \Pi_{S^p_+}(\overline{\Gamma}-c g(\overline{x}))$. Then there exist two
matrices $Q\in {\cal O}(F(\overline x)$ with $ P\in {\cal O}(g(\overline x))$
and $\Delta_{1/c}$ satisfying (\ref{Shc}) such that (\ref{eq:w1h}) holds,
$\Theta_c \in {\cal S}^p$ satisfying (\ref{Sh}) such that
(\ref{wh}) holds. Let $A(Q,P)$ have the singular value decomposition
as in (\ref{singvalu}), i.e.,
\begin{equation}\label{eq:singvalu-rep}
A(Q,P)=U[\Sigma (Q,P)\,\,\,\,\, 0]R^T\, .
\end{equation}

  Let $\overline{y}:=(\overline Y, \overline{\mu}, \overline{\Gamma})$. Then  we
have the  following result for ${\cal A}_c(\overline{y}, W_1,W_2)$.

\begin{lemma}\label{le2} Let $c> c_0$  and $W_1\in \partial_B [{\rm D} \theta_c]^*(F(\overline x)+\overline Y/c)$ and $W_2 \in
\partial_B \Pi_{{\cal S}^p_+}(\overline{\Gamma}-c g(\overline{x}))$. Suppose that Assumptions
(sdnop-A1) and (sdnop-A2) are satisfied.   Then we have
\begin{equation}\label{ac-2-101}
\begin{array}{l}
{\cal A}_c(\overline{y},W_1,W_2)^{-1} \preceq R\left [
\begin{array}{cc}
\Sigma^{-1}U^T\Big({\underline{\sigma}}\underline{\eta}I_{n_2}+(c-c_0)D_c\Big)^{-1}U \Sigma^{-1} & 0\\
 0 & {\underline{\sigma}}^{-1}\underline{\eta}^{-1} I_{n_3}
\end{array}\right
]R^T,
\end{array}
\end{equation}
\begin{equation}\label{ac-1-101}
\begin{array}{l}
{\cal A}_c(\overline{y},W_1,W_2)^{-1} \succeq R\left [
\begin{array}{cc}
\Sigma^{-1}U^T \Big(\overline{\sigma}\overline{\eta}I_{n_2}+(c-c_0)D_c\Big)^{-1}U\Sigma^{-1} & 0\\
 0 & {\overline{\sigma}}^{-1}\overline{\eta}^{-1} I_{n_3}
\end{array}\right]R^T\, ,
\end{array}
\end{equation}
and
\begin{equation}\label{ac-1-101a}
\|{\cal A}_c(\overline{y},W_1,W_2)^{-1}A^TD_c u\| \leq \sqrt{2}\left(
\overline{\sigma} + ( \underline{\sigma} \underline{\eta} )^{-2}
(\overline{\sigma} \overline{\eta} )^2\right) \|u\|/(c-c_0),\ \,
\forall\, u \in \Re^{n_2},
\end{equation}
 where $\Sigma: =\Sigma(Q,P)$.
\end{lemma}

\noindent {\bf Proof.}  Let
 $\hat{c}:=c-c_0$.
By (\ref{acyy}), (\ref{eq:B-form}), and the singular value
decomposition (\ref{eq:singvalu-rep}) of $A:=A(P)$, we have
\begin{equation}\label{calnlp}
\begin{array}[b]{l}
\quad {\cal A}_c(\overline{y},W_1,W_2)^{-1}=\Big({\cal B}_{c_0,
c}(\overline{y},W_1,W_2)+\hat{c}A^TD_cA \Big)^{-1}
\\[2mm]
=\left({\cal B}_{c_0,c}(\overline{y},W_1,W_2)+\hat{c}R[\Sigma \quad
0]^TU^TD_cU[\Sigma \quad 0]R^T\right)^{-1}
\\[2mm]
=R \left(R^T{\cal B}_{c_0,c}(\overline{y},W_1,W_2)R+\hat{c}\left [
\begin{array}{cc}
\Sigma & 0\\
0 & I_{n_3}
\end{array}\right
]\left [
\begin{array}{cc}
U^T D_cU  & 0
\\
0 & 0
\end{array}\right
]\left [
\begin{array}{cc}
\Sigma & 0\\
0 & I_{n_3}
\end{array}\right
]\right)^{-1}R^T
\\[2mm]
=R \left [
\begin{array}{cc}
\Sigma^{-1} & 0
\\
0 & I_{n_3}
\end{array}\right
] \left({\cal G}_{c_0,c}(\overline{y},W_1,W_2)+\hat{c} \left [
\begin{array}{cc}
U^TD_cU  & 0\\
0 & 0
\end{array}\right
]\right)^{-1} \left [
\begin{array}{cc}
\Sigma^{-1} & 0\\
0 & I_{n_3}
\end{array}\right
]R^T,
\end{array}
\end{equation}
where
\[
{\cal G}_{c_0,c}(\overline{y},W_1,W_2):=\left [
\begin{array}{cc}
\Sigma^{-1} & 0\\
0 & I_{n_3}
\end{array}\right
]R^T{\cal B}_{c_0,c}(\overline{y},W_1,W_2)R\left [
\begin{array}{cc}
\Sigma^{-1} & 0\\
0 & I_{n_3}
\end{array}\right
]\, .
\]
It follows from  Proposition \ref{pronlsdp}, the definitions of
${\underline{\sigma}}$ and $\overline{\sigma}$, and
(\ref{mu1nlsdp}) that
\begin{equation}\label{eq:G-lower-bound}
{\cal G}_{c_0,c}(\overline{y},W_1,W_2) \succeq  \underline{ \eta} \left [
\begin{array}{cc}
\Sigma^{-1} & 0\\
0 & I_{n_3}
\end{array}\right
]^2 \succeq {\underline{\sigma}}\underline{\eta} I_{n}
\end{equation}
and
\begin{equation}\label{eq:G-upper-bound}
{\cal G}_{c_0,c}(\overline{y},W_1,W_2) \preceq   \overline{\eta} \left [
\begin{array}{cc}
\Sigma^{-1} & 0\\
0 & I_{n_3}
\end{array}\right
]^2 \preceq \overline{\sigma} \overline{\eta} I_{n}\, .
\end{equation}
Therefore, (\ref{ac-2-101}) and (\ref{ac-1-101}) follow from
(\ref{calnlp}).

\vskip 15 true pt
 Now we turn to the proof  of (\ref{ac-1-101a}).
Let \[ \overline{{\cal G}}_{c_0,c}(\overline{y},W_1,W_2):= \left [
\begin{array}{cc}
U & 0\\
0 & I_{n_3} \end{array}\right ] {\cal G}_{c_0,c}(\overline{y},W_1,W_2)
\left [
\begin{array}{cc}
U^T & 0\\
0 & I_{n_3} \end{array}\right ]
\]
and
\[
 {\overline {\cal
H}}_{c_0}(\overline{y},W_1,W_2) :=  \overline{{\cal
G}}_{c_0}(\overline{y},W_1,W_2)^{-1}\, .
\]
 Partition $ {\overline{\cal H}}_{c_0,c}(\overline{y},W_1,W_2)$ as
$$ \overline{{\cal H}}_{c_0,c}(\overline{y},W_1,W_2)=\left [
\begin{array}{cc}
H_1(W_1,W_2) & H_{2}(W_1,W_2)^T\\[2mm]
H_{2}(W_1,W_2) & H_{3}(W_1,W_2)
\end{array}
\right ] \]
 with  $H_{1}(W_1,W_2)\in {\cal S}^{n_2}$, $H_{2}(W_1,W_2)\in
\Re^{n_3 \times n_2}$, and $H_{3}(W_1,W_2)\in {\cal S}^{n_3}.$
 Then, it follows from (\ref{eq:G-lower-bound}) and
 (\ref{eq:G-upper-bound}) that
\begin{equation}\label{eq:H12-form}
\begin{array}{l}
\|H_{1}(W_1,W_2)\|_2 \le (\underline{\sigma}\underline{\eta})^{-1}\, , \
\ \|H_{1}(W_1,W_2)^{-1}\|_2 \le \overline{\sigma}\overline{\eta}\\[2mm]
 {\rm and} \ \ \|H_2(W_1,W_2)H_{1}(W_1,W_2)^{-1}\|_2 \leq (\underline{\sigma}
\underline{\eta}  )^{-1} \overline{\sigma}\overline{\eta}\, .
\end{array}
\end{equation}
For any $\varepsilon>0$, let
$$
D_{c,\varepsilon}:=D_c+\varepsilon I_{n_2}, \quad {\cal
A}_{c,\varepsilon}(\overline{y},W_1,W_2):={\cal
B}_{c_0,c}(\overline{y},W_1,W_2)+\hat{c}A^TD_{c,\varepsilon}A.$$
Let
$\varepsilon >0$. By referring to (\ref{calnlp}), we obtain
\[
\begin{array}{l}
 {\cal A}_{c, \varepsilon}(\overline{y},W_1,W_2)^{-1}
  \\[2mm]
  =R
  \left [
\begin{array}{cc}
\Sigma^{-1}U^T & 0
\\
0 & I_{n_3}
\end{array}\right]
\left(\overline{{\cal G}}_{c_0,c}(\overline{y},W_1,W_2)+\hat{c} \left [
\begin{array}{cc}
D_{c, \varepsilon}  & 0\\
0 & 0
\end{array}\right
]\right)^{-1} \left [
\begin{array}{cc}
U\Sigma^{-1} & 0\\
0 & I_{n_3} \end{array}\right ]R^T\, ,
\end{array}
\]
which, together with (\ref{eq:singvalu-rep})  and the
Sherman-Morrison-Woodbury formula (cf. \cite[Section
2.1]{GVanLoan96}), implies
$$
\begin{array}{l}
 {\cal A}_{c,\varepsilon}(\overline{y},W_1,W_2)^{-1}A^TD_{c,\varepsilon}\\[2mm]
 =
 R \left [
\begin{array}{ll}
\Sigma^{-1}U^T &0\\[2mm]
0& I_{n_3}
\end{array}
\right ] \left [
\begin{array}{c}
\left ( {H}_{1}(W_1,W_2)^{-1}+\hat{c}D_{c,\varepsilon} \right)^{-1}D_{c,\varepsilon}\\[2mm]
{H}_{2}(W_1,W_2) {H}_{1}(W_1,W_2)^{-1}\left
({H}_{1}(W_1,W_2)^{-1}+\hat{c}D_{c,\varepsilon}\right )^{-1}
D_{c,\varepsilon}\end{array} \right ]\, .
\end{array}
$$
Since, it follows from the Sherman-Morrison-Woodbury formula  that
$$
\begin{array}{l}
\left({H}_{1}(W_1,W_2)^{-1}+\hat{c}D_{c,\varepsilon}\right)^{-1}D_{c,\varepsilon}\\[1.2mm]
=\left(\hat{c}I_{n_2}+ D_{c,\varepsilon}^{-1}
{H}_{1}(W_1,W_2)^{-1}\right)^{-1}\\[1.2mm]
=\hat{c}^{-1}I_{n_2}-\hat{c}^{-2} D_{c,\varepsilon}^{-1} \left
(I_{n_2}+\hat{c}^{-1} {H}_{1}(W_1,W_2)^{-1} D_{c,\varepsilon}^{-1}
\right)^{-1}{H}_{1}(W_1,W_2)^{-1}\\[1.2mm]
=\hat{c}^{-1}I_{n_2}-\hat{c}^{-1}\left(\hat{c}D_{c,\varepsilon}+
{H}_{1}(W_1,W_2)^{-1}\right)^{-1}{H}_{1}(W_1,W_2)^{-1},
\end{array}
$$
we have
$$
\begin{array}{l}
\quad {\cal A}_c(\overline{y},W_1,W_2)^{-1}A^TD_c =\displaystyle
\lim_{\varepsilon \downarrow 0} {\cal
A}_{c,\epsilon}(\overline{y},W_1,W_2)^{-1}A^TD_{c,\epsilon}\\[1.2mm]
= R \left [
\begin{array}{c}
\Sigma^{-1}U^T \\
 {H}_{2}(W_1,W_2) {H}_{1}(W_1,W_2)^{-1}
\end{array}
\right ]\left(
\hat{c}^{-1}I_{n_2}-\hat{c}^{-1}\left(\hat{c}D_c+{H}_{1}(W_1,W_2)^{-1}\right)^{-1}{H}_{1}(W_1,W_2)^{-1}
\right )\, .
\end{array}
$$
Therefore, from the definition of $\overline{\sigma}$ and
(\ref{eq:H12-form}) we have for any $u \in \Re^{n_2}$ that
$$
\begin{array}{l}
\quad \|{\cal A}_c(\overline{y},W_1,W_2)^{-1}A^TD_cu\|^2\\[1.5mm]
\leq \left( \overline{\sigma} + ( \underline{\sigma}
\underline{\eta} )^{-2} (\overline{\sigma} \overline{\eta} )^2
 \right)
\left\|\left(\hat{c}^{-1}I_{n_2}-\hat{c}^{-1}\left(\hat{c}D_c+{H}_{1}(W_1,W_2)^{-1}\right)^{-1}{H}_{1}(W_1,W_2)^{-1}\right)u\right\|^2\\[1.5mm]
\leq  \left( \overline{\sigma} + ( \underline{\sigma}
\underline{\eta} )^{-2} (\overline{\sigma} \overline{\eta} )^2
 \right) \left(
\hat{c}^{-1}\|u\|+\hat{c}^{-1}\left\| \left(
\hat{c}D_c+{H}_{1}(W_1,W_2)^{-1} \right )^{-1}\right \|_2
\|{H}_{1}(W_1,W_2)^{-1}\|_2  \|u\|\right)^2
\\[1.5mm]
\leq \left( \overline{\sigma} + ( \underline{\sigma}
\underline{\eta} )^{-2} (\overline{\sigma} \overline{\eta} )^2
 \right)\hat{c}^{-2}\left(1+\|{H}_{1}(W_1,W_2)\|_2
 \|{H}_{1}(W_1,W_2)^{-1}\|_2\right)^2\|u\|^2\\[1.5mm]
\leq \left( \overline{\sigma} + ( \underline{\sigma}
\underline{\eta} )^{-2} (\overline{\sigma} \overline{\eta} )^2
 \right) \hat{c}^{-2} \left(1+ ( \underline{\sigma}
\underline{\eta} )^{-1} (\overline{\sigma} \overline{\eta} )
\right)^2\|u\|^2\, ,
\end{array}
$$
which, together with the fact that $\overline{\sigma}\ge 1$, proves
(\ref{ac-1-101a}).
 \qed

\vskip 7 true pt
Let
\begin{equation}\label{eq:barc-sdp}
\overline{c}:=\max\left \{(2+\sqrt{2})c_0, \, (\overline{\sigma}
\overline{\eta}-c_0)^2/c_0, \ (\underline{\sigma}
\underline{\eta}/2-c_0)^2/c_0\right\}  \,
\end{equation}
and

\begin{equation}\label{varrho0-sdp}
 \varrho_0:= \left( \overline{\nu} \overline{\sigma}{\underline{\sigma}}^{-2}
\underline{\eta}^{-2}\max \left
\{8\overline \nu_1^2,16\overline \nu_2^2,32\overline \nu_3^2,64\overline \nu_4^2,128\overline \nu_5^2,128\overline{\nu}_0^2, \ 4
\kappa_0^2\right\}\right)^{1/2}
 \, .
\end{equation}
where ($\overline X=F(\overline x)$)
\[
\begin{array}{l}
\overline \nu_1=\displaystyle \left (\max_{i \in a, j \in \{1,\ldots,|b_S|\}}\displaystyle \frac{1-(w_{b_S})_j}{\lambda_i(\overline X)}\right)\\[2mm]
\overline \nu_2=\displaystyle \left (\max_{i \in a, j \in b_L}\displaystyle \frac{2}{\lambda_i(\overline X)}\right)\\[2mm]
\overline \nu_3=\displaystyle \left (\max_{i \in a, j \in c}\displaystyle \frac{2}{\lambda_i(\overline X)-\lambda_j(\overline X)}\right)\\[2mm]
\overline \nu_4=\displaystyle \left (\max_{i \in b_U, j \in c}\displaystyle \frac{2}{-\lambda_j(\overline X)}\right)\\[2mm]
\overline \nu_5=\displaystyle \left (\max_{i \in \{1,\ldots,|b_S|\},j \in c}\displaystyle \frac{(w_{b_S})_i+1}{-\lambda_j(\overline X)}\right)
\end{array}
\]
and
\[ \kappa_0:=\sqrt{2}\left( \overline{\sigma} + (
\underline{\sigma} \underline{\eta} )^{-2} (\overline{\sigma}
\overline{\eta} )^2\right)\, .
\]

\begin{prop}\label{thdcnlsdp} Suppose that Assumptions (sdnop-A1) and
(sdnop-A2) are satisfied.  
  Then there exists a positive number $\mu_0$ such that for any $c \geq {\overline{c}}$ and $\Delta y
\in {\cal S}^q \times \Re^m \times {\cal S}^p$,
 \begin{equation}\label{eq:direction-added-sdp}
\|(x_c)^\prime(\overline{y};\Delta y)\|\le \mu_0 \|\Delta
y\|/c
\end{equation}
 and
 \begin{equation}\label{import11}
\left \langle V(\Delta y)+c^{-1}\Delta y,  \Delta y \right \rangle
\in \mu_0[-1, 1] \| \Delta y\|^2/c^{2},\quad \forall \, V(\Delta y)
\in \overline{{\cal V}}_c(\Delta y)\, .
\end{equation}
\end{prop}

\noindent {\bf Proof.} Let $c\ge {\overline{c}}$. Let $\Delta y:=
(\Delta Y,\Delta\mu, \Delta\Gamma) \in {\cal S}^q \times  \Re^m \times {\cal S}^p$. From the
proof of Proposition  \ref{comthta} we know that there exist
$W_1\in \partial_B [{\rm D} \theta_c]^*(F(\overline x)+\overline Y/c)$  and
$W_2 \in
\partial_B \Pi_{S^p_+}(\overline{\Gamma}-c g(\overline{x}))$ such that
\begin{equation}\label{comxderrii-sdp}
(x_c)'(\overline y;\Delta y)={\cal A}_c(\overline y,W_1,W_2 )^{-1}\left(-{\rm D}F(\overline x)^* W_1(\Delta Y/c)-{\cal J}
h(\overline x)^T (\Delta\mu)+{\rm D}g(\overline x)^*
W_2(\Delta\Gamma)\right).
\end{equation}
For this $W_1\in \partial_B [{\rm D} \theta_c]^*(F(\overline x)+\overline Y/c)$, there exist
matrices $Q\in {\cal O}(F(\overline x))$ and $\Delta_{1/c} \in {\cal S}^q$ satisfying (\ref{eq:dDeltac}) such that
\[
W_1(H_1)= Q\left( \Delta_{1/c} \circ (Q^TH_1Q)\right )Q^T,\quad \forall \, H_1
\in {\cal S}^q\, .
\]
For this $W_2 \in
\partial_B \Pi_{{\cal S}^p_+}(\overline{\Gamma}-cg(\overline{x}))$,
 there exist two matrices $P \in {\cal O}(\overline X)$ and $\Theta_c\in {\cal S}^p$ satisfying
  (\ref{Sh}) such that
\[
W_2(H_2) = P\left( \Theta_c \circ (P^TH_2P)\right )P^T,\quad \forall \, H_2
\in {\cal S}^p\, .
\]
Let $A: =A(Q,P)$ have the singular value decomposition as in
(\ref{singvalu}), i.e.,
\begin{equation}\label{singvalue-in}
 A=U[\Sigma \,\,\,\,\, 0]R^T\, ,
\end{equation}
where $\Sigma:=\Sigma(Q,P)$.\\
For any two
index sets $\chi, \chi' \in \{b_U,b_S,b_L\}$, let
\[
\xi_{(\chi,\chi')}:={\rm vec}(Q_{\chi}^{\, T}\Delta Y
Q_{\chi'})\, , \quad {\widehat \xi}_{(\chi,\chi)}:={\rm
svec}(Q_{\chi}^{\, T}\Delta Y Q_{\chi})\, .
\]

 For any two
index sets $\chi, \chi' \in \{\alpha, \beta, \gamma\}$, let
\[
\omega_{(\chi,\chi')}:={\rm vec}(P_{\chi}^{\, T}\Delta\Gamma
P_{\chi'})\, , \quad {\widehat \omega}_{(\chi,\chi)}:={\rm
svec}(P_{\chi}^{\, T}\Delta\Gamma P_{\chi})\, .
\]
Define
\[
  \Delta d_0:= \left ( \begin{array}{l}
 \Delta\mu\\
 {\widehat \xi}_{(b_U,b_U)}\\
 \xi_{(b_U,b_S)}\\
 \xi_{(b_U,b_L)}\\
 {\widehat \xi}_{(b_S,b_S)}\\
 \xi_{(b_S,b_L)}\\
 {\widehat \xi}_{(b_L,b_L)}\\
  {\widehat \omega} _{(\alpha,\alpha)}\\
{\widehat \omega}_{(\beta,\beta)}\\
 \omega _{(\alpha,\beta)}
\end{array}\right)\, , \quad
\Delta d: =\left ( \begin{array}{l} \Delta d_0\\
\xi_{(a,b_S)}\\
  \xi_{(a,b_L)}\\
   \xi_{(a,c)}\\
   \xi_{(c,b_U)}\\
    \xi_{(c,b_S)}\\
\omega_{(\alpha,\gamma)}
\end{array}
\right )\, .
\]
Then,  {}from (\ref{comxderrii-sdp}), we have
\[
\begin{array}{ll}
(x_c)^\prime(\overline{y};\Delta y) =&
-{\cal
A}_c(\overline{y},W_1,W_2)^{-1}
[A^T D_c\Delta d_0+2 B_{(a,b_S)}^{\, T}(\Delta_{1/c})_{(a,b_S)}\xi_{(a,b_S)}\\
&+2 B_{(a,b_L)}^{\, T}(\Delta_{1/c})_{(a,b_L)}\xi_{(a,b_L)}
 +2 B_{(a,c)}^{\, T}(\Delta_{1/c})_{(a,c)}\xi_{(a,c)}\\
 &+2 B_{(c,b_U)}^{\, T}(\Delta_{1/c})_{(c,b_U)} \xi_{(c,b_U)}
 + 2 B_{(c,b_S)}^{\, T}(\Delta_{1/c})_{(c,b_S)}\xi_{(c,b_S)}\\
 & -2C_{(\alpha,\gamma)}^{\,
T}(\Theta_c)_{(\alpha,\gamma)}\omega_{(\alpha,\gamma)}]
\end{array}
\]
and
\begin{equation}\label{Ahh-0}
\begin{array}[b]{l}
\quad  \left \langle (x_c)^\prime(\overline{y};\Delta y),
(x_c)^\prime(\overline{y};\Delta y)\right\rangle
\\[2mm]
\le2 \left \langle A^TD_c\Delta d_0, {\cal
A}_c(\overline{y},W_1,W_2)^{-2}A^TD_c \Delta d_0\right \rangle
\\[2mm]
+16 \left \langle B_{(a,b_S)}^{\, T}(\Delta_{1/c})_{(a,b_S)}\xi_{(a,b_S)}, {\cal
A}_c(\overline{y},W_1,W_2)^{-2} B_{(a,b_S)}^{\, T}(\Delta_{1/c})_{(a,b_S)}\xi_{(a,b_S)} \right \rangle\\[2mm]
+32\left \langle B_{(a,b_L)}^{\, T}(\Delta_{1/c})_{(a,b_L)}\xi_{(a,b_L)}, {\cal
A}_c(\overline{y},W_1,W_2)^{-2} B_{(a,b_L)}^{\, T}(\Delta_{1/c})_{(a,b_L)}\xi_{(a,b_L)} \right \rangle\\[2mm]
+64 \left \langle B_{(a,c)}^{\, T}(\Delta_{1/c})_{(a,c)}\xi_{(a,c)}, {\cal
A}_c(\overline{y},W_1,W_2)^{-2} B_{(a,c)}^{\, T}(\Delta_{1/c})_{(a,c)}\xi_{(a,c)} \right \rangle\\[2mm]
+128\left \langle B_{(c,b_U)}^{\, T}(\Delta_{1/c})_{(c,b_U)}\xi_{(c,b_U)}, {\cal
A}_c(\overline{y},W_1,W_2)^{-2} B_{(c,b_U)}^{\, T}(\Delta_{1/c})_{(c,b_U)}\xi_{(c,b_U)} \right \rangle\\[2mm]
+256 \left \langle B_{(c,b_S)}^{\, T}(\Delta_{1/c})_{(c,b_S)}\xi_{(c,b_S)}, {\cal
A}_c(\overline{y},W_1,W_2)^{-2} B_{(c,b_S)}^{\, T}(\Delta_{1/c})_{(c,b_S)}\xi_{(c,b_S)} \right \rangle\\[2mm]
 +256 \left \langle C^{\, T}_{(\alpha,\gamma)}(\Theta_c)_{(\alpha,\gamma)}\omega_{(\alpha,\gamma)}, {\cal
A}_c(\overline{y},W_1,W_2)^{-2}C^{\,
T}_{(\alpha,\gamma)}(\Theta_c)_{(\alpha,\gamma)}\omega_{(\alpha,\gamma)}\right
\rangle\, .
\end{array}
\end{equation}
{}From (\ref{ac-1-101a}), we have for $c \geq \overline{c} \quad
(\geq (2+\sqrt{2})c_0)$ that
\begin{equation}\label{anup-1-add0}
\begin{array}{l}
\quad \left\langle A^TD_c \Delta d_0, {\cal
A}_c(\overline{y},W_1,W_2)^{-2}A^TD_c \Delta d_0 \right\rangle
 \\[2mm]
=\|{\cal A}_c(\overline{y},W_1,W_2)^{-1}A^TD_c \Delta d_0\|^2
\\[2mm]
\le \kappa_0^2{\hat c}^{-2} \left(\|\Delta d_0 \right\|)^2
\\[2mm]
\le \kappa_0^2{\hat c}^{-2} (\|(\Delta\mu, \xi_{(b_U,b_U)},
 \xi_{(b_S,b_S)},
 \xi_{(b_L,b_L)}, { \omega}_{(\alpha,\alpha)}, {
\omega}_{(\beta,\beta)})\|^2\\[2mm]
\quad + 2  \|( \xi_{(b_U,b_S)}, \xi_{(b_S,b_L)},
 \xi_{(b_U,b_L)}, \omega_{(\alpha,\beta)})\|^2
)
\\[2mm]
\le \displaystyle \frac{1}{2}\varrho_0^2
c^{-2}(\|(\Delta\mu, \xi_{(b_U,b_U)},
 \xi_{(b_S,b_S)},
 \xi_{(b_L,b_L)}, { \omega}_{(\alpha,\alpha)}, {
\omega}_{(\beta,\beta)})\|^2\\[2mm]
\quad  + 2\|( \xi_{(b_U,b_S)}, \xi_{(b_S,b_L)},
 \xi_{(b_U,b_L)}, \omega_{(\alpha,\beta)})\|^2 )\, .
\end{array}
\end{equation}
Let
\[
\underline{{\cal E}}_c:= \left(\overline{\sigma} \overline{\eta}
I_{n_2} +(c-c_0)D_c\right )^{-1}\, , \quad \overline{{\cal
E}}_c:=\left(\underline{\sigma}\underline{\eta}I_{n_2}
+(c-c_0)D_c\right)^{-1}
\]
and
\begin{equation}\label{eu1-A}
 \underline{{\cal H}}_c := \left[
\begin{array}{cc}
 \underline{{\cal E}}_c & 0
\\
0 &{\overline{\sigma}}^{-1} \overline{\eta}^{-1}I_{n_3}
\end{array}\right]\, , \quad
\overline{{\cal H}}_c := \left[
\begin{array}{cc}
\overline{{\cal E}}_c & 0
\\
0 &{\underline{\sigma}}^{-1} \underline{\eta}^{-1}I_{n_3}
\end{array}\right]\, .
\end{equation}

 We know  from Lemma \ref{le2}, (\ref{eu1-A}),
(\ref{eq:nu-lower-upper}), and (\ref{eq:two-nus}) that
\[
\begin{array}[b]{l}
\quad \left\langle C^{\,
T}_{(\alpha,\gamma)}(\Theta_c)_{(\alpha,\gamma)}\omega_{(\alpha,\gamma)},
{\cal A}_c(\overline{y},W_1,W_2)^{-2}C^{\,
T}_{(\alpha,\gamma)}(\Theta_c)_{(\alpha,\gamma)}\omega_{(\alpha,\gamma)}
\right\rangle
\\[2mm]
=\left\langle {\cal A}_c(\overline{y},W_1,W_2)^{-1/2}C^{\,
T}_{(\alpha,\gamma)}(\Theta_c)_{(\alpha,\gamma)}\omega_{(\alpha,\gamma)},
{\cal A}_c(\overline{y},W_1,W_2)^{-1}{\cal
A}_c(\overline{y},W_1,W_2)^{-1/2}C^{\,
T}_{(\alpha,\gamma)}(\Theta_c)_{(\alpha,\gamma)}\omega_{(\alpha,\gamma)}
\right\rangle
\\[2mm]
\leq \left\langle  {\cal A}_c(\overline{y},W_1,W_2)^{-1/2}C^{\,
T}_{(\alpha,\gamma)}(\Theta_c)_{(\alpha,\gamma)}\omega_{(\alpha,\gamma)},
{\widetilde{R}} \overline{{\cal H}}_c {\widetilde{R}} ^T{\cal
A}_c(\overline{y},W_1,W_2)^{-1/2}C^{\,
T}_{(\alpha,\gamma)}(\Theta_c)_{(\alpha,\gamma)}\omega_{(\alpha,\gamma)}
\right\rangle\\[2mm]
\leq \overline{\sigma} {\underline{\sigma}}^{-1}
\underline{\eta}^{-1} \left\langle C^{\,
T}_{(\alpha,\gamma)}(\Theta_c)_{(\alpha,\gamma)}\omega_{(\alpha,\gamma)},
{\cal A}_c(\overline{y},W_1,W_2)^{-1}C^{\,
T}_{(\alpha,\gamma)}(\Theta_c)_{(\alpha,\gamma)}\omega_{(\alpha,\gamma)}
\right\rangle
\\[2mm]
\leq  \overline{\sigma}{\underline{\sigma}}^{-2}
\underline{\eta}^{-2} \left\langle {\widetilde{C}}^{\,
T}_{(\alpha,\gamma)}
(\Theta_c)_{(\alpha,\gamma)}\omega_{(\alpha,\gamma)},
{\widetilde{C}}^{\, T}_{(\alpha,\gamma)}
(\Theta_c)_{(\alpha,\gamma)} \omega_{(\alpha,\gamma)}
 \right\rangle\\[2mm]
 \leq
 \overline{\nu}  \overline{\sigma}{\underline{\sigma}}^{-2} \underline{\eta}^{-2}
\|(\Theta_c)_{(\alpha,\gamma)}\omega_{(\alpha,\gamma)}\|^2
\\[2mm]
\leq \overline{\nu} \overline{\sigma}{\underline{\sigma}}^{-2}
\underline{\eta}^{-2} \displaystyle \left (\max_{i \in \alpha, j \in
\gamma}\lambda_i/(\lambda_i+c|\lambda_j|)\right)^2 \|
\omega_{(\alpha,\gamma)}\|^2
\end{array}
\]
\begin{equation}\label{ineq2-add0}
\begin{array}{l}
\leq \overline{\nu}\overline{\sigma} {\underline{\sigma}}^{-2}
\underline{\eta}^{-2} \overline{\nu}_0^2 (\overline{\nu}_0+c)^{-2}
\| \omega_{(\alpha,\gamma)}\|^2\,
\\[2mm]
\leq \overline{\nu} \overline{\sigma}{\underline{\sigma}}^{-2}
\underline{\eta}^{-2} \overline{\nu}_0^2 c^{-2} \|
\omega_{(\alpha,\gamma)}\|^2\, 
\\[2mm]
\leq \displaystyle \frac{1}{256}\varrho_0^2c^{-2} (2\|
\omega_{(\alpha,\gamma)}\|^2).\, 
\end{array}
\end{equation}
Similarly, we obtain
\begin{equation}\label{eq:est1}
\begin{array}{l}
 \left \langle B_{(a,b_S)}^{\, T}(\Delta_{1/c})_{(a,b_S)}\xi_{(a,b_S)}, {\cal
A}_c(\overline{y},W_1,W_2)^{-2} B_{(a,b_S)}^{\, T}(\Delta_{1/c})_{(a,b_S)}\xi_{(a,b_S)} \right \rangle\\[2mm]
\leq
 \overline{\nu}  \overline{\sigma}{\underline{\sigma}}^{-2} \underline{\eta}^{-2}
\|(\Delta_{1/c})_{(a,b_S)}\xi_{(a,b_S)}\|^2\\[2mm]
\leq \overline{\nu}  \overline{\sigma}{\underline{\sigma}}^{-2} \underline{\eta}^{-2}
\overline \nu_1=\displaystyle \left (\max_{i \in a, j \in \{1,\ldots,|b_S|\}}\displaystyle \frac{1-(w_{b_S})_j}{c\lambda_i(\overline X)}\right)^2 \|\xi_{(a,b_S)}\|^2\\[2mm]
\leq
\displaystyle \frac{1}{16}\varrho_0^2c^{-2} (2 \|\xi_{(a,b_S)}\|^2);\,
\end{array}
\end{equation}
\begin{equation}\label{eq:est2}
\begin{array}{l}
\left \langle B_{(a,b_L)}^{\, T}(\Delta_{1/c})_{(a,b_L)}\xi_{(a,b_L)}, {\cal
A}_c(\overline{y},W_1,W_2)^{-2} B_{(a,b_L)}^{\, T}(\Delta_{1/c})_{(a,b_L)}\xi_{(a,b_L)} \right \rangle\\[2mm]
\leq
 \overline{\nu}  \overline{\sigma}{\underline{\sigma}}^{-2} \underline{\eta}^{-2}
\|(\Delta_{1/c})_{(a,b_L)}\xi_{(a,b_L)}\|^2\\[2mm]
  \leq \overline{\nu}  \overline{\sigma}{\underline{\sigma}}^{-2} \underline{\eta}^{-2}
\overline \nu_2=\displaystyle \left (\max_{i \in a, j \in b_L\}}\displaystyle \frac{2}{c\lambda_i(\overline X)}\right)^2 \|\xi_{(a,b_L)}\|^2\\[2mm]
\leq
\displaystyle \frac{1}{32}\varrho_0^2c^{-2} (2 \|\xi_{(a,b_L)}\|^2);
\end{array}
\end{equation}
\begin{equation}\label{eq:est3}
\begin{array}{l}
 \left \langle B_{(a,c)}^{\, T}(\Delta_{1/c})_{(a,c)}\xi_{(a,c)}, {\cal
A}_c(\overline{y},W_1,W_2)^{-2} B_{(a,c)}^{\, T}(\Delta_{1/c})_{(a,c)}\xi_{(a,c)} \right \rangle\\[2mm]
\leq
 \overline{\nu}  \overline{\sigma}{\underline{\sigma}}^{-2} \underline{\eta}^{-2}
\|(\Delta_{1/c})_{(a,c)}\xi_{(a,c)}\|^2\\[2mm]
  \leq \overline{\nu}  \overline{\sigma}{\underline{\sigma}}^{-2} \underline{\eta}^{-2}
\overline \nu_3=\displaystyle \left (\max_{i \in a, j \in c\}}\displaystyle \frac{2}{c(\lambda_i(\overline X)-\lambda_j(\overline X))}\right)^2 \|\xi_{(a,c)}\|^2\\[2mm]
\leq
\displaystyle \frac{1}{64}\varrho_0^2c^{-2} (2  \|\xi_{(a,c)}\|^2);
\end{array}
\end{equation}
\begin{equation}\label{eq:est4}
\begin{array}{l}
\left \langle B_{(c,b_U)}^{\, T}(\Delta_{1/c})_{(c,b_U)}\xi_{(c,b_U)}, {\cal
A}_c(\overline{y},W_1,W_2)^{-2} B_{(c,b_U)}^{\, T}(\Delta_{1/c})_{(c,b_U)}\xi_{(c,b_U)} \right \rangle\\[2mm]
\leq
 \overline{\nu}  \overline{\sigma}{\underline{\sigma}}^{-2} \underline{\eta}^{-2}
\|(\Delta_{1/c})_{(c,b_U)}\xi_{(c,b_U)}\|^2\\[2mm]
  \leq \overline{\nu}  \overline{\sigma}{\underline{\sigma}}^{-2} \underline{\eta}^{-2}
\overline \nu_4=\displaystyle \left (\max_{i \in b_U, j \in c\}}\displaystyle \frac{2}{-c\lambda_j(\overline X)}\right)^2 \|\xi_{(b_U,c)}\|^2\\[2mm]
\leq
\displaystyle \frac{1}{128}\varrho_0^2c^{-2} (2  \|\xi_{(b_U,c)}\|^2)
\end{array}
\end{equation}
and
\begin{equation}\label{eq:est5}
\begin{array}{l}
\left \langle B_{(c,b_S)}^{\, T}(\Delta_{1/c})_{(c,b_S)}\xi_{(c,b_S)}, {\cal
A}_c(\overline{y},W_1,W_2)^{-2} B_{(c,b_S)}^{\, T}(\Delta_{1/c})_{(c,b_S)}\xi_{(c,b_S)} \right \rangle\\[2mm]
\leq
 \overline{\nu}  \overline{\sigma}{\underline{\sigma}}^{-2} \underline{\eta}^{-2}
\|(\Delta_{1/c})_{(c,b_S)}\xi_{(c,b_S)}\|^2\\[2mm]
\leq \overline{\nu}  \overline{\sigma}{\underline{\sigma}}^{-2} \underline{\eta}^{-2}
\overline \nu_5=\displaystyle \left (\max_{i \in \{1,\ldots,|b_S|\},j \in c}\displaystyle \frac{(w_{b_S})_i+1}{-c\lambda_j(\overline X)}\right)^2 \|\xi_{(b_S,c)}\|^2\\[2mm]
\leq
\displaystyle \frac{1}{256}\varrho_0^2c^{-2} (2  \|\xi_{(b_S,c)}\|^2).
\end{array}
\end{equation}
Combining (\ref{eq:est1})-(\ref{eq:est5}) with (\ref{Ahh-0}) and (\ref{anup-1-add0}), we obtain
\[
  \left \langle (x_c)^\prime(\overline{y};\Delta y),
(x_c)^\prime(\overline{y};\Delta y)\right\rangle \le
\varrho_0^2\|\Delta y\|^2 /c^2.
\]
 Thus (\ref{eq:direction-added-sdp}) holds for $\mu_0\geq \varrho_0$.

Now we prove (\ref{import11}) for some $\mu_0\geq \varrho_0$. Let
$V(\Delta y) \in \overline{{\cal V}}_c(\Delta y)$. Then from the
definition of $\overline{{\cal V}}_c(\Delta y)$,  there exist $W_1\in \partial_B [{\rm D} \theta_c]^*(F(\overline x)+\overline Y/c)$ and $W_2 \in
\partial_B \Pi_{{\cal S}^p_+}(\overline{\Gamma}-c g(\overline{x}))$
such that
\[
\begin{array}{l}
V(\Delta y)\\[6pt]
= \left
 [
 \begin{array}{c}
 c^{-1}W_1{\rm D}F(\overline x)\\
{\cal J} h (\overline{x})\\
 -W_2{\rm D}g(\overline{x})
\end{array}
\right
 ]{\cal A}_c(\overline{y},W_1,W_2)^{-1}\left[-c^{-1}{\rm D}F(\overline x)^*W_1(\Delta Y)-{\cal J}
h(\overline{x})^T\Delta\mu+{\rm D}g(\overline{x})^*W_2(\Delta\Gamma)\right]\\[8mm]
 \quad \quad \quad \quad \quad \quad\quad \quad \quad \quad \quad \quad \quad \quad \quad  +\left
 (
 \begin{array}{c}
 -c^{-1}\Delta Y+c^{-2}W_1\Delta Y\\
 0\\
-c^{-1}\Delta \Gamma+ c^{-1}W_2( \Delta \Gamma)
\end{array}
 \right
 ).
 \end{array}
\]
For notational convenience, we assume that $(W_1,W_2) \in \partial_B [{\rm D} \theta_c]^*(F(\overline x)+\overline Y/c) \times
\partial_B \Pi_{{\cal S}^p_+}(\overline{\Gamma}-cg(\overline{x}))$ is
the same as in (\ref{comxderrii-sdp}).
 After direct calculations, we obtain
\begin{equation}\label{Ahh}
\begin{array}{l}
\quad -\left\langle V(\Delta y), \Delta y \right \rangle=
\\[2mm]
[A^T D_c\Delta d_0+2 B_{(a,b_S)}^{\, T}(\Delta_{1/c})_{(a,b_S)}\xi_{(a,b_S)}\\[2mm]
+2 B_{(a,b_L)}^{\, T}(\Delta_{1/c})_{(a,b_L)}\xi_{(a,b_L)}
 +2 B_{(a,c)}^{\, T}(\Delta_{1/c})_{(a,c)}\xi_{(a,c)}\\[2mm]
 +2 B_{(c,b_U)}^{\, T}(\Delta_{1/c})_{(c,b_U)} \xi_{(c,b_U)}
 + 2 B_{(c,b_S)}^{\, T}(\Delta_{1/c})_{(c,b_S)}\xi_{(c,b_S)}\\[2mm]
  -2C_{(\alpha,\gamma)}^{\,
T}(\Theta_c)_{(\alpha,\gamma)}\omega_{(\alpha,\gamma)}]^{\,
T}{\cal
A}_c(\overline{y},W_1,W_2)^{-1}
[A^T D_c\Delta d_0+2 B_{(a,b_S)}^{\, T}(\Delta_{1/c})_{(a,b_S)}\xi_{(a,b_S)}\\[2mm]
+2 B_{(a,b_L)}^{\, T}(\Delta_{1/c})_{(a,b_L)}\xi_{(a,b_L)}
 +2 B_{(a,c)}^{\, T}(\Delta_{1/c})_{(a,c)}\xi_{(a,c)}\\[2mm]
 +2 B_{(c,b_U)}^{\, T}(\Delta_{1/c})_{(c,b_U)} \xi_{(c,b_U)}
 + 2 B_{(c,b_S)}^{\, T}(\Delta_{1/c})_{(c,b_S)}\xi_{(c,b_S)}
  -2C_{(\alpha,\gamma)}^{\,
T}(\Theta_c)_{(\alpha,\gamma)}\omega_{(\alpha,\gamma)}]\\[4mm]
+c^{-1}\|\Delta Y\|^2-c^{-2}\left \langle W_1\Delta Y, \Delta Y \right \rangle
+ {c}^{-1}\|\Delta\Gamma\|^2 -c^{-1}\left\langle
\Delta\Gamma, W_2(\Delta\Gamma) \right \rangle\,
\end{array}
\end{equation}

Next, we estimate the lower and upper bounds of the right hand side
of (\ref{Ahh}).
 By using (\ref{singvalue-in}) and Lemma
\ref{le2} we obtain
\[
\underline{{\cal E}}_c \preceq
 A {\cal A}_c(\overline{y},W_1,W_2)^{-1}A^T
 \preceq\overline{{\cal E}}_c\, .
\]
Thus, for $l_U=|b_U|(|b_U|+1)/2$, $l_L=|b_L|(|b_L|+1)/2$ and $l_\beta= |\beta|(|\beta|+1)/2$, we have
\[
\begin{array}[b]{l}
\quad \left\langle A^TD_c \Delta d_0, {\cal
A}_c(\overline{y},W_1,W_2)^{-1}A^TD_c \Delta d_0 \right\rangle \geq
\left\langle D_c \Delta d_0,\underline{{\cal E}}_c D_c \Delta d_0
\right\rangle
\\[2mm]
\ge
\left({\overline{\sigma}}\overline{\eta}+(c-c_0)\right)^{-1}\|(\Delta\mu,\xi_{(b_S,b_S)},{
\omega}_{(\alpha,\alpha)})\|^2\\[2mm]
\quad + {4}\left(
{\overline{\sigma}}\overline{\eta}+2(c-c_0)\right)^{-1} \|(\xi_{(b_U,b_S)},\xi_{(b_U,b_L)},\xi_{(b_S,b_L)},\omega_{(\alpha,\beta)})\|^2
\\[2mm]
\, +\left\langle ({\widehat \Delta_{1/c}})_{(b_U,b_U)}{\widehat
\xi}_{(b_U,b_U)}, \left({\overline{\sigma}}\overline{\eta}
I_{l_U}+(c-c_0) ({\widehat
\Delta_{1/c}})_{(b_U,b_U)} \right )^{-1} ({\widehat
\Delta_{1/c}})_{(b_U,b_U)}{\widehat \xi}_{(b_U,b_U)}
\right\rangle\\[2mm]
\, +\left\langle ({\widehat \Delta_{1/c}})_{(b_L,b_L)}{\widehat
\xi}_{(b_L,b_L)}, \left({\overline{\sigma}}\overline{\eta}
I_{l_L}+(c-c_0) ({\widehat
\Delta_{1/c}})_{(b_L,b_L)} \right )^{-1} ({\widehat
\Delta_{1/c}})_{(b_L,b_L)}{\widehat \xi}_{(b_L,b_L)}
\right\rangle\\[2mm]
\, +\left\langle ({\widehat \Theta_c})_{(\beta, \beta)}{\widehat
\omega}_{(\beta,\beta)}, \left({\overline{\sigma}}\overline{\eta}
I_{l_\beta}+(c-c_0) ({\widehat
\Theta_c})_{(\beta,\beta)} \right )^{-1} ({\widehat
\Theta_c})_{(\beta,\beta)}{\widehat \omega}_{(\beta,\beta)}
\right\rangle
\end{array}
\]
\begin{equation}\label{anup-opposite}
\begin{array}{l}
\ge
\left({\overline{\sigma}}\overline{\eta}+(c-c_0)\right)^{-1}\|(\Delta\mu,\xi_{(b_S,b_S)},{
\omega}_{(\alpha,\alpha)})\|^2\\[2mm]
\quad + {4}\left(
{\overline{\sigma}}\overline{\eta}+2(c-c_0)\right)^{-1} \|(\xi_{(b_U,b_S)},\xi_{(b_U,b_L)},\xi_{(b_S,b_L)},\omega_{(\alpha,\beta)})\|^2
\\[2mm]
\, +\left\langle ({\Delta_{1/c}})_{(b_U,b_U)}{
\xi}_{(b_U,b_U)}, \left({\overline{\sigma}}\overline{\eta}
I_{|b_U|}+(c-c_0) ({
\Delta_{1/c}})_{(b_U,b_U)} \right )^{-1} ({
\Delta_{1/c}})_{(b_U,b_U)}{\xi}_{(b_U,b_U)}
\right\rangle\\[2mm]
\, +\left\langle ({ \Delta_{1/c}})_{(b_L,b_L)}{
\xi}_{(b_L,b_L)}, \left({\overline{\sigma}}\overline{\eta}
I_{|b_L|}+(c-c_0) ({
\Delta_{1/c}})_{(b_L,b_L)} \right )^{-1} ({
\Delta_{1/c}})_{(b_L,b_L)}{ \xi}_{(b_L,b_L)}
\right\rangle\\[2mm]
\, +\left\langle ({ \Theta_c})_{(\beta, \beta)}{
\omega}_{(\beta,\beta)}, \left({\overline{\sigma}}\overline{\eta}
I_{|\beta|}+(c-c_0) ({ \Theta_c})_{(\beta,\beta)} \right )^{-1} ({
\Theta_c})_{(\beta,\beta)}{ \omega}_{(\beta,\beta)}
\right\rangle\,
\end{array}
\end{equation}
and
\begin{equation}\label{anup-1}
\begin{array}[b]{l}
\quad \left\langle A^TD_c \Delta d_0, {\cal
A}_c(\overline{y},W_1,W_2)^{-1}A^TD_c \Delta d_0 \right\rangle \leq
\left\langle D_c \Delta d_0,\overline{{\cal E}}_c D_c \Delta
d_0\right\rangle
\\[2mm]
\le
\left({\underline{\sigma}}\underline{\eta}/2+(c-c_0)\right)^{-1}\|(\Delta\mu,\xi_{(b_S,b_S)},{
\omega}_{(\alpha,\alpha)})\|^2\\[2mm]
\quad + {4}\left(
{\underline{\sigma}}\underline{\eta}+2(c-c_0)\right)^{-1} \|(\xi_{(b_U,b_S)},\xi_{(b_U,b_L)},\xi_{(b_S,b_L)},\omega_{(\alpha,\beta)})\|^2
\\[2mm]
\, +\left\langle ({\widehat \Delta_{1/c}})_{(b_U,b_U)}{\widehat
\xi}_{(b_U,b_U)}, \left({\underline{\sigma}}\underline{\eta}
I_{l_U}+(c-c_0) ({\widehat
\Delta_{1/c}})_{(b_U,b_U)} \right )^{-1} ({\widehat
\Delta_{1/c}})_{(b_U,b_U)}{\widehat \xi}_{(b_U,b_U)}
\right\rangle\\[2mm]
\, +\left\langle ({\widehat \Delta_{1/c}})_{(b_L,b_L)}{\widehat
\xi}_{(b_L,b_L)}, \left({\underline{\sigma}}\underline{\eta}
I_{l_L}+(c-c_0) ({\widehat
\Delta_{1/c}})_{(b_L,b_L)} \right )^{-1} ({\widehat
\Delta_{1/c}})_{(b_L,b_L)}{\widehat \xi}_{(b_L,b_L)}
\right\rangle\\[2mm]
\, +\left\langle ({\widehat \Theta_c})_{(\beta, \beta)}{\widehat
\omega}_{(\beta,\beta)}, \left({\underline{\sigma}}\underline{\eta}
I_{l_\beta}+(c-c_0) ({\widehat
\Theta_c})_{(\beta,\beta)} \right )^{-1} ({\widehat
\Theta_c})_{(\beta,\beta)}{\widehat \omega}_{(\beta,\beta)}
\right\rangle
\\[2mm]
\le
\left({\underline{\sigma}}\underline{\eta}/2+(c-c_0)\right)^{-1}\|(\Delta\mu,\xi_{(b_S,b_S)},{
\omega}_{(\alpha,\alpha)})\|^2\\[2mm]
\quad + {4}\left(
{\underline{\sigma}}\underline{\eta}+2(c-c_0)\right)^{-1} \|(\xi_{(b_U,b_S)},\xi_{(b_U,b_L)},\xi_{(b_S,b_L)},\omega_{(\alpha,\beta)})\|^2
\\[2mm]
\, +\left\langle ({\Delta_{1/c}})_{(b_U,b_U)}{
\xi}_{(b_U,b_U)}, \left({\underline{\sigma}}\underline{\eta}
I_{|b_U|}+(c-c_0) ({
\Delta_{1/c}})_{(b_U,b_U)} \right )^{-1} ({
\Delta_{1/c}})_{(b_U,b_U)}{\xi}_{(b_U,b_U)}
\right\rangle\\[2mm]
\, +\left\langle ({ \Delta_{1/c}})_{(b_L,b_L)}{
\xi}_{(b_L,b_L)}, \left({\underline{\sigma}}\underline{\eta}
I_{|b_L|}+(c-c_0) ({
\Delta_{1/c}})_{(b_L,b_L)} \right )^{-1} ({
\Delta_{1/c}})_{(b_L,b_L)}{ \xi}_{(b_L,b_L)}
\right\rangle\\[2mm]
\, +\left\langle ({ \Theta_c})_{(\beta, \beta)}{
\omega}_{(\beta,\beta)}, \left({\underline{\sigma}}\underline{\eta}
I_{|\beta|}+(c-c_0) ({ \Theta_c})_{(\beta,\beta)} \right )^{-1} ({
\Theta_c})_{(\beta,\beta)}{ \omega}_{(\beta,\beta)}
\right\rangle.
\end{array}
\end{equation}
By recalling  that
\[
\widetilde{C}=C{\widetilde
R}\quad  {\rm and}\quad  {\widetilde R} =R \left [
\begin{array}{cc}
\Sigma^{-1}U^T & 0\\
 0 & I_{n_3}
\end{array}
\right ] \]
and
\[
\underline{\nu} \|s\|^2 \leq \max\left\{ \left\langle s,
\widetilde{C}\widetilde{C}^{\,
T} s \right\rangle, \left\langle s,
C C^Ts \right\rangle
\right\}
 \le \overline{\nu}\|s\|^2,\, \forall s,
\]
from Lemma \ref{le2}, (\ref{eu1-A}),
(\ref{eq:nu-lower-upper}), and (\ref{eq:two-nus}) we know  that
\[
\begin{array}[b]{l}
\displaystyle \langle  [B_{(a,b_S)}^{\, T}(\Delta_{1/c})_{(a,b_S)}\xi_{(a,b_S)}
 +B_{(a,b_L)}^{\, T}(\Delta_{1/c})_{(a,b_L)}\xi_{(a,b_L)}
 + B_{(a,c)}^{\, T}(\Delta_{1/c})_{(a,c)}\xi_{(a,c)}\\[2mm]
 \quad + B_{(c,b_U)}^{\, T}(\Delta_{1/c})_{(c,b_U)} \xi_{(c,b_U)}
 +  B_{(c,b_S)}^{\, T}(\Delta_{1/c})_{(c,b_S)}\xi_{(c,b_S)}
  -C_{(\alpha,\gamma)}^{\,
T}(\Theta_c)_{(\alpha,\gamma)}\omega_{(\alpha,\gamma)}],\\[2mm]
{\cal
A}_c(\overline{y},W_1,W_2)^{-1}
[ B_{(a,b_S)}^{\, T}(\Delta_{1/c})_{(a,b_S)}\xi_{(a,b_S)}
+B_{(a,b_L)}^{\, T}(\Delta_{1/c})_{(a,b_L)}\xi_{(a,b_L)}\\[2mm]
 \quad \quad \quad \quad \quad \quad \quad \quad \quad \quad+ B_{(a,c)}^{\, T}(\Delta_{1/c})_{(a,c)}\xi_{(a,c)}
 + B_{(c,b_U)}^{\, T}(\Delta_{1/c})_{(c,b_U)} \xi_{(c,b_U)}\\[2mm]
\quad \quad \quad \quad \quad \quad \quad \quad \quad \quad+ B_{(c,b_S)}^{\, T}(\Delta_{1/c})_{(c,b_S)}\xi_{(c,b_S)}
  -C_{(\alpha,\gamma)}^{\,
T}(\Theta_c)_{(\alpha,\gamma)}\omega_{(\alpha,\gamma)}] \rangle
\end{array}
\]
\begin{equation}\label{ineq1}
\begin{array}{l}
\geq
\displaystyle \langle  [\widetilde B_{(a,b_S)}^{\, T}(\Delta_{1/c})_{(a,b_S)}\xi_{(a,b_S)}
+\widetilde B_{(a,b_L)}^{\, T}(\Delta_{1/c})_{(a,b_L)}\xi_{(a,b_L)}
 + \widetilde B_{(a,c)}^{\, T}(\Delta_{1/c})_{(a,c)}\xi_{(a,c)}\\[2mm]
 \quad  \, + \widetilde B_{(c,b_U)}^{\, T}(\Delta_{1/c})_{(c,b_U)} \xi_{(c,b_U)}
 +  \widetilde B_{(c,b_S)}^{\, T}(\Delta_{1/c})_{(c,b_S)}\xi_{(c,b_S)}
  -\widetilde C_{(\alpha,\gamma)}^{\,
T}(\Theta_c)_{(\alpha,\gamma)}\omega_{(\alpha,\gamma)}],\\[2mm]
\underline{{\cal H}}_c
[ \widetilde B_{(a,b_S)}^{\, T}(\Delta_{1/c})_{(a,b_S)}\xi_{(a,b_S)}
+\widetilde B_{(a,b_L)}^{\, T}(\Delta_{1/c})_{(a,b_L)}\xi_{(a,b_L)}
 + \widetilde B_{(a,c)}^{\, T}(\Delta_{1/c})_{(a,c)}\xi_{(a,c)}\\[2mm]
\quad  \, + \widetilde B_{(c,b_U)}^{\, T}(\Delta_{1/c})_{(c,b_U)} \xi_{(c,b_U)}
 + \widetilde B_{(c,b_S)}^{\, T}(\Delta_{1/c})_{(c,b_S)}\xi_{(c,b_S)}
  -\widetilde C_{(\alpha,\gamma)}^{\,
T}(\Theta_c)_{(\alpha,\gamma)}\omega_{(\alpha,\gamma)}] \rangle\\[2mm]
\geq
\left(
{\overline{\sigma}}\overline{\eta}+2(c-c_0)\right)^{-1}
\displaystyle \langle  [\widetilde B_{(a,b_S)}^{\, T}(\Delta_{1/c})_{(a,b_S)}\xi_{(a,b_S)}
+\widetilde B_{(a,b_L)}^{\, T}(\Delta_{1/c})_{(a,b_L)}\xi_{(a,b_L)}\\[2mm]
\quad  \, + \widetilde B_{(a,c)}^{\, T}(\Delta_{1/c})_{(a,c)}\xi_{(a,c)}
 + \widetilde B_{(c,b_U)}^{\, T}(\Delta_{1/c})_{(c,b_U)} \xi_{(c,b_U)}
 +  \widetilde B_{(c,b_S)}^{\, T}(\Delta_{1/c})_{(c,b_S)}\xi_{(c,b_S)}\\[2mm]
 \quad  \, -\widetilde C_{(\alpha,\gamma)}^{\,T}(\Theta_c)_{(\alpha,\gamma)}\omega_{(\alpha,\gamma)}],
[ \widetilde B_{(a,b_S)}^{\, T}(\Delta_{1/c})_{(a,b_S)}\xi_{(a,b_S)}
+\widetilde B_{(a,b_L)}^{\, T}(\Delta_{1/c})_{(a,b_L)}\xi_{(a,b_L)}\\[2mm]
 \quad  \,+ \widetilde B_{(a,c)}^{\, T}(\Delta_{1/c})_{(a,c)}\xi_{(a,c)}
 + \widetilde B_{(c,b_U)}^{\, T}(\Delta_{1/c})_{(c,b_U)} \xi_{(c,b_U)}
 + \widetilde B_{(c,b_S)}^{\, T}(\Delta_{1/c})_{(c,b_S)}\xi_{(c,b_S)}\\[2mm]
 \quad  \, -\widetilde C_{(\alpha,\gamma)}^{\,
T}(\Theta_c)_{(\alpha,\gamma)}\omega_{(\alpha,\gamma)}] \rangle\\[2mm]
\geq
\underline{\nu}\left({\overline{\sigma}}\overline{\eta}+2(c-c_0)\right)^{-1}
\|[(\Delta_{1/c})_{(a,b_S)}\xi_{(a,b_S)},
(\Delta_{1/c})_{(a,b_L)}\xi_{(a,b_L)},(\Delta_{1/c})_{(a,c)}\xi_{(a,c)},\\[2mm]
\,\, \quad \quad \quad \quad \quad \quad
 (\Delta_{1/c})_{(c,b_U)} \xi_{(c,b_U)},(\Delta_{1/c})_{(c,b_S)}\xi_{(c,b_S)},(\Theta_c)_{(\alpha,\gamma)}\omega_{(\alpha,\gamma)}]\|^2
\\[2mm]
\ge \underline{\nu}\left({\overline{\sigma}}\overline{\eta}+2(c-c_0)\right)^{-1} \left [
\displaystyle \left (\min_{i \in a,1\leq j \leq |b_S|}\displaystyle \frac{1-(w_{b_S})_j}{c\lambda_i(\overline X)+(1-(w_{b_S})_j)}\right)^2 \|\xi_{(a,b_S)}\|^2\right.
\\[2mm]
+
\displaystyle \left (\min_{i \in a}\displaystyle \frac{2}{c\lambda_i(\overline X)+2}\right)^2 \|
\xi_{(a,b_L)}\|^2
+
\displaystyle \left (\min_{i \in a,j\in c}\displaystyle \frac{2}{c[\lambda_i(\overline X)-\lambda_i(\overline X)]+2}\right)^2 \|\xi_{(a,c)}\|^2
\\[2mm]
+
\displaystyle \left (\min_{i \in c, 1\leq j\leq |b_S|}\displaystyle \frac{(w_{b_S})_j+1}{((w_{b_S})_j+1)-c\lambda_i(\overline X )}\right)^2 \|\xi_{(c,b_S)}\|^2
\\[2mm]
\left. +
\displaystyle \left (\min_{i \in c}\displaystyle \frac{2}{-c\lambda_i(\overline X)+2}\right)^2 \|\xi_{(c,b_U)}\|^2+
\displaystyle \left (\min_{i \in \alpha, j \in
\gamma}\lambda_i/(\lambda_i+c|\lambda_j|)\right)^2 \|
\omega_{(\alpha,\gamma)}\|^2\right]
\\[2mm]
\ge
\underline{\nu}
\left({\overline{\sigma}}\overline{\eta}+2(c-c_0)\right)^{-1}
\left[ \underline{\nu}_{a,b_S}^2 (\underline{\nu}_{a,b_S}+c)^{-2} \|
\xi_{(a,b_S)}\|^2+
\underline{\nu}_{a,b_L}^2 (\underline{\nu}_{a,b_L}+c)^{-2} \|
\xi_{(a,b_L)}\|^2\right.\\[2mm]
\quad\, +\underline{\nu}_{a,c}^2 (\underline{\nu}_{a,c}+c)^{-2} \|
\xi_{(a,c)}\|^2
+\underline{\nu}_{c,b_U}^2 (\underline{\nu}_{c,b_U}+c)^{-2} \|
\xi_{(c,b_U)}\|^2\\[2mm]
\quad \, \left.+\underline{\nu}_{c,b_S}^2 (\underline{\nu}_{c,b_S}+c)^{-2} \|
\xi_{(c,b_S)}\|^2
+\underline{\nu}_{\alpha,\gamma}^2 (\underline{\nu}_{\alpha,\gamma}+c)^{-2} \|
\omega_{(\alpha,\gamma)}\|^2\right]\\[2mm]
\geq \underline{\nu}
\left({\overline{\sigma}}\overline{\eta}+2(c-c_0)\right)^{-1}
\underline{\nu}_0^2 (\underline{\nu}_0+c)^{-2}
  \|(\xi_{(a,b_S)},\xi_{(a,b_L)},
\xi_{(a,c)},\xi_{(c,b_S)},
\omega_{(\alpha,\gamma)})\|^2.
\end{array}
\end{equation}
Similarly, we get
\[
\begin{array}[b]{l}
 \displaystyle \langle  [B_{(a,b_S)}^{\, T}(\Delta_{1/c})_{(a,b_S)}\xi_{(a,b_S)}
+B_{(a,b_L)}^{\, T}(\Delta_{1/c})_{(a,b_L)}\xi_{(a,b_L)}
  + B_{(a,c)}^{\, T}(\Delta_{1/c})_{(a,c)}\xi_{(a,c)}\\[2mm]
\quad \,  + B_{(c,b_U)}^{\, T}(\Delta_{1/c})_{(c,b_U)} \xi_{(c,b_U)}
 +  B_{(c,b_S)}^{\, T}(\Delta_{1/c})_{(c,b_S)}\xi_{(c,b_S)}
  -C_{(\alpha,\gamma)}^{\,
T}(\Theta_c)_{(\alpha,\gamma)}\omega_{(\alpha,\gamma)}],\\[2mm]
{\cal
A}_c(\overline{y},W_1,W_2)^{-1}
[ B_{(a,b_S)}^{\, T}(\Delta_{1/c})_{(a,b_S)}\xi_{(a,b_S)}
+B_{(a,b_L)}^{\, T}(\Delta_{1/c})_{(a,b_L)}\xi_{(a,b_L)}\\[2mm]
 \quad \,\quad \,\quad \,\quad \, \quad \,\quad \,\quad \,\quad+ B_{(a,c)}^{\, T}(\Delta_{1/c})_{(a,c)}\xi_{(a,c)}
 + B_{(c,b_U)}^{\, T}(\Delta_{1/c})_{(c,b_U)} \xi_{(c,b_U)}\\[2mm]
 \quad \,\quad \, \quad \,\quad \,\quad \,\quad \,\quad \,\quad + B_{(c,b_S)}^{\, T}(\Delta_{1/c})_{(c,b_S)}\xi_{(c,b_S)}
  -C_{(\alpha,\gamma)}^{\,
T}(\Theta_c)_{(\alpha,\gamma)}\omega_{(\alpha,\gamma)}] \rangle
\end{array}
\]
\begin{equation}\label{ineq2}
\begin{array}{l}
\leq
\displaystyle \langle  [\widetilde B_{(a,b_S)}^{\, T}(\Delta_{1/c})_{(a,b_S)}\xi_{(a,b_S)}
+\widetilde B_{(a,b_L)}^{\, T}(\Delta_{1/c})_{(a,b_L)}\xi_{(a,b_L)}
 + \widetilde B_{(a,c)}^{\, T}(\Delta_{1/c})_{(a,c)}\xi_{(a,c)}\\[2mm]
\quad \, + \widetilde B_{(c,b_U)}^{\, T}(\Delta_{1/c})_{(c,b_U)} \xi_{(c,b_U)}
 +  \widetilde B_{(c,b_S)}^{\, T}(\Delta_{1/c})_{(c,b_S)}\xi_{(c,b_S)}
  -\widetilde C_{(\alpha,\gamma)}^{\,
T}(\Theta_c)_{(\alpha,\gamma)}\omega_{(\alpha,\gamma)}],\\[2mm]
\overline{{\cal H}}_c
[ \widetilde B_{(a,b_S)}^{\, T}(\Delta_{1/c})_{(a,b_S)}\xi_{(a,b_S)}
+\widetilde B_{(a,b_L)}^{\, T}(\Delta_{1/c})_{(a,b_L)}\xi_{(a,b_L)}
 + \widetilde B_{(a,c)}^{\, T}(\Delta_{1/c})_{(a,c)}\xi_{(a,c)}\\[2mm]
 \quad \, + \widetilde B_{(c,b_U)}^{\, T}(\Delta_{1/c})_{(c,b_U)} \xi_{(c,b_U)}
 + \widetilde B_{(c,b_S)}^{\, T}(\Delta_{1/c})_{(c,b_S)}\xi_{(c,b_S)}
  -\widetilde C_{(\alpha,\gamma)}^{\,
T}(\Theta_c)_{(\alpha,\gamma)}\omega_{(\alpha,\gamma)}] \rangle\\[2mm]
\leq
\underline{\sigma}^{-1}\underline{\eta}^{-1}
\displaystyle \langle  [\widetilde B_{(a,b_S)}^{\, T}(\Delta_{1/c})_{(a,b_S)}\xi_{(a,b_S)}
+\widetilde B_{(a,b_L)}^{\, T}(\Delta_{1/c})_{(a,b_L)}\xi_{(a,b_L)}\\[2mm]
 \quad \,+ \widetilde B_{(a,c)}^{\, T}(\Delta_{1/c})_{(a,c)}\xi_{(a,c)}
 + \widetilde B_{(c,b_U)}^{\, T}(\Delta_{1/c})_{(c,b_U)} \xi_{(c,b_U)}
 +  \widetilde B_{(c,b_S)}^{\, T}(\Delta_{1/c})_{(c,b_S)}\xi_{(c,b_S)}\\[2mm]
 \quad \, -\widetilde C_{(\alpha,\gamma)}^{\,T}(\Theta_c)_{(\alpha,\gamma)}\omega_{(\alpha,\gamma)}],
[ \widetilde B_{(a,b_S)}^{\, T}(\Delta_{1/c})_{(a,b_S)}\xi_{(a,b_S)}
+\widetilde B_{(a,b_L)}^{\, T}(\Delta_{1/c})_{(a,b_L)}\xi_{(a,b_L)}\\[2mm]
 \quad \,+ \widetilde B_{(a,c)}^{\, T}(\Delta_{1/c})_{(a,c)}\xi_{(a,c)}
 + \widetilde B_{(c,b_U)}^{\, T}(\Delta_{1/c})_{(c,b_U)} \xi_{(c,b_U)}
 + \widetilde B_{(c,b_S)}^{\, T}(\Delta_{1/c})_{(c,b_S)}\xi_{(c,b_S)}\\[2mm]
 \quad \, -\widetilde C_{(\alpha,\gamma)}^{\,
T}(\Theta_c)_{(\alpha,\gamma)}\omega_{(\alpha,\gamma)}] \rangle\\[2mm]
\leq
\overline{\nu}\underline{\sigma}^{-1}\underline{\eta}^{-1}
\|[(\Delta_{1/c})_{(a,b_S)}\xi_{(a,b_S)},
(\Delta_{1/c})_{(a,b_L)}\xi_{(a,b_L)},(\Delta_{1/c})_{(a,c)}\xi_{(a,c)},\\[2mm]
\,\, \quad \quad \quad \quad \quad \quad
 (\Delta_{1/c})_{(c,b_U)} \xi_{(c,b_U)},(\Delta_{1/c})_{(c,b_S)}\xi_{(c,b_S)},(\Theta_c)_{(\alpha,\gamma)}\omega_{(\alpha,\gamma)}]\|^2
\\[2mm]
\leq \overline{\nu}\underline{\sigma}^{-1}\underline{\eta}^{-1}
\displaystyle \left (\max_{i \in a,1\leq j \leq |b_S|}\displaystyle \frac{1-(w_{b_S})_j}{c\lambda_i(F(\overline x))+(1-(w_{b_S})_j)}\right)^2 \|
\xi_{(a,b_S)}\|^2
\\[2mm]
\quad \,+\overline{\nu}\underline{\sigma}^{-1}\underline{\eta}^{-1}
\displaystyle \left (\max_{i \in a}\displaystyle \frac{2}{c\lambda_i(F(\overline x))+2}\right)^2 \|
\xi_{(a,b_L)}\|^2
\\[2mm]
\quad \,+\overline{\nu}\underline{\sigma}^{-1}\underline{\eta}^{-1}
\displaystyle \left (\max_{i \in a,j\in c}\displaystyle \frac{2}{c[\lambda_i(F(\overline x))-\lambda_i(F(\overline x))]+2}\right)^2 \|
\xi_{(a,c)}\|^2
\\[2mm]
\quad \,+\overline{\nu}\underline{\sigma}^{-1}\underline{\eta}^{-1}
\displaystyle \left (\max_{i \in c}\displaystyle \frac{2}{-c\lambda_i(F(\overline x))+2}\right)^2 \|
\xi_{(c,b_U)}\|^2
\\[2mm]
\quad \,+\overline{\nu}\underline{\sigma}^{-1}\underline{\eta}^{-1}
\displaystyle \left (\max_{i \in c, 1\leq j\leq |b_S|}\displaystyle \frac{(w_{b_S})_j+1}{((w_{b_S})_j+1)-c\lambda_i(F(\overline x))}\right)^2 \|
\xi_{(c,b_S)}\|^2
\\[2mm]
\quad \,+
\overline{\nu}\underline{\sigma}^{-1}\underline{\eta}^{-1}
\displaystyle \left (\max_{i \in \alpha, j \in
\gamma}\lambda_i/(\lambda_i+c|\lambda_j|)\right)^2 \|
\omega_{(\alpha,\gamma)}\|^2
\\[2mm]
\leq
\overline{\nu}\underline{\sigma}^{-1}\underline{\eta}^{-1}
\overline{\nu}_{a,b_S}^2 (\overline{\nu}_{a,b_S}+c)^{-2} \|
\xi_{(a,b_S)}\|^2
+\overline{\nu}\underline{\sigma}^{-1}\underline{\eta}^{-1}
\overline{\nu}_{a,b_L}^2 (\overline{\nu}_{a,b_L}+c)^{-2} \|
\xi_{(a,b_L)}\|^2\\[2mm]
\quad \, +\overline{\nu}\underline{\sigma}^{-1}\underline{\eta}^{-1}
\overline{\nu}_{a,c}^2 (\overline{\nu}_{a,c}+c)^{-2} \|
\xi_{(a,c)}\|^2
+\overline{\nu}\underline{\sigma}^{-1}\underline{\eta}^{-1}
\overline{\nu}_{c,b_U}^2 (\overline{\nu}_{c,b_U}+c)^{-2} \|
\xi_{(c,b_U)}\|^2\\[2mm]
\quad \, +\overline{\nu}\underline{\sigma}^{-1}\underline{\eta}^{-1}
\overline{\nu}_{c,b_S}^2 (\overline{\nu}_{c,b_S}+c)^{-2} \|
\xi_{(c,b_S)}\|^2
+\overline{\nu}\underline{\sigma}^{-1}\underline{\eta}^{-1}
\overline{\nu}_{\alpha,\gamma}^2 (\overline{\nu}_{\alpha,\gamma}+c)^{-2} \|
\omega_{(\alpha,\gamma)}\|^2\\[2mm]
\leq \overline{\nu}\underline{\sigma}^{-1}\underline{\eta}^{-1}
\overline{\nu}_0^2 (\overline{\nu}_0+c)^{-2}  \|(\xi_{(a,b_S)},\xi_{(a,b_L)},
\xi_{(a,c)},\xi_{(c,b_S)},
\omega_{(\alpha,\gamma)})\|^2.
\end{array}
\end{equation}
and
\[
\begin{array}[b]{l}
 \displaystyle \|  B_{(a,b_S)}^{\, T}(\Delta_{1/c})_{(a,b_S)}\xi_{(a,b_S)}
+B_{(a,b_L)}^{\, T}(\Delta_{1/c})_{(a,b_L)}\xi_{(a,b_L)}
 + B_{(a,c)}^{\, T}(\Delta_{1/c})_{(a,c)}\xi_{(a,c)}\\
\quad \, + B_{(c,b_U)}^{\, T}(\Delta_{1/c})_{(c,b_U)} \xi_{(c,b_U)}
 +  B_{(c,b_S)}^{\, T}(\Delta_{1/c})_{(c,b_S)}\xi_{(c,b_S)}
  -C_{(\alpha,\gamma)}^{\,
T}(\Theta_c)_{(\alpha,\gamma)}\omega_{(\alpha,\gamma)}\|^2
\\[2mm]
\leq
\overline{\nu}
\|[(\Delta_{1/c})_{(a,b_S)}\xi_{(a,b_S)},
(\Delta_{1/c})_{(a,b_L)}\xi_{(a,b_L)},(\Delta_{1/c})_{(a,c)}\xi_{(a,c)},\\[3mm]
\,\, \quad \quad \quad \quad \quad \quad
 (\Delta_{1/c})_{(c,b_U)} \xi_{(c,b_U)},(\Delta_{1/c})_{(c,b_S)}\xi_{(c,b_S)},(\Theta_c)_{(\alpha,\gamma)}\omega_{(\alpha,\gamma)}]\|^2
\\[3mm]
\leq \overline{\nu}
\displaystyle \left (\max_{i \in a,1\leq j \leq |b_S|}\displaystyle \frac{1-(w_{b_S})_j}{c\lambda_i(F(\overline x))+(1-(w_{b_S})_j)}\right)^2 \|
\xi_{(a,b_S)}\|^2
\end{array}
\]
\begin{equation}\label{ineq2-added}
\begin{array}{l}
\, \quad +\overline{\nu}
\displaystyle \left (\max_{i \in a}\displaystyle \frac{2}{c\lambda_i(F(\overline x))+2}\right)^2 \|
\xi_{(a,b_L)}\|^2
\\[3mm]
\, \quad +\overline{\nu}
\displaystyle \left (\max_{i \in a,j\in c}\displaystyle \frac{2}{c[\lambda_i(F(\overline x))-\lambda_i(F(\overline x))]+2}\right)^2 \|
\xi_{(a,c)}\|^2
\\[3mm]
\, \quad +\overline{\nu}
\displaystyle \left (\max_{i \in c}\displaystyle \frac{2}{-c\lambda_i(F(\overline x))+2}\right)^2 \|
\xi_{(c,b_U)}\|^2
\\[3mm]
\, \quad +\overline{\nu}
\displaystyle \left (\max_{i \in c, 1\leq j\leq |b_S|}\displaystyle \frac{(w_{b_S})_j+1}{((w_{b_S})_j+1)-c\lambda_i(F(\overline x))}\right)^2 \|
\xi_{(c,b_S)}\|^2
\\[3mm]
\, \quad +
\overline{\nu}
\displaystyle \left (\max_{i \in \alpha, j \in
\gamma}\lambda_i/(\lambda_i+c|\lambda_j|)\right)^2 \|
\omega_{(\alpha,\gamma)}\|^2
\\[3mm]
\leq
\overline{\nu}
\overline{\nu}_{a,b_S}^2 (\overline{\nu}_{a,b_S}+c)^{-2} \|
\xi_{(a,b_S)}\|^2
+\overline{\nu}
\overline{\nu}_{a,b_L}^2 (\overline{\nu}_{a,b_L}+c)^{-2} \|
\xi_{(a,b_L)}\|^2\\[2mm]
+\overline{\nu}
\overline{\nu}_{a,c}^2 (\overline{\nu}_{a,c}+c)^{-2} \|
\xi_{(a,c)}\|^2
+\overline{\nu}
\overline{\nu}_{c,b_U}^2 (\overline{\nu}_{c,b_U}+c)^{-2} \|
\xi_{(c,b_U)}\|^2\\[2mm]
+\overline{\nu}
\overline{\nu}_{c,b_S}^2 (\overline{\nu}_{c,b_S}+c)^{-2} \|
\xi_{(c,b_S)}\|^2
+\overline{\nu}
\overline{\nu}_{\alpha,\gamma}^2 (\overline{\nu}_{\alpha,\gamma}+c)^{-2} \|
\omega_{(\alpha,\gamma)}\|^2\\[2mm]
\leq \overline{\nu}
\overline{\nu}_0^2 c^{-2}
  \|(\xi_{(a,b_S)},\xi_{(a,b_L)},
\xi_{(a,c)},\xi_{(c,b_U)},\xi_{(c,b_S)},
\omega_{(\alpha,\gamma)})\|^2.
\end{array}
\end{equation}
By using  (\ref{anup-1-add0}) and (\ref{ineq2-added}) we have
\begin{equation}\label{eq:inequality003}
\begin{array}{l}
\displaystyle \Big|\Big\langle A^TD_c \Delta d_0, {\cal
A}_c(\overline{y},W_1,W_2)^{-1}[B_{(a,b_S)}^{\, T}(\Delta_{1/c})_{(a,b_S)}\xi_{(a,b_S)}\\[3mm]
\quad \, +B_{(a,b_L)}^{\, T}(\Delta_{1/c})_{(a,b_L)}\xi_{(a,b_L)}
 + B_{(a,c)}^{\, T}(\Delta_{1/c})_{(a,c)}\xi_{(a,c)}
 + B_{(c,b_U)}^{\, T}(\Delta_{1/c})_{(c,b_U)} \xi_{(c,b_U)}\\[3mm]
 \quad \quad \quad \quad \quad \quad +  B_{(c,b_S)}^{\, T}(\Delta_{1/c})_{(c,b_S)}\xi_{(c,b_S)}
  -C_{(\alpha,\gamma)}^{\,
T}(\Theta_c)_{(\alpha,\gamma)}\omega_{(\alpha,\gamma)}]
\Big\rangle \Big|
\\[3mm]
 \leq \| {\cal A}_c(\overline{y},W_1,W_2)^{-1}A^TD_c\Delta
d_0\|\, \|[B_{(a,b_S)}^{\, T}(\Delta_{1/c})_{(a,b_S)}\xi_{(a,b_S)}\\[3mm]
\quad \, +B_{(a,b_L)}^{\, T}(\Delta_{1/c})_{(a,b_L)}\xi_{(a,b_L)}
  + B_{(a,c)}^{\, T}(\Delta_{1/c})_{(a,c)}\xi_{(a,c)}\\[3mm]
 \quad \,+ B_{(c,b_U)}^{\, T}(\Delta_{1/c})_{(c,b_U)} \xi_{(c,b_U)}
  +  B_{(c,b_S)}^{\, T}(\Delta_{1/c})_{(c,b_S)}\xi_{(c,b_S)}
  -C_{(\alpha,\gamma)}^{\,
T}(\Theta_c)_{(\alpha,\gamma)}\omega_{(\alpha,\gamma)}]
\|
\\[3mm]
\leq  \displaystyle \frac{\varrho_0}{\sqrt{2}}c^{-1}
(\|(\Delta\mu, \xi_{(b_U,b_U)},
 \xi_{(b_S,b_S)},
 \xi_{(b_L,b_L)}, { \omega}_{(\alpha,\alpha)}, {
\omega}_{(\beta,\beta)})\|^2\\[3mm]
\quad \quad \,  + 2\|( \xi_{(b_U,b_S)}, \xi_{(b_S,b_L)},
 \xi_{(b_U,b_L)}, \omega_{(\alpha,\beta)})\|^2
)^{1/2}\\[3mm]
\quad \quad \, \times\left( \overline{\nu}_0 \sqrt{\overline{\nu}} c^{-1}  \|(\xi_{(a,b_S)},\xi_{(a,b_L)},
\xi_{(a,c)},\xi_{(c,b_L)},\xi_{(c,b_S)},
\omega_{(\alpha,\gamma)})\| \right)\,
\\[3mm]
\leq   \displaystyle
\frac{\varrho_0\overline{\nu}_0\sqrt{\overline{\nu}}}{4}c^{-2}
(\|(\Delta\mu, \xi_{(b_U,b_U)},
 \xi_{(b_S,b_S)},
 \xi_{(b_L,b_L)}, { \omega}_{(\alpha,\alpha)}, {
\omega}_{(\beta,\beta)})\|^2\\[3mm]
\quad \quad \, \quad \, \quad \, \quad \,\quad \,+ 2\|( \xi_{(b_U,b_S)}, \xi_{(b_S,b_L)},
 \xi_{(b_U,b_L)}, \omega_{(\alpha,\beta)})\|^2\\[3mm]
 \, \quad \quad \,\quad \,\quad \,\quad \, \quad \,+2
 \|(\xi_{(a,b_S)},\xi_{(a,b_L)},
\xi_{(a,c)},\xi_{(c,b_U)},\xi_{(c,b_S)},
\omega_{(\alpha,\gamma)})\|^2)\, .
\end{array}
\end{equation}

By direct calculations we have
\[
\begin{array}{l}
\|\Delta Y\|^2-c^{-1}\left \langle W_1\Delta Y, \Delta Y \right \rangle\\[2mm]
=(\|\xi_{(a,a)}\|^2+2\|\xi_{(a,b_U)}\|^2+2\|\xi_{(c,b_L)}\|^2+\|\xi_{(c,c)}\|^2)\\[2mm]
+2(\|\xi_{(a,b_S)}\|^2-\langle \xi_{(a,b_S)}, (\Delta_{1/c})_{(a,b_S)}\xi_{(a,b_S)} \rangle)\end{array}
\]
\begin{equation}\label{ht2-1}
\begin{array}{l}
+2(\|\xi_{(a,b_L)}\|^2-\langle \xi_{(a,b_L)}, (\Delta_{1/c})_{(a,b_L)}\xi_{(a,b_L)} \rangle)\\[2mm]
+2(\|\xi_{(a,c)}\|^2-\langle \xi_{(a,c)}, (\Delta_{1/c})_{(a,c)}\xi_{(a,c)} \rangle)\\[2mm]
+2(\|\xi_{(c,b_U)}\|^2-\langle \xi_{(c,b_U)}, (\Delta_{1/c})_{(c,b_U)}\xi_{(c,b_U)} \rangle)\\[2mm]
+2(\|\xi_{(c,b_S)}\|^2-\langle \xi_{(c,b_S)}, (\Delta_{1/c})_{(c,b_S)}\xi_{(c,b_S)} \rangle)\\[2mm]
+(\|\xi_{(b_U,b_U)}\|^2-\langle \xi_{(b_U,b_U)}, (\Delta_{1/c})_{(b_U,b_U)}\xi_{(b_U,b_U)} \rangle)\\[2mm]
+(\|\xi_{(b_L,b_L)}\|^2-\langle \xi_{(b_L,b_L)}, (\Delta_{1/c})_{(b_L,b_L)}\xi_{(b_L,b_L)} \rangle)
\end{array}
\end{equation}
and
\begin{equation}\label{ht2-2}
\begin{array}{l}
 \|\Delta\Gamma \|^2-\langle \Delta\Gamma , W_2(\Delta\Gamma
)\rangle  \\[2mm]
=  \left(\|\omega_{(\gamma,\gamma)}\|^2
+2\|\omega_{(\beta,\gamma)}\|^2\right)
+2\left(\|\omega_{(\alpha,\gamma)}\|^2 -\langle
\omega_{(\alpha,\gamma)},(\Theta_c)_{(\alpha,\gamma)}\omega_{(\alpha,\gamma)}\rangle\right)
\\[2mm]
\quad \,  + \left(\|\omega_{(\beta,\beta)}\|^2 -\langle
\omega_{(\beta,\beta)},(\Theta_c)_{(\beta,\beta)}\omega_{(\beta,\beta)}\rangle
\right)\, .
\end{array}
\end{equation}
Now we are ready to estimate  the lower and upper bounds of
$-\langle V (\Delta y), \Delta y \rangle$.
 In light of (\ref{Ahh}), (\ref{anup-opposite}), (\ref{ineq1}),
 (\ref{eq:inequality003}),
  (\ref{ht2-1}) and (\ref{ht2-2}), we have
\begin{equation}\label{lb}
\begin{array}{lcl}
-\langle V(\Delta y), \Delta y \rangle  &\geq&
 {c}^{-1}(\|\xi_{(a,a)}\|^2+2\|\xi_{(a,b_U)}\|^2+2\|\xi_{(c,b_L)}\|^2+\|\xi_{(c,c)}\|^2)\\[2mm]
&& +{c}^{-1}\left(\|\omega_{(\gamma,\gamma)}\|^2+2
\|\omega_{(\beta,\gamma)}\|^2\right)
\\[2mm]
& &  +\underline\kappa_1(c)\|(\Delta\mu,\xi_{(b_S,b_S)},\omega_{(\alpha,\alpha)})\|^2
+\underline\kappa_2(c)\|\xi_{(a,b_S)}\|^2\\[2mm]
&&+\underline \kappa_3(c)\|\xi_{(a,b_L)}\|^2 +\underline \kappa_4(c)\|\xi_{(a,c)}\|^2+\underline\kappa_5(c)\|\xi_{(b_U,b_U)}\|^2\\[2mm]
&& +\underline \kappa_6(c)\|\xi_{(b_U,b_S)}\|^2+\underline\kappa_7(c)\|\xi_{(b_U,b_L)}\|^2+\underline \kappa_8(c)\|\xi_{(b_U,c)}\|^2\\[2mm]
&& +\underline \kappa_9(c)\|\xi_{(b_S,b_L)}\|^2 +\underline \kappa_{10}(c)\|\xi_{(b_S,c)}\|^2+\underline\kappa_{11}(c)\|\xi_{(b_L,b_L)}\|^2\\[2mm]
& &
+\underline\kappa_{12}(c)\|\omega_{(\alpha,\beta)}\|^2+\underline\kappa_{13}(c)\|\omega_{(\alpha,\gamma)}\|^2+
\underline\kappa_{14}(c)\|\omega_{(\beta,\beta)}\|^2,
\end{array}
\end{equation}
where
$$
\begin{array}{l}
\underline \kappa_1(c):=\left( \overline{\sigma}\overline{\eta}+(c-c_0)\right )^{-1}
  -{\varrho_0 \overline{\nu}_0 \sqrt{\overline{\nu}}} c^{-2}\\[2mm]
    \underline \kappa_2(c):=\underline{\nu}
\left({\overline{\sigma}}\overline{\eta}+2(c-c_0)\right)^{-1}
\underline{\nu}_0^2 (\underline{\nu}_0+c)^{-2}-{\varrho_0\overline{\nu}_0\sqrt{\overline{\nu}}}c^{-2}
+\displaystyle \frac{\min_{i\in a}\lambda_i(F(\overline x))}{[c \min_{i\in a}\lambda_i(F(\overline x))+2]}\\[2mm]
\underline \kappa_3(c):=\underline{\nu}
\left({\overline{\sigma}}\overline{\eta}+2(c-c_0)\right)^{-1}
\underline{\nu}_0^2 (\underline{\nu}_0+c)^{-2}-{\varrho_0\overline{\nu}_0\sqrt{\overline{\nu}}}c^{-2}
+\displaystyle \frac{\min_{i\in a}\lambda_i(F(\overline x))}{[c \min_{i\in a}\lambda_i(F(\overline x))+2]}\\[2mm]
\underline \kappa_4(c):=\underline{\nu}
\left({\overline{\sigma}}\overline{\eta}+2(c-c_0)\right)^{-1}
\underline{\nu}_0^2 (\underline{\nu}_0+c)^{-2}-{\varrho_0\overline{\nu}_0\sqrt{\overline{\nu}}}c^{-2}\\[2mm]
\quad \quad \quad \quad \quad \quad \quad \quad \quad +\displaystyle \frac{(\min_{i\in a}\lambda_i(F(\overline x))-\max_{j \in c}\lambda_j(F(\overline x)))}{ [c (\min_{i\in a}\lambda_i(F(\overline x))-\max_{j \in c}\lambda_j(F(\overline x)))+2]}\\[4mm]
\underline \kappa_{6}(c):=2\left(
\overline{\sigma}\overline{\eta}/2+(c-c_0)\right )^{-1}
 -{\varrho_0\overline{\nu}_0\sqrt{\overline{\nu}}}c^{-2}
 \\[3mm]
 \underline \kappa_{7}(c):=\underline \kappa_{6}(c)\\[2mm]
 \underline \kappa_{8}(c):=-{\varrho_0\overline{\nu}_0\sqrt{\overline{\nu}}}c^{-2}
+\displaystyle \frac{-\max_{j\in c}\lambda_j(F(\overline x))}{[-c \max_{j\in c}\lambda_j(F(\overline x))+2]}\\[2mm]
\underline \kappa_{9}(c):=\underline \kappa_{6}(c)
 \end{array}
 $$
$$
\begin{array}{l}
\underline \kappa_{10}(c):=\underline{\nu}
\left({\overline{\sigma}}\overline{\eta}+2(c-c_0)\right)^{-1}
\underline{\nu}_0^2 (\underline{\nu}_0+c)^{-2}-{\varrho_0\overline{\nu}_0\sqrt{\overline{\nu}}}c^{-2}
+\displaystyle \frac{-\max_{j\in c}\lambda_j(F(\overline x))}{[-c \max_{j\in c}\lambda_j(F(\overline x))+2]}\\[3mm]
\underline \kappa_{12}(c):=\underline \kappa_{6}(c)
\\[3mm]
\underline \kappa_{13}(c):=2 c^{-1}[1 -\overline{\nu}_0(\overline{\nu}_0+
c)^{-1}] +2\underline{\nu}
\left({\overline{\sigma}}\overline{\eta}+2(c-c_0)\right)^{-1}
{\underline{\nu}}_0^2 (\underline{\nu}_0+c)^{-2}
-{\varrho_0\overline{\nu}_0\sqrt{\overline{\nu}}}c^{-2} \, ,
\end{array}
$$
and
\[
  \underline \kappa_{5}(c)=\underline \kappa_{11}(c)=\underline \kappa_{14}(c): = \displaystyle \min_{t \in [0,1]}
 \psi(t;c,a_c,b_c,c_0)
\]
with $\psi(\cdot; \cdot)$ being defined as (\ref{phicom}) in Lemma
\ref{phiineq} and
\[
a_c:={c}^{-1}-{\varrho_0\overline{\nu}_0\sqrt{\overline{\nu}}}c^{-2}\,
, \quad b_c:=\overline{\sigma}\overline{\eta}\, .
\]
 It follows from (\ref{k4com}) in Lemma
\ref{phiineq} that for $c\ge {\overline{c}}$,
\[
\underline \kappa_{14}(c)=c^{-1}
-{\rho_0\overline{\nu}_0\sqrt{\overline{\nu}}}c^{-2} -
\frac{\overline{\sigma}\overline{\eta}}{c(\sqrt{c}+\sqrt{c_0})^2}
\, .
\]
Thus, there exists a positive number $\underline \epsilon_1$ such that
 for $c\ge {\overline{c}}$ we have
\[
\min \left\{ \frac{1}{2} \min_{i \notin  \{1,5,11,14\}}  \{\underline \kappa_i(c) \}, \min_{i \in  \{1,5,11,14\}}  \{\underline \kappa_i(c) \} \right\}\ge
c^{-1}-\underline \epsilon_1{c^{-2}}\, .
\]
Therefore, from (\ref{lb}) we have
\begin{equation}\label{lbe}
-\langle V(\Delta y), \Delta y \rangle \geq ( c^{-1}-\underline \epsilon_1c^{-2}) \|\Delta y\|^2\, .
\end{equation}
On the other hand, in light of (\ref{Ahh}),
(\ref{anup-1}), (\ref{ineq2}),  (\ref{eq:inequality003}),
  (\ref{ht2-1}) and (\ref{ht2-2}), we have

  \begin{equation}\label{ub}
\begin{array}{lcl}
-\langle V(\Delta y), \Delta y \rangle  &\leq&
 {c}^{-1}(\|\xi_{(a,a)}\|^2+2\|\xi_{(a,b_U)}\|^2+2\|\xi_{(c,b_L)}\|^2+\|\xi_{(c,c)}\|^2)\\[2mm]
&& +{c}^{-1}\left(\|\omega_{(\gamma,\gamma)}\|^2+2
\|\omega_{(\beta,\gamma)}\|^2\right)
\\[2mm]
& &  +\overline\kappa_1(c)\|(\Delta\mu,\xi_{(b_S,b_S)},\omega_{(\alpha,\alpha)})\|^2
+\overline\kappa_2(c)\|\xi_{(a,b_S)}\|^2\\[2mm]
&&+\overline \kappa_3(c)\|\xi_{(a,b_L)}\|^2 +\overline \kappa_4(c)\|\xi_{(a,c)}\|^2+\overline\kappa_5(c)\|\xi_{(b_U,b_U)}\|^2\\[2mm]
&& +\overline \kappa_6(c)\|\xi_{(b_U,b_S)}\|^2+\overline\kappa_7(c)\|\xi_{(b_U,b_L)}\|^2+\overline \kappa_8(c)\|\xi_{(b_U,c)}\|^2\\[2mm]
&& +\overline \kappa_9(c)\|\xi_{(b_S,b_L)}\|^2 +\overline \kappa_{10}(c)\|\xi_{(b_S,c)}\|^2+\overline\kappa_{11}(c)\|\xi_{(b_L,b_L)}\|^2\\[2mm]
& &
+\overline\kappa_{12}(c)\|\omega_{(\alpha,\beta)}\|^2+\overline\kappa_{13}(c)\|\omega_{(\alpha,\gamma)}\|^2+
\overline\kappa_{14}(c)\|\omega_{(\beta,\beta)}\|^2,
\end{array}
\end{equation}
where
$$
\begin{array}{l}
\overline \kappa_1(c):=\left(
\underline{\sigma}\underline{\eta}/2+(c-c_0)\right )^{-1}
  +{\varrho_0\overline{\nu}_0\sqrt{\overline{\nu}}}c^{-2}\\
\overline \kappa_2(c):=2c^{-1}+\overline{\nu}\underline{\sigma}^{-1}\underline{\eta}^{-1}
\overline{\nu}_0^2 (\overline{\nu}_0+c)^{-2}+{\varrho_0\overline{\nu}_0\sqrt{\overline{\nu}}}c^{-2}\\
\overline \kappa_3(c)=\overline \kappa_4(c):=\overline \kappa_2(c)\\
\overline \kappa_{6}(c):=2\left(
\underline{\sigma}\underline{\eta}/2+(c-c_0)\right )^{-1}
 +{\varrho_0\overline{\nu}_0\sqrt{\overline{\nu}}}c^{-2}
\\
\overline \kappa_{7}(c):=\overline \kappa_{6}(c)\\
 \overline \kappa_{8}(c):=2c^{-1}+{\varrho_0\overline{\nu}_0\sqrt{\overline{\nu}}}c^{-2}\\
\overline \kappa_{9}(c):=\overline \kappa_{6}(c)
\\
\overline \kappa_{10}(c):=\overline{\nu}\underline{\sigma}^{-1}\underline{\eta}^{-1}
\overline{\nu}_0^2 (\overline{\nu}_0+c)^{-2}+{\varrho_0\overline{\nu}_0\sqrt{\overline{\nu}}}c^{-2}\\
\overline \kappa_{12}(c):=\overline \kappa_{6}(c)
\\
\overline \kappa_{13}(c):=\overline \kappa_{2}(c)
\end{array}
$$
and
\[
  \overline \kappa_{5}(c)=\overline \kappa_{11}(c)=\overline \kappa_{14}(c): = \displaystyle \max_{t \in [0,1]} \psi(t;c,a'_c,b'_c,c_0)
\]
with
\[
a^\prime_c:={c}^{-1}+{\varrho_0\overline{\nu}_0\sqrt{\overline{\nu}}}c^{-2}
 \, , \quad b^\prime_c:=\underline{\sigma}\underline{\eta}/2\, .
\]
It follows from (\ref{u4com}) in Lemma \ref{phiineq} that for
$c\ge {\overline{c}}$,
\begin{equation}\label{u4c}
\begin{array}[b]{lcl}
\overline \kappa_{14}(c)&=&\max\{\psi(0;c,a^\prime_c,b^\prime_c,c_0),\psi(1;c,a^\prime_c,b^\prime_c,c_0)\}
\\[2mm]
&=& {\varrho_0\overline{\nu}_0\sqrt{\overline{\nu}}}c^{-2}
+\max\{{c}^{-1}, \left(
{\underline{\sigma}}\underline{\eta}/2+(c-c_0)\right)^{-1}\}\, .
\end{array}
\end{equation}
Thus, there exists a positive number $\mu_0\geq \max \{\varrho_0,\underline \epsilon_1\}$ such that
 for $c\ge {\overline{c}}$ we have
\[
\max \left\{ \frac{1}{2} \max_{i \notin  \{1,5,11,14\}}  \{\overline \kappa_i(c) \}, \min_{i \in  \{1,5,11,14\}}  \{\overline \kappa_i(c) \} \right\}\le
c^{-1}+\mu_0 {c^{-2}}\, .
\]
Therefore, from (\ref{ub}) we have
\begin{equation}\label{ube}
-\langle V(\Delta y), \Delta y \rangle \leq ( c^{-1}+\mu_0c^{-2}) \|\Delta y\|^2\, .
\end{equation}
By (\ref{lbe}) and  (\ref{ube}), noting that $\mu_0\ge \underline \epsilon_1$,  we obtain that
\[
\mu_0 c^{-2}\|\Delta y\|^2 \ge \langle V(\Delta y)+c^{-1} \Delta
y, \Delta y\rangle \ge -\mu_0 c^{-2}\|\Delta y\|^2 \, .
\]
This shows that (\ref{import11}) holds.  The proof is completed.
\qed

\vskip 7 true pt

Now we are ready to state our main  result on the rate of
convergence of  the augmented Lagrangian method for nonlinear
semidefinite nuclear norm composite optimization.

\begin{theorem}\label{thconv1nlsdp}
Suppose that Assumptions (sdnop-A1) and (sdnop-A2) are satisfied.
Let ${c_0}$ and $\underline {\eta}$ be two positive numbers
obtained by {\rm Proposition \ref{pronlsdp}}.
Let  $\overline{\eta}$, $\overline{c}$,  and $\varrho_0$ be
defined as in {\rm (\ref{mu1nlsdp})},  {\rm (\ref{eq:barc-sdp})}, and
{\rm (\ref{varrho0-sdp})}, respectively. Let $\mu_0$ be obtained
by {\rm Proposition \ref{thdcnlsdp}}. Define
\[
\varrho_1:=2\varrho_0 \quad {\rm and} \quad \varrho_2:=4\mu_0.
\]
 Then for any  $c \geq {\overline{c}}$,  there exist two
positive numbers $\varepsilon$ and  $\delta $ $($both depending on
$c$$)$ such that for any $(Y,\mu, \Gamma)
\in\mathbb{B}_{\delta}(\overline{Y},\overline \mu, \overline{\Gamma})$, the
problem
 \[
  \min \ L_c(x, Y,\mu, \Gamma)
\quad {\rm s.t.} \ x\in {\mathbb {B}}_{\varepsilon}(\overline{x})\,
\]
has a unique solution denoted $x_c(Y,\mu, \Gamma)$. The function
$x_c(\cdot, \cdot, \cdot)$ is locally Lipschitz continuous on
$\mathbb{B}_{\delta}(\overline{Y},\overline \mu, \overline{\Gamma})$ and is
semismooth at any point in $\mathbb{B}_{\delta}(\overline{Y},\overline \mu, \overline{\Gamma})$,
 and for any $(\zeta, \Xi) \in
\mathbb{B}_{\delta}(\overline{Y},\overline \mu, \overline{\Gamma})$, we have
\[
\|x_c(Y,\mu, \Gamma)-\overline{x}\| \leq \varrho_1 \|(Y,\mu,
\Gamma)-(\overline{Y},\overline \mu, \overline{\Gamma})\|/{{c}}
\]
and
\[
\|(Y_c(Y,\mu, \Gamma), {\mu}_c(Y,\mu, \Gamma),\Gamma_c(Y,\mu, \Gamma))- (\overline{Y},\overline \mu, \overline{\Gamma})\| \leq \varrho_2 \|(Y,\mu,
\Gamma)-(\overline{Y},\overline \mu, \overline{\Gamma})\|/{{c}},
\]
where $Y_c(Y,\mu, \Gamma)$, ${\mu}_c(Y,\mu, \Gamma)$ and $ {\Gamma}_c(Y,\mu, \Gamma)$ are
defined as
\[
\begin{array}{rl}
Y_c(Y,\mu,\Gamma)&:={\rm D}\theta_c(F(x_c(Y,\mu,\Gamma))+Y/c)^*,\\[3mm]
\mu_c(Y,\mu,\Gamma)&:= \mu+ch(x_c(Y,\mu,\Gamma))\quad \rm{and} \\[2mm]
\quad  \Gamma_c(Y,\mu,\Gamma)&:=\Pi_{{\cal S}^p_+}(\Gamma-c g(x_c(Y,\mu,\Gamma)))\, .
\end{array}
\]
\end{theorem}

\noindent {\bf Proof.} If Assumptions (sdnop-A1) and (sdnop-A2)
are satisfied, then from Propositions \ref{pronlsdp} and
\ref{thdcnlsdp} we know that both  Assumption B1 and  Assumption
B2 (with  $\gamma =2$) made in Section \ref{general-discussions} are
satisfied. Then the conclusions in this theorem  follow from
Theorem \ref{comthconv1}.
 \qed

\section{Conclusions}\label{final-section}
\setcounter{equation}{0}
 This  paper provides an analysis on   the
rate of convergence of the augmented Lagrangian method for solving the
nonlinear semidefinite nuclear norm optimization problem. By assuming that $K$ is a closed convex
 cone, and that ${\rm D}\theta_c(\cdot)$ and $\Pi_{K^*}(\cdot)$ are semismooth everywhere,
we first establish a general result on the  rate of convergence of
the augmented Lagrangian method for a class of  general
composite optimization problems.
Then we apply this general result to the nonlinear semidefinite
nuclear norm optimization problem under the constraint nondegeneracy condition and
the strong second order sufficient condition. The methodology  suggests us that
we may verify Assumptions B1 and B2 to obtain the rate of convergence of the augmented Lagrange method for other optimization problems.


\end{document}